\documentclass[11pt]{article}
\usepackage[top = 1in, bottom = 1in, 
left = 1in, right = 1in]{geometry}

\usepackage{amsmath}
\usepackage{amsfonts}
\usepackage{amssymb}
\usepackage{amsrefs}
\usepackage{mathtools}

 \usepackage[disable]{todonotes} 

\usepackage{cases}
\usepackage{xcolor}

\newcommand{\dt}{\partial_t}
\newcommand{\dv}{\mathrm{div}}
\newcommand{\dr}{\partial_r}
\newcommand{\initial}{\text{initial data}}
\newcommand{\subeqref}[2]{$\eqref{#1}_{#2}$}
\newcommand{\abs}[2]{\bigl| #1 \bigr|^{#2}}
\newcommand{\norm}[2]{\bigl\Arrowvert #1 \bigr\Arrowvert_{#2}}

\newcommand{\stnorm}[2]{L_t^{#1}L_x^{#2}}

\newtheorem{thm}{Theorem}[section]
\newtheorem{lm}{Lemma}

\newenvironment{pf}{\paragraph{Proof}}{\hfill$\square$}

\newenvironment{rmk}{\paragraph{Remark}}{}

\numberwithin{equation}{section}
\allowdisplaybreaks
\title{On the Expanding Configurations of Viscous Radiation Gaseous Stars: Isentropic and Thermodynamic Models}

\author{Xin Liu}

\begin{document}
\allowdisplaybreaks
\maketitle

\abstract{In this work, we study the stability of the expanding configurations for
radiation gaseous stars. Such expanding configurations exist for the isentropic
model and the thermodynamic model. In particular, we focus on the effect of the
monatomic gas viscosity tensor on the expanding configurations. With respect to
small perturbation, it is shown the self-similarly expanding homogeneous
solutions of the isentropic model is unstable, while the linearly expanding
homogeneous solutions of the isentropic model and the thermodynamic model are
stable. This is an extensive study of the result in arXiv:1605.08083 by Hadzic
and Jang.}
\\

\noindent{{\bf MSC class:} 35A01, 35Q30, 35Q35, 35Q85, 76E20}

\tableofcontents

\section{Introduction}

The evolution of a viscous gaseous star can be described by the system 
\begin{equation}\label{EGS}
	\begin{cases}
		\dt \rho + \dv (\rho \vec u) = 0 & \text{in} ~ \Omega(t), \\
		\dt (\rho \vec u) + \dv (\rho \vec u \otimes \vec u) + \nabla p + \rho\nabla\psi = \dv \mathbb S & \text{in} ~ \Omega(t), \\
		\Delta \psi = \rho & \text{in} ~ \mathbb R^3,
	\end{cases}
\end{equation}
where $ \rho, \vec u, \psi, \Omega(t) $ denote the density, the velocity field, the self-gravitation potential and the evolving occupied domain. $ \mathbb S = \mu (\nabla \vec{u} + \nabla \vec{u}^\top ) + \lambda \dv\vec u \, \mathbb I_3 $ is the Newtonian viscous tensor. The viscosity coefficients are taken as the ones for the monatomic gas in this work. That is,
\begin{equation}\label{weakviscosity}
	\mu = \text{constant} > 0, 2\mu + 3 \lambda = 0.
\end{equation}
This is indicated by the kinetic theory (see \cite[pp 3, (1,11)]{Lions1996}). $  p $ is the pressure potential. In this work, we will consider two cases: the isentropic case
\begin{equation}\label{eqstt:isentropic}
	p = \rho^{4/3},
\end{equation}
and the non-isentropic(thermodynamic) case
\begin{equation}\label{eqstt:non-isentropic}
	p = K\rho \theta ~~ ( K = \text{constant} > 0),
\end{equation}
where $ \theta $ is the temperature distribution. 

The system \eqref{EGS} with \eqref{eqstt:isentropic} is completed with the following kinetic boundary conditions
\begin{equation*}
	\begin{aligned}
		 \vec u \cdot \vec{n} = \mathcal V ~~~~~~~~ & \text{on} ~ \Gamma(t) : = \partial \Omega(t), \\
		 p\mathbb I_3 \vec{n} - \mathbb S \vec{n} = 0 ~~~~~~~~ & \text{on} ~ \Gamma(t), 
	\end{aligned}
\end{equation*}
where $ \mathcal V $ is the normal velocity of the moving surface $ \Gamma(t) $. Such a system is referred to as {\bf the isentropic  model for radiation gaseous stars} in this work.

On the other hand, when the equation of state is taken as \eqref{eqstt:non-isentropic}, the following equation for the temperature is integrated into the system,
\begin{equation}\label{Theta-Eq}
	\dt (c_\nu \rho\theta) + \dv( c_\nu \rho\theta \vec u) + p_\theta \dv \vec{u} - \Delta \theta = \epsilon \rho + \mathbb S: \nabla\vec{u} ~~~~~~~~ \text{in} ~ \Omega(t),
\end{equation}
where $ c_\nu > 0, \epsilon > 0 $ denote constants for the specific heat coefficient and the rate of generation of energy. The full system consisted of \eqref{EGS} and \eqref{Theta-Eq} is imposed with the boundary conditions, 
\begin{equation*}
	\begin{aligned}
		\vec u \cdot \vec{n} = \mathcal V ~~~~~~~~ & \text{on} ~ \Gamma(t) : = \partial \Omega(t), \\
		 p\mathbb I_3 \vec{n} - \mathbb S \vec{n} = 0 ~~~~~~~~ & \text{on} ~ \Gamma(t), \\
		 \theta = 0 ~~~~~~~~ & \text{on} ~ \Gamma(t).
	\end{aligned}
\end{equation*}
Such a system is referred to as {\bf the thermodynamic model for radiation gaseous stars } in this work. In fact, we will show the existence of expanding configurations when $ c_\nu = 3K $ in the following and the stability of such expanding configurations is investigated. 

The gaseous star problem has been studied in a huge number of literatures. To name a few, in \cite{AuchmutyBeals1971}, Auchmuty and Beals studied the variational solutions for the rotating gaseous stars. That is to find the solutions (equilibrium states) to the following system
\begin{equation}\label{RotatingGS}
	\begin{cases}
		\dv(\rho \vec{u} ) = 0 & \text{in} ~ \Omega : = \lbrace \rho > 0 \rbrace \subset \mathbb R^3 , \\
		\dv(\rho \vec{u} \otimes \vec u ) + \nabla p(\rho) = - \rho \nabla\psi & \text{in} ~ \Omega, \\
		\Delta \psi = \rho & \text{in} ~ \mathbb R^3,
	\end{cases}
\end{equation}
where $ p(\rho) $ is a given function of the density $ \rho $. Under some conditions on the angular momentum or the angular velocity, and the equation of state, it is shown that there exists at least one compactly supported solution (with $ \text{supp} \, \rho \subset B_R $ for some $ 0 < R < + \infty $). In \cite{Caffarelli1980}, Caffarelli and Friedman showed that such solutions have at most a finite number of rings (i.e. the support of the density distribution is consisted of a finite number of connected components). Also, the regularity of the solutions was studied. See also in \cite{Friedman1981,Li1991,Chanillo1994,Luo2004,Luo2008a,Luo2009b,Chanillo2012,Luo2014,Wu2013,Deng2002} for more discussions on the related problems. 

On the other hand, when $ \vec u \equiv 0 $ in \eqref{RotatingGS} (referred to as the non-rotating gaseous stars), the existence of equilibrium was established by Chandrasekhar in \cite{Chandrasekhar1958}. In particular, Lieb and Yau showed that the solution has to be spherically symmetric in \cite{Lieb1987}. In the case when $ p = \rho^\gamma $, such solutions are called the Lane-Emden solutions for the non-rotating gaseous stars. 

The stability and instability theory of these solutions was studied in \cite{Lin1997}, where Lin studied the eigenvalue of the linearization of the associated evolutionary problem. A conditional nonlinear stability theory was given by Rein in \cite{Rein2003a}.

The compactly supported solutions for the non-rotating gaseous star admit the physical vacuum boundary. That is, the sound speed is only $1/2$-H\"older continuous across the gas-vacuum interface. Such singularity makes the corresponding evolutionary problem challenging (see \cite{Liu1996}). It is only recently that the local well-posedness in smooth functional spaces is studied by Coutand, Lindblad and Shkoller \cite{Coutand2010,Coutand2011a,Coutand2012}, Jang and Masmoudi \cite{Jang2009b,Jang2015}, Luo, Xin and Zeng \cite{LuoXinZeng2014} in the setting of one spatial dimension, three spatial dimension and spherical symmetry with or without self-gravitation. See \cite{Jang2010} for a viscous flow. We refer to \cite{ZengHH2015a,ZengHH2015,ZengHH2017,Liu1996} and the references therein for other discussions on the vacuum-interface problems.

With the local well-posedness theory, it is able to analyse the nonlinear stability and instability of the Lane-Emden solutions for the non-rotating gaseous stars. In particular, the works from Jang and Tice \cite{Jang2008a,Jang2014,Jang2013a} show that for $ 6/5 \leq \gamma < 4/3 $, the Lane-Emden solutions are unstable, and additional viscosities can not reduce such instability. When it comes to the case $ 4/3 < \gamma < 2 $, the asymptotic stability theory is first studied in \cite{LuoXinZeng2016} in the viscous case. See \cite{LuoXinZeng2015} for the model with a degenerate viscosity.

Meanwhile, the isentropic model for radiation gaseous stars (that is, the gaseous star problem with \eqref{eqstt:isentropic}. See \cite{Chandrasekhar1958,Hadzic2016,Deng2002}.) is rather complicated. On one hand, only for a specific total mass $ \int \rho\, dx =  M_* > 0 $, called the critical mass, the non-rotating gaseous star problem admits solutions with compact support. On the other hand, by studying a class of self-similar solutions in \cite{Fu1998}, the authors showed the existence of expanding and collapsing solutions with variant total mass. Such phenomena are far from fully understood. Recently, in \cite{Hadzic2016}, the authors show the asymptotic stability of the expanding solutions. The damping structure in the Lagrangian coordinates and the conservation of energy for smooth solutions are the two important ingredients in their work. When considering the heat generation and heat conductivity in the gaseous star, i.e. the thermodynamic model for radiation gaseous stars, in \cite{XL2017} the author shows that for $ 1/6 < \epsilon K < 1 $, there exist infinitely many self-similar equilibrium states for the non-rotating radiation gaseous star problem. Also, there exist self-similar expanding and collapsing configurations just as in the isentropic model. 

\subsection{The expanding configurations with spherically symmetric motions}\label{section:introduction_of_expanding_sol}

We will study the aforementioned expanding solutions of the isentropic model and the thermodynamic model for the radiation gaseous stars. Considering the spherically symmetric motion, by denoting, $ \vec u(\vec x,t) = u(r,t) \cdot \dfrac{\vec x}{r}, \rho(\vec x,t) = \rho(r,t), \theta(\vec x, t) = \theta(r,t), p(\vec x, t) = p(r,t), \Omega(t) = B_{R(t)} $ where $ r = \abs{\vec x}{} \in [0,R(t) ), \vec{x} \in \mathbb R^3 $, the system \eqref{EGS} can be written as, 
\begin{equation}\label{evol-RGS}
	\begin{cases}
		\dt(r^2\rho) + \dr(r^2 \rho u) = 0 ,\\
		\dt(r^2 \rho u) + \dr(r^2 \rho u^2) + r^2 \dr p + \rho \int_0^r s^2 \rho(t,s)\,ds \\
		~~~~~~~~~~  = ( 2\mu +\lambda)  r^2 \bigl\lbrack \dfrac{\dr(r^2 u)}{r^2} 		\bigr\rbrack_r = \dfrac{4\mu}{3} r^2 \bigl\lbrack \dfrac{\dr(r^2 u)}{r^2} 		\bigr\rbrack_r ,\\
		\dt(c_\nu r^2 \rho \theta) + \dr(c_\nu r^2 \rho \theta u) + K\rho\theta \dr( r^2 u) - \dr(r^2 \dr \theta) - \varepsilon r^2 \rho  \\
		~~~~~~~~~~ = 2 \mu r^2 \bigl( (\dr u)^2 + 2\bigl(\dfrac{u}{r}\bigr)^2 \bigr) + \lambda r^2 \bigl( \dr u + 2 \dfrac{u}{r}\bigr)^2\\
		~~~~~~~~~~ 
		= \dfrac{4\mu}{3} r^2 \bigl( \dr u - \dfrac{u}{r} \bigr)^2 ~~ \text{(thermodynamic model)},
	\end{cases}
\end{equation}
where
\begin{equation*}
	p = \rho^{4/3} ~ \text{(isentropic model)} ~~ \text{or} ~~ K\rho\theta ~ \text{(thermodynamic model)}.
\end{equation*}
The boundary conditions then can be written as
\begin{equation}\label{boundarycondition}
	\begin{gathered}
	u(R(t),t) = \dt R(t), ~ u(0,t) = 0, \\
	p - \mathfrak B \bigr|_{r = R(t)} = 0,
	\\
	\theta(R(t),t) = 0 ~ \text{(thermodynamic model)},
\end{gathered}
\end{equation}
where $ \mathfrak B = (2\mu + \lambda) \dr u + 2 \lambda \dfrac{u}{r} = \dfrac{4\mu}{3} \bigl( \dr u - \dfrac{u}{r} \bigr) $.

The expanding configurations are obtained by searching solutions to \eqref{evol-RGS} with the following self-similar ansatz 
\begin{equation}\label{HomogeneousSol}
	\begin{aligned}
		r = r_\alpha & = \alpha(t) y, \\
		\rho = \rho_\alpha(t,r) & = \alpha^{-3}(t) \bar\rho (y),\\
		\theta = \theta_\alpha(t,r) & = \alpha^{-1}(t) \bar\theta(y) ~ \text{(thermodynamic model)}, \\
		u = u_\alpha (t,r) & = \alpha'(t) y,
	\end{aligned}
\end{equation}
for $ y \in [0, R_0) $ with some $ 0 < R_0 < +\infty $. 

Notice, when the viscosity coefficients are assumed to be the ones for the monatomic gases \eqref{weakviscosity} as above, the viscous tensor $ \mathbb S $ automatically vanishes and so does $ \mathfrak B $ on the boundary for the self-similar ansatz \eqref{HomogeneousSol}.

\subsubsection*{The isentropic model: self-similarly and linearly homogeneous solutions} 

Considering the isentropic model for the radiation gaseous stars, i.e. $ p = \rho^{4/3} $, after plugging in the self-similar ansatz \eqref{HomogeneousSol}, $ \alpha, \bar\rho $ will satisfy the following equations: for some constant $ \delta $, 
\begin{gather}
		\partial_y^2 \bar\rho^{1/3} + \dfrac{2}{y}\partial_y\bar\rho^{1/3} + \dfrac{1}{4} \bar\rho + \dfrac{3}{4} \delta  = 0, \label{HomogeneousSol02}\\
		\alpha^2(t) \alpha''(t) = \delta, ~ t \geq 0, \label{HomogeneousSol03}
\end{gather}
with the conditions
\begin{equation*}
	\begin{gathered}
		\alpha(0) = a_0 > 0, ~ \alpha'(0) = a_1, \\
		\bar\rho(y) > 0, ~ y \in [0, R_0) ~~ \text{and} ~~ \bar\rho(0) = 1, ~\partial_y \bar\rho(0) = 0, ~ \bar\rho(R_0) = 0.
	\end{gathered}
\end{equation*}
Notice, for a fixed $ \delta $, $ R_0 = R_0(\delta) $ is the first zero of the solution $ \bar\rho $ to the initial value problem \eqref{HomogeneousSol02}; $ \alpha $ is the solution to the initial value problem \eqref{HomogeneousSol03}; the initial data for the initial value problems is given in the conditions above. 

The properties of the solutions to \eqref{HomogeneousSol02}, \eqref{HomogeneousSol03} are well-known. See \cite{Fu1998,Hadzic2016}. We summarize the results in the following.

Considering \eqref{HomogeneousSol02}, for $ \delta \geq \delta^* $ with some $ - \infty < \delta^* < 0 $, there exists a solution $ \bar\rho $ with first zero $ R_0 < +\infty $. $ \bar\rho $ satisfies the following behavior near the origin $ y = 0 $ and the boundary $ y = R_0 $: 
\begin{enumerate}
	\item[(1)] $ \bar\rho $ admits the physical vacuum property near the vacuum boundary, i.e., $$ -\infty < \partial_y \bar\rho^{1/3}(R_0) < -C < 0; $$
	\item[(2)] $ \bar\rho^{1/3} \in C^\infty(0,R_0) $ and is analytic near $ y = 0 $ as well as $ y = R_0 $; 
	\item[(3)] $ \bar\rho^{1/3}(y) = A_1 + A_2 y^2 + O(y^4) $, $ z \sim 0 $ for some constants $ A_1, A_2 $ and $ \partial_y^{2k+1} \bar\rho^{1/3}(0) = 0 $ for any nonnegative integer $ k \geq 0 $. 
\end{enumerate}
Notice, the physical vacuum property of $ \bar \rho $ above means near the boundary $ y = R_0 $, 
\begin{equation}\label{densitydistance}
	\bar\rho^{1/3}(y) \simeq R_0 - y. 
\end{equation}

Considering \eqref{HomogeneousSol03}, we have the following properties of the solutions to the initial value problem:
\begin{enumerate}
	\item[(1)] If $ \delta > 0 $, then $ \alpha(t) > 0 $ for all $ t > 0 $, and $ \lim\limits_{t\rightarrow \infty} \alpha(t) = \infty $. Moreover, there exist $ c_1, c_2 > 0 $ such that
	\begin{equation}\label{LinearlyExpanding} \alpha(t) \sim_{t\rightarrow \infty} c_1(1+c_2t); \end{equation}
	\item[(2)] if $ \delta = 0 $, then $ \alpha(t) = a_0 + a_1 t $. In particular, it admits \eqref{LinearlyExpanding};
	\item[(3)] if $ \delta < 0 $, define
	\begin{equation}\label{first-cosmology-speed} a_1^* = \sqrt{\dfrac{2|\delta|}{a_0}}. \end{equation}
	\begin{enumerate}
		\item[(a)] If $ a_1 = a_1^* > 0  $, then $ \alpha(t) > 0 $ for all $ t > 0 $ and it is explicitly given by the formula: 
		\begin{equation}\label{Self-similarExpanding} \alpha(t) = \bigl( a_0^{3/2}+\dfrac{3}{2} a_{0}^{1/2} a_1 t \bigr)^{2/3}, ~ t \geq 0; \end{equation}
		\item[(b)] if $ a_1 > a_1^* $, then $ \alpha(t) > 0 $ for all $ t > 0 $, $ \lim\limits_{t\rightarrow \infty}\alpha(t) = \infty $, and there exist $ c_1, c_2 > 0 $ such that $ \alpha(t) $ satisfies the asymptotic behavior \eqref{LinearlyExpanding};
		\item[(c)] if $ a_1 < a_1^* $, then there is a time $ 0 < T < \infty $ such that $ \alpha(t) > 0 $ in $ ( 0, T )$ and $ \alpha(t) \rightarrow 0 $ as $ t \rightarrow T^- $. Moreover, there exists a constant $ k_1 > 0 $ such that
		$$ \alpha(t) \sim_{t\rightarrow T^-} k_1(T-t)^{2/3}. $$
	\end{enumerate}
\end{enumerate}
\begin{rmk}
	The critical constant $ a_1^* $ corresponding to the escape velocity of the gaseous star. Also, the solutions given by (3)(a), (3)(b) and (3)(c) are the ones with velocity equal to, larger than and smaller than the escape velocity respectively.
\end{rmk}

In the following, we will refer to the constants $ c_1c_2, a_1 $ as {\bf the expanding rate } for the expanding solutions. Also, we will adopt the terminologies used in \cite{Hadzic2016}. An expanding solution is called a {\bf self-similarly expanding homogeneous solution} if it is given by \eqref{Self-similarExpanding}; it is called a {\bf linearly expanding homogeneous solution} if it admits \eqref{LinearlyExpanding}. 


Before we move to the expanding solutions of the thermodynamic model, let me write down the corresponding energy identity for the system \eqref{evol-RGS} with $ p = \rho^{4/3} $. Temporally, it is assumed $ \rho, u, R(t) $ exist and is regular enough. Moreover, $ \rho(R(t),t) = 0 $. Multiply \subeqref{evol-RGS}{2} with $ u $ and integrate the resulting expression. After integration by parts, it yields 
\begin{equation*}
		\dfrac{d}{dt} E + D = 0, 
\end{equation*}
where
\begin{equation}\label{energy-eulercoordinate}
	\begin{aligned}
		& E = \dfrac{1}{2} \int_0^{R(t)} r^2 \rho u^2 \,dr + 3 \int_0^{R(t)} r^2 \rho^{4/3} \,dr - \int_0^{R(t)} r \rho \int_0^r s^2 \rho \,ds \,dr,\\
		& D = \int_0^{R(t)} 2 \mu (r^2 u_r^2 + 2  u^2 ) + \lambda (r u_r + 2 u)^2 \,dr = \int_0^{R(t)}\dfrac{4}{3}\mu \bigl( r u_r - u \bigr)^2 \,dr \geq 0 . 
	\end{aligned}
\end{equation}
Hence the total energy $ E $ is non-increasing over time.

For the expanding solutions in the inviscid case as well as the case when the viscosity coefficients are assumed to be \eqref{weakviscosity}, the energy for $ (\rho_\alpha, u_\alpha, I_\alpha) $ is given by (\cite[Lemma 2.8]{Hadzic2016})
\begin{equation}\label{energy-inviscid}
	E_{\alpha} = \biggl( a_1^2 +\dfrac{2\delta}{a_0}  \biggr) \int_0^{R_0} s^4 \bar\rho(s) \,ds.
\end{equation}

In the inviscid case, the stability of such expanding solutions is given by Had\v{z}i\'{c} and Jang in \cite{Hadzic2016}. The main observation in \cite{Hadzic2016} is that in the Lagrangian coordinate, the perturbation variable admits a non-linear wave equation with damping. Such a structure provides the dissipation estimate and therefore controls the nonlinearity. Moreover, the conservation of the physical energy $ E $ defined above will gives the energy coercivity together with the spectral gap of the associated linearized operator. An additional assumption for the stable expanding solutions in their work is that $ \delta > - \epsilon $ for some small enough $ \epsilon > 0 $ in order to control some extra forcing terms. 

It should be emphasized that the conservation of energy plays important role in Had\v{z}i\'{c} and Jang's work. Roughly speaking, the spectral gap of the associated linearized operator gives a weighted Poincar\'e type inequality in the space of functions with vanishing weighted mean value. The conservation of energy gives the control the remaining weighted mean value of the perturbation variable. Indeed, the difference of the perturbation variable and its mean value is of order $ \mathcal O(2) $, and therefore can be treated as nonlinearity and controlled by the total energy. 

The goal of current work is to analysis the effect of viscosities on such expanding solutions. To be more precise, notice that when the viscosity coefficients are taken as in \eqref{weakviscosity}, 
the expanding solutions $(\rho_\alpha, u_\alpha, I_\alpha(t))$ are also solutions to \eqref{evol-RGS}. 
Indeed, when the viscosity coefficients satisfy \eqref{weakviscosity}, 
$ D(t) \equiv 0 $ if and only if $ ru_r - u = 0 $, i.e. $ u(r,t) = a(t) r $ for some $ a(t) > 0 $. Then from the continuity equation \subeqref{evol-RGS}{1}, $ \rho(r,t) = e^{-3\int_0^t a(s)\,ds}\bar\rho(e^{-\int_0^ta(s)\,ds} r) $. Then it follows that, after denoting $ \alpha(t) = e^{\int_0^t a(s)\,ds} $, the solutions to \eqref{evol-RGS} with $ D(t) \equiv 0 $ for  $ t > 0 $ are given by the self-similar solutions. This work is devoted to study the stability of such solutions. The main difficulty here is the loss of conservation of energy due to the viscosities. In fact, the loss of energy conservation will lead to the instability of the self-similarly expanding homogeneous solutions, i.e., the one given by \eqref{Self-similarExpanding}. 

To study the stability of the linearly expanding homogeneous solutions, another difficulty is due to the choice of viscosity coefficients \eqref{weakviscosity}. As one will see (and have already observed from the above analysis of $ D(t) $), under such a condition, the dissipation $ D(t) $ loses the strictly positive definiteness in the Lagrangian coordinates (compared with the one in \cite{LuoXinZeng2015,LuoXinZeng2016}). Therefore it seems hard to get the full regularity. 

Our strategy is to notice that the extra damping structure in the Lagrangian coordinates (just as in Had\v{z}i\'{c} and Jang's work) will provide extra regularity information of the solutions, and therefore eventually leads to the regularity analysis. However, the story is not so simple. In fact, the extra integrability of the solutions provided by the damping structure has different temporal weights with the one provided by the viscosity tensor. This will break the pattern of energy estimates as one will see. In particular, we will fail to get the decay estimate of higher order energy. We perform delicate energy estimates with negative temporal weights in this work, which will control the growth of higher order energy. In the end, we are able to make use of the imbalance of the temporal weight in the damping term and the viscosity tensor to recover the control of the perturbation and show the stability of the linear expanding solutions.


\subsubsection*{The thermodynamic model: the linearly homogeneous solutions}

In \cite{XL2017}, the following linearly homogeneous solution is given for \eqref{evol-RGS} for the thermodynamic model. Let $ ( \bar \rho, \bar \theta ) $ be a equilibrium state for \eqref{evol-RGS} with $ 1/6 < \epsilon K < 1 $. That is, it satisfies the following
\begin{equation}\label{HomogenerousSol04}
	\begin{cases}
		y^2 ( K \bar\rho\bar\theta )_y + \bar\rho \int_0^y s^2 \bar\rho(s)\,ds = 0, \\
		- ( y^2 \bar\theta_y )_y = \epsilon y^2 \bar\rho,
	\end{cases}
\end{equation}
with $\bar\theta, \bar\rho > 0 ~ \text{in} ~ [0,R_0) , ~ \text{and} ~ \bar\theta ( R_0 ) = \bar\rho ( R_0 ) = 0 $. Then after plugging the ansatz \eqref{HomogeneousSol} into the system \eqref{evol-RGS}, it holds the following
\begin{equation*}
		\alpha''(t)\alpha(t) = 0, ~ 
		(3K - c_\nu) \alpha^{-2}(t) \alpha'(t) = 0. 
\end{equation*}
Such a system admits non-trivial solutions if and only if $$ 3K - c_\nu = 0. $$ Therefore the linearly expanding solution of the form \eqref{HomogeneousSol} is given when $ 3K = c_\nu $ and it admits 
\begin{equation}\label{Rene:expandingmodule}
	\alpha(t) = a_0 + a_1 t. 
\end{equation}
Notice it is the same as the case when $ \delta = 0 $ in the isentropic model. Similarly, we will refer to the constant $ a_1 $ as {\bf the expanding rate} of the linearly expanding solutions.

Also, near the boundary $ y = R_0 $, $ \bar\rho, \bar \theta $ admits the following properties
\begin{equation*}
	\begin{gathered}
		-\infty < \bigl( \bar\rho^{\frac{\epsilon K}{1-\epsilon K}} \bigr)_y , \bigl( \bar\theta \bigr)_y \leq -C < 0, \\
		\abs{\partial_y^k \bar\rho}{} \leq O(\bar\rho^{\frac{1-(1+k)\epsilon K}{1-\epsilon K}}), ~ \abs{\partial_y^k \bar\theta}{} \leq O(1) + O(\bar\theta^{\frac{1-\epsilon K}{\epsilon K} - k + 2}),  ~ k \in \mathbb Z^+.
	\end{gathered}
\end{equation*}
Notice, this implies that near the boundary $ y = R_0 $,
\begin{equation}\label{Rene:densitydistance}
	\bar\rho^{\frac{\epsilon K}{1-\epsilon K}}(y), ~ \bar\theta(y) \simeq R_0 - y.
\end{equation}

While the difficulty caused by the loss of conservation of energy still exists, unlike the isentropic case, the physical energy for the thermodynamic model is not monotone. In fact, it holds 
\begin{equation*}
	\dfrac{d}{dt} E = R^2(t)\dr \theta(R(t),t) + \epsilon \int_0^{R(t)} r^2 \rho \,dr,
\end{equation*}
where
\begin{equation*}
	\begin{aligned}
		E = \dfrac{1}{2} \int_0^{R(t)} r^2 \rho u^2 \,dr + c_\nu \int_0^{R(t)} r^2 \rho \theta \,dr - \int_0^{R(t)} r \rho \int_0^r s^2 \rho \,ds \,dr.
	\end{aligned}
\end{equation*}
In the Lagrangian coordinates, this means that there are several linear forcing terms with no determined sign, which will cause troubles to close the energy estimates. To overcome these difficulties, while doing the temporal weight estimates, we simply integrate the forcing terms in the temporal variable and it turns out that such calculations will give an extra denominator  proportioned to the expanding rate $ a_1 $. Therefore, by letting $ a_1 $ be large enough, we will be able to close the energy estimates. 


\subsection{Main results}

We will state the main results of this work in this section. 

\todo{formally write down the main results here}
\begin{thm}[The Isentropic Model: Self-similarly Expanding Solutions]\label{theorem1}
The self-similarly expanding solutions of the isentropic model for the radiation gaseous stars given by \eqref{HomogeneousSol} to \eqref{evol-RGS} with 
\eqref{HomogeneousSol02} and \eqref{Self-similarExpanding} is unstable in the following sense. For any small perturbation to the self-similarly expanding solution with a total kinetic and potential energy satisfying any of the conditions given in Section \ref{Section:instability}, there is a time $ 0 < \bar T < \infty  $ such that the perturbation will grow large beyond this time.
\end{thm}

\begin{thm}[The Isentropic Model: Linearly Expanding Solutions]\label{theorem2}
The linearly expanding solutions of the isentropic model for the radiation gaseous stars given by \eqref{HomogeneousSol} to \eqref{evol-RGS} with \eqref{HomogeneousSol02} and \eqref{LinearlyExpanding} is stable with respect to small perturbation if 
\begin{equation}\label{thm1.2-001} \delta > - \dfrac{a_0a_1^2}{8}.\end{equation}
\end{thm}

\begin{thm}[The Thermodynamic Model: Linearly Expanding Solutions]\label{theorem3}
The linearly expanding solutions of the isentropic model for the radiation gaseous stars given by \eqref{HomogeneousSol} to \eqref{evol-RGS} with \eqref{HomogenerousSol04} and \eqref{Rene:expandingmodule} is stable with respect to small perturbation if the expanding rate of the expanding solutions $ a_1 $ is large enough. 
\end{thm}

\section{Lagrangian Formulations}

In this section, we will write down the system \eqref{evol-RGS} in a fixed domain (in the Lagrangian coordinates). Denote the Lagrangian spatial and temporal variables as $ (x,t) $. Similarly as in \cite{LuoXinZeng2016}, the Lagrangian unknown $ r(x,t) $ is defined by
\begin{equation}\label{Lagrangian001}
	\int_{0}^{r(x,t)} s^2 \rho(s,t) \,ds  = \int_{0}^{\alpha(t)x} s^2 \rho_\alpha(s,t) \,ds = \int_0^x s^2 \bar\rho(s)\,ds, ~ \text{and} ~ r(x,0) = x, ~ x \in [0, R_0),
\end{equation}
where $ 0 < R_0 < +\infty $ is the first zero of $ \bar\rho $ (and $ \bar\theta $).
After applying spatial and temporal derivative to \eqref{Lagrangian001}, it follows,
\begin{equation}\label{Lagrangian002}
		\dt r(x,t) = u(r(x,t),t), ~~ 
		\rho(r(x,t),t) = \dfrac{x^2 \bar\rho(x)}{r^2r_x},
\end{equation}
where the continuous equation \subeqref{evol-RGS}{1} is applied. Then the expanding solution \eqref{HomogeneousSol} is given by
\begin{equation*}
	\begin{aligned}
		r_\alpha = \alpha(t) x, ~ \rho_\alpha = \alpha^{-3}(t) \bar\rho(x), ~ \theta_\alpha = \alpha^{-1}(t) \bar\theta(t) \text{(thermodynamic model)}, ~ u_\alpha = \alpha'(t) x,
	\end{aligned}
\end{equation*}
where $ \alpha, \bar\rho, \bar\theta $ are given in Section \ref{section:introduction_of_expanding_sol}. Moreover, in the Lagrangian coordinates, the unknowns $ \lbrace \rho,u,\theta\rbrace $ in the moving domain $ [0,R(t)) $ are now replaced by $ \lbrace r, r_t ,\theta \rbrace $ in the fixed domain $ [0,R_0) $. 

Notice, we use the same notation $ \theta $ for the temperature in the original problem \eqref{evol-RGS} and in the Lagrangian coordinates.
The system \eqref{evol-RGS} can be written in terms of the Lagrangian variables in the Lagrangian coordinates as follows, for $ x \in [0,R_0) $,
\begin{equation}\label{eq:Lagrangian}
	\begin{cases}
		\bigl(\dfrac{x}{r}\bigr)^2 \bar\rho \dt^2 r + P_x + \bigl(\dfrac{x}{r}\bigr)^4  \dfrac{\bar\rho}{x^2} \int_0^x s^2 \bar\rho(s)\,ds = \dfrac{4}{3}\mu \bigl\lbrack  \dfrac{(r^2r_t)_x}{r^2r_x} \bigr\rbrack_x ,\\
		3K x^2 \bar \rho \dt \theta + K\dfrac{x^2\bar\rho\theta}{r^2r_x} (r^2r_t)_x - \bigl\lbrack \dfrac{r^2}{r_x} \theta_x \bigr\rbrack_x - \epsilon x^2 \bar\rho = \dfrac{4}{3}\mu r^2 r_x \bigl( \dfrac{r_{xt}}{r_x}- \dfrac{r_t}{r} \bigr)^2 
		~~ \text{(thermodynamic model)}, 
	\end{cases}
\end{equation}
where
\begin{equation*}
	P = \bigl( \dfrac{x^2\bar\rho}{r^2r_x} \bigr)^{4/3} ~ \text{(isentropic model)} ~~ \text{or} ~~ K\dfrac{x^2\bar\rho\theta}{r^2r_x} ~ \text{(thermodynamic model)}.
\end{equation*}
Also, the boundary conditions \eqref{boundarycondition} can be written as
\begin{equation}
	\begin{gathered}
		r(0,t) = 0,~
		\mathfrak B = \dfrac{4}{3}\mu\bigl(\dfrac{r_{xt}}{r_x} - \dfrac{r_t}{r} \bigr)\biggr|_{x=R_0} = 0, ~ \theta(R_0, t) = 0 ~ \text{(thermodynamic model)}.
	\end{gathered}
\end{equation}

What to do next is to define the perturbation unknowns and write down the system satisfied by the perturbation variables. As one will see, the aforementioned damping structure will appear naturally, but with a temporal weight. In order to avoid dealing with the temporal weight, in \cite{Hadzic2016}, the authors introduce some new temporal variables for the self-similarly expanding solutions and the linearly expanding solutions respectively. 
We will adopt the same strategy in the following and write down the corresponding systems for the self-similarly expanding solution of the isentropic model, the linearly expanding solutions of the isentropic model and the thermodynamic model.

\subsection{Self-similarly expanding solutions of the isentropic model in the perturbation variable}

Considering the isentropic model for radiation gaseous stars, we define the perturbation variable as 
\begin{equation}\label{ptbv:isentropic}
 	\eta = \dfrac{r(x,t)}{r_\alpha(x,t)} - 1 = \dfrac{r(x,t)}{\alpha(t) x} - 1.
\end{equation}
Then after substituting \eqref{HomogeneousSol02} into \eqref{eq:Lagrangian} and using \eqref{ptbv:isentropic}, 
$ \eta $ satisfies the following equation,
\begin{equation}\label{eq:perturbation}
\begin{aligned}
	& \dfrac{\alpha^3(t)}{(1+\eta)^2} x \bar\rho \eta_{tt} + \dfrac{2\alpha^2(t)\alpha'(t)}{(1+\eta)^2} x \bar\rho \eta_t + \bigl( \dfrac{1}{1+\eta} - \dfrac{1}{(1+\eta)^4}\bigr) \delta x \bar\rho \\
	& ~~~~~~~~ + \bigl\lbrack\bigl(\dfrac{\bar\rho}{(1+\eta)^2(1+\eta+x\eta_x)}\bigr)^{4/3}\bigr\rbrack_x - \dfrac{(\bar\rho^{4/3})_x}{(1+\eta)^4} \\
	&  = \dfrac{4}{3}\mu \alpha^4(t) \bigl\lbrack \dfrac{(\alpha^3(t)x^3(1+\eta)^3)_{xt}}{(\alpha^3(t)x^3(1+\eta)^3)_x} \bigr\rbrack_x = \alpha^4(t) \mathfrak{B}_x + 4\mu\alpha^4(t)\bigl(\dfrac{\eta_t}{1+\eta}\bigr)_x,
\end{aligned}
\end{equation}
with the boundary condition
\begin{equation*}
	\mathfrak{B}(R_0,t) = 0, 
\end{equation*}
where
\begin{equation*}
	\mathfrak{B} := \dfrac{4}{3} \mu \bigl( \dfrac{\eta_t+x\eta_{xt}}{1+\eta + x \eta_{x}} - \dfrac{\eta_t}{1+\eta} \bigr).
\end{equation*}

Inspired by \cite{Hadzic2016}, we adopt the following formulation. For the self-similarly expanding solution (that is given by \eqref{Self-similarExpanding}), define the self-similar temporal variable $ s $ by
\begin{equation}\label{def:self-similar-time-isentropic}
	s = s(t) := \int_0^t \dfrac{1}{\alpha^{3/2}(\sigma)} \,d\sigma.
\end{equation}
Since we have, by \eqref{first-cosmology-speed} and \eqref{Self-similarExpanding},
$$ \delta < 0, ~ a_1 = \sqrt{\dfrac{2|\delta|}{a_0}} , ~ \alpha(t) = \bigl( a_0^{3/2} + 3 \sqrt{\dfrac{|\delta|}{2}}t \bigr)^{2/3},  $$
direct calculation yields
\begin{equation*}
	s(t) = \dfrac{\ln \alpha(t)}{\sqrt{2|\delta|}} - \dfrac{\ln a_0}{\sqrt{2|\delta|}}.
\end{equation*}
Denote
\begin{equation}\label{idi:ss=growingrate}
	\phi(x,s(t)) := \eta(x,t), ~~ \bar\alpha(s(t))  : = \alpha(t) = a_0 e^{\sqrt{2|\delta|}s}.
\end{equation}
Then \eqref{eq:perturbation} can be written as
\begin{equation}{\label{eq:perturbation=self-similar}}
	\begin{aligned}
		& \dfrac{1}{(1+\phi)^2}\bigl( x\bar\rho \phi_{ss} + \dfrac{b}{2} x\bar\rho \phi_s \bigr) + \bigl( \dfrac{1}{(1+\phi)^4} - \dfrac{1}{1+\phi} \bigr)|\delta| x\bar\rho \\
		& ~~~~~~ + \bigl\lbrack \bigl(\dfrac{\bar\rho}{(1+\phi)^2(1+\phi+x\phi_x)}\bigr)^{4/3}\bigr\rbrack_x - \dfrac{(\bar\rho^{4/3})_x}{(1+\phi)^4} \\
		&  = a_0^{5/2} e^{5/2 b s} \bigl(  \mathfrak{\bar B}_x + 4\mu \bigl(\dfrac{\phi_s}{1+\phi}\bigr)_x \bigr),
	\end{aligned}
\end{equation}
where $ b = \sqrt{2|\delta|} $ and
\begin{equation*}
	\mathfrak{\bar B} := \dfrac{4}{3} \mu \bigl( \dfrac{\phi_s+x\phi_{xs}}{1+\phi + x \phi_{x}} - \dfrac{\phi_s}{1+\phi} \bigr).
\end{equation*}
The boundary condition is given by
\begin{equation}\label{bc:self-similar}
	\mathfrak{\bar B}(R_0,s) = 0.
\end{equation}

\subsection{Linearly expanding solutions of the isentropic model in the perturbation variable}

Considering the linearly expanding solution (that is $ \alpha $ satisfies the asymptotic behavior \eqref{LinearlyExpanding}), define the linearly temporal variable $ \tau $ by
\begin{equation}\label{def:linearly-temporal-variable}
	\tau = \tau(t) = \int_0^t \dfrac{1}{\alpha(\sigma)}\,d\sigma.
\end{equation}
Denote
\begin{equation*}
	\vartheta(x,\tau(t)) : = \eta(x,t), ~~ \tilde \alpha(\tau(t)) : = \alpha(t).
\end{equation*}
Then $ \tilde \alpha $ satisfies 
\begin{equation}\label{def:growing-rate-isentropic}
	\tilde \alpha_\tau = \tilde \alpha \sqrt{a_1^2 + \dfrac{2\delta}{a_0} - \dfrac{2\delta}{\tilde\alpha}}, ~~ \text{and} ~ \tilde\alpha_{\tau\tau} = \delta + \tilde\alpha^{-1} \tilde\alpha_\tau^2 .
\end{equation}
Notice, $ \tilde\alpha(0) = a_0, \tilde\alpha_\tau(0) = a_0a_1 $.
It is easy to verify, (\cite[Sec 2.3.2]{Hadzic2016})
\begin{equation}\label{idi:le=growingrate}
	\begin{gathered}
		a_0 e^{\beta_1 \tau} \leq |\tilde \alpha| \leq a_0 e^{\beta_2 \tau}, ~ \beta_1 |\tilde\alpha| \leq |\tilde\alpha_\tau| \leq \beta_2 |\tilde\alpha|, \\
		\beta_1^2| \tilde \alpha| + \delta \leq |\tilde\alpha_{\tau\tau}| \leq |\delta| + \beta_2^2 |\tilde \alpha|, 
	\end{gathered}
\end{equation}
where \begin{equation}\label{le:growingspeed}
\beta _1 = \min\lbrace a_1, \sqrt{a_1^2 + \dfrac{2\delta}{a_0}} \rbrace  > 0  , \beta_2 = \max\lbrace a_1, \sqrt{a_1^2 + \dfrac{2\delta}{a_0}} \rbrace > 0. \end{equation}
In particular, when $ \delta = 0 $, $ \beta_1 = \beta_2 = a_1 $ and therefore,
\begin{equation*}
	e^{a_1 \tau} \lesssim |\tilde \alpha|+|\tilde\alpha_\tau| + |\tilde\alpha_{\tau\tau}| \lesssim e^{a_1 \tau}.
\end{equation*}
Then \eqref{eq:perturbation} can be written as
\begin{equation}\label{eq:perturbationrbation=linearly}
	\begin{aligned}
		& \dfrac{1}{(1+\vartheta)^2}\bigl( \tilde\alpha x\bar\rho \vartheta_{\tau\tau} + \tilde\alpha_\tau x\bar\rho \vartheta_\tau \bigr) + \bigl( \dfrac{1}{1+\vartheta} - \dfrac{1}{(1+\vartheta)^4} \bigr)\delta x\bar\rho \\
		& ~~~~ + \bigl\lbrack \bigl(\dfrac{\bar\rho}{(1+\vartheta)^2(1+\vartheta+x\vartheta_x)}\bigr)^{4/3}\bigr\rbrack_x - \dfrac{(\bar\rho^{4/3})_x}{(1+\vartheta)^4} 
		= \tilde\alpha^3 \bigl( \mathfrak{\tilde B}_x + 4 \mu \bigl( \dfrac{\vartheta_\tau}{1+\vartheta}\bigr)_x \bigr),
	\end{aligned}
\end{equation}
where
$$ \mathfrak{\tilde B} := \dfrac{4}{3} \mu \bigl( \dfrac{\vartheta_\tau+x\vartheta_{x\tau}}{1+\vartheta + x \vartheta_{x}} - \dfrac{\vartheta_\tau}{1+\vartheta} \bigr). $$
Similarly, the boundary condition is given by
\begin{equation}\label{BC:isentropic,linearlyexpanding}
\mathfrak{ \tilde B }(R_0,\tau) = 0.	
\end{equation}

\subsection{Linearly expanding solutions of the thermodynamic model in the perturbation variables}
Considering the linearly expanding solution of the thermodynamic model, again, denote the perturbation variables as
\begin{equation}
	\eta = \dfrac{r(x,t)}{r_\alpha(x,t)} - 1 = \dfrac{r(x,t)}{\alpha(t)x} - 1, ~ \varsigma = \alpha (\theta - \theta_\alpha) = \alpha(t) \theta - \bar\theta.
\end{equation}
Then similar as in the isentropic case, after employing the equations \eqref{HomogenerousSol04}, the system \eqref{eq:Lagrangian} in terms of the perturbation variables $ \lbrace \eta, \varsigma \rbrace $ is in the following form,
\begin{equation}\label{eq:RGS-perterbation}
	\begin{cases}
		\dfrac{\alpha^3(t)}{(1+\eta)^2} x\bar\rho \eta_{tt} + \dfrac{2\alpha^2(t)\alpha'(t)}{(1+\eta)^2} x\bar\rho \eta_t + \bigl\lbrack \dfrac{K\bar\rho(\varsigma + \bar\theta)}{(1+\eta)^2(1+\eta+x\eta_x)} \bigr\rbrack_x - \dfrac{(K\bar\rho\bar\theta)_x}{(1+\eta)^4} \\
		~~~~~~ = \dfrac{4}{3}\mu \alpha^4(t) \bigl\lbrack \dfrac{(\alpha^3(t)x^3(1+\eta)^3)_{xt}}{(\alpha^3(t)x^3(1+\eta)^3)_x} \bigr\rbrack_x = \alpha^4(t) \mathfrak{B}_x + 4\mu\alpha^4(t)\bigl(\dfrac{\eta_t}{1+\eta}\bigr)_x, \\
		3K  x^2 \bar\rho \varsigma_t + K  \dfrac{\bar\rho(\varsigma + \bar\theta) (x^3(1+\eta)^2 \eta_t)_x}{(1+\eta)^2(1+\eta+x\eta_x)} - \alpha(t) \bigl\lbrack \dfrac{(1+\eta)^2}{1+\eta+x\eta_x} x^2 \varsigma_x + \bigl( \dfrac{(1+\eta)^2}{1+\eta+x\eta_x}\\
		~~~~~~ - 1 \bigr)x^2 \bar\theta_x \bigr\rbrack_x = \alpha^2(t) x^2 (1+\eta)^2 (1+\eta+x\eta_x) \cdot \mathfrak F(\eta), 
		\end{cases}
\end{equation}
where
\begin{equation*}
	\begin{aligned}
	\mathfrak B & := \dfrac{4}{3}\mu \bigl( \dfrac{\eta_t + x \eta_{xt}}{1+\eta+x\eta_x}  - \dfrac{\eta_t}{1+\eta}\bigr) , \\
	\mathfrak F(\eta) 
	& := \dfrac{4}{3}\mu \bigl\lbrack \dfrac{\eta_t+x\eta_{xt}}{1+\eta+x\eta_x} - \dfrac{\eta_t}{1+\eta} \bigr\rbrack^2.
	\end{aligned}
\end{equation*}
The associated boundary condition is then given by
\begin{equation*}
	\mathfrak B(R_0,t) = 0, ~ \varsigma(R_0,t) = 0. 
\end{equation*}

Notice, $ \alpha $ is growing linearly over time. Similar as before, it will be convenient to work with the linearly temporal variable $ \tau $ defined by \eqref{def:linearly-temporal-variable}
\begin{equation*}\label{def:ln-tm-variable}  \tag{\ref{def:linearly-temporal-variable}}
	\tau = \tau(t) = \int_0^t \dfrac{1}{\alpha(\sigma)} \, d\sigma = \dfrac{\ln{(1+a_1t/a_0)}}{a_1}.
\end{equation*}
Consequently, 
$$ t(\tau) = \dfrac{a_0}{a_1}(e^{a_1\tau} -1 ), ~ \dfrac{d}{d\tau} t = a_0 e^{a_1\tau} = \tilde \alpha(\tau), ~ \tilde \alpha(\tau(t)) : = \alpha(t) = a_0 e^{a_1\tau} . $$
Here we employ the same notation for the new temporal variable as in the linearly expanding solution for the isentropic model since asymptotically they grow at the same rate. By denoting the unknowns in the linearly temporal variable as
\begin{equation}\label{def:pt-variable-thermo}
	\xi(x,\tau(t)) := \eta(x,t), ~ \zeta(x,\tau(t)) : = \varsigma(x,t),
\end{equation}
the equations \eqref{eq:RGS-perterbation} take the form,
\begin{equation}\label{eq:RGS-ptb-ln-tm-vb}
\begin{cases}
	\dfrac{1}{(1+\xi)^2} \bigl(\tilde\alpha x\bar\rho \xi_{\tau\tau} + \tilde\alpha_\tau x\bar\rho\xi_\tau\bigr) + \bigl\lbrack \dfrac{K\bar\rho(\zeta + \bar\theta)}{(1+\xi)^2(1+\xi+x\xi_x)} \bigr\rbrack_x - \dfrac{(K\bar\rho\bar\theta)_x}{(1+\xi)^4} = \tilde\alpha^3 \bigl( \mathfrak{\hat B}_x + 4 \mu \bigl( \dfrac{\xi_\tau}{1+\xi}\bigr)_x \bigr), \\
	3K x^2 \bar\rho \zeta_\tau + K \dfrac{\bar\rho(\zeta+\bar\theta)(x^3(1+\xi)^2\xi_\tau)_x}{(1+\xi)^2(1+\xi+x\xi_x)} - \tilde\alpha^2 \bigl\lbrack \dfrac{(1+\xi)^2}{1+\xi+x\xi_x}x^2\zeta_x + \bigl(\dfrac{(1+\xi)^2}{1+\xi+x\xi_x} - 1\bigr) x^2 \bar\theta_x \bigr\rbrack_x \\
	~~~~~~ = \tilde\alpha x^2(1+\xi)^2(1+\xi+x\xi_x)\cdot\mathfrak{\hat F}(\xi),
\end{cases}	
\end{equation}
where 
\begin{equation}\label{BC:thermodynamic}
\begin{aligned} 
	\mathfrak{\hat B} & := \dfrac{4}{3}\mu \bigl( \dfrac{\xi_\tau+x\xi_{x\tau}}{1+\xi + x \xi_{x}} - \dfrac{\xi_{\tau}}{1+\xi}\bigr), \\
	\mathfrak{\hat F}(\xi) & : = \dfrac{4}{3}\mu \bigl\lbrack \dfrac{\xi_\tau+x\xi_{x\tau}}{1+\xi+x\xi_x} - \dfrac{\xi_\tau}{1+\xi}\bigr\rbrack^2.
\end{aligned}	
\end{equation}
Now the boundary conditions take the forms
$$ \mathfrak{\hat B}(R_0,\tau) = 0, ~ \zeta(R_0,\tau) = 0. $$

\subsection{Comments and methodology}

For the self-similarly expanding solutions, we shall analyse the energy. Since the energy is exactly zero for such solutions, a small perturbation and the non-increasing property of the physical energy will lead to the instability.

As for the linearly expanding solution of the isentropic model \eqref{eq:perturbationrbation=linearly}, it is worth to look at some linear model equations. First, consider the following linear structure associated with \eqref{eq:perturbationrbation=linearly},
\begin{equation*}
	e^{k\tau}(f_{\tau\tau} + k f_{\tau}) + \delta f - (\bar\rho^{4/3}f_x)_x = e^{3k\tau} f_{xx\tau},
\end{equation*}
where $ k > 0 $ is a constant. Then the $ L^2 $-estimate of this model equation is of the form
\begin{align*}
	& \dfrac{e^{k\tau}}{2} \norm{f_\tau}{\stnorm{\infty}{2}}^2 + \dfrac{1}{2} \norm{\bar\rho^{2/3}f_x}{\stnorm{\infty}{2}}^2 +  \dfrac{\delta}{2} \norm{f}{\stnorm{\infty}{2}}^2 \\
	& ~~~~   + \int \dfrac{ke^{k\tau}}{2} \norm{f_\tau}{L_x^2}^2 \,d\tau + \int e^{3k\tau} \norm{f_{x\tau}}{L_x^2}^2 \,d\tau \leq \text{initial data}.
\end{align*}
If $ \delta < 0 $, we use the fundamental theorem of calculus to derive that
\begin{equation*}
	\norm{f}{\stnorm{\infty}{2}} \leq \int \norm{f_\tau}{L_x^2}\,d\tau + \cdots \leq (\int e^{-k\tau}\,d\tau)^{1/2} \cdot (\int e^{k\tau}\norm{f_\tau}{L_x^2}^2 \,d\tau)^{1/2} + \cdots,
\end{equation*}
and therefore
\begin{align*}
	\cdots + \dfrac{\delta}{2}\norm{f}{\stnorm{\infty}{2}}^2 + \int \dfrac{ke^{k\tau}}{2} \norm{f_\tau}{L_x^2}^2 \,d\tau \geq  (\delta \int e^{-k\tau}\,d\tau + \dfrac{k}{2}) \int e^{k\tau} \norm{f_\tau}{L_x^2}^2 \,d\tau + \cdots.
\end{align*}
Then $ \delta $ has to been chosen such that the coefficient on the right is positive. This is how the constraint \eqref{thm1.2-001} comes out. Also, the damping structure in the above equation will imply a faster decay. On the other hand, the following linear structure is associated with the temporal derivative of \eqref{eq:perturbationrbation=linearly},
\begin{equation*}
	e^{k\tau}(f_{\tau\tau\tau} + k f_{\tau\tau} ) = e^{3k\tau} f_{xx\tau\tau} + e^{3k\tau} f_{x\tau} + \cdots.
\end{equation*}
Then the $ L^2 $-estimate will yields, 
\begin{align*}
	& \dfrac{e^{k\tau}}{2}\norm{f_{\tau\tau}}{\stnorm{\infty}{2}}^2 + \int \dfrac{ke^{k\tau}}{2} \norm{f_{\tau\tau}}{L_x^2}^2 \,d\tau + \int \dfrac{e^{3k\tau}}{2} \norm{f_{x\tau\tau}}{L_x^2}^2\,d\tau \\
	& ~~~ \leq \int \dfrac{e^{3k\tau}}{2}\norm{f_\tau}{L_x^2}^2 \,d\tau + \cdots,
\end{align*}
where we can not bound the first integral on the right since we only have the bound of $ \int e^{k\tau}\norm{f_\tau}{L_x^2}^2 \,d\tau $ from the previous analysis. Therefore we shall apply the negative temporal weight to manipulate the estimate. Such estimates will control the growth of the higher order norms. Eventually, as long as the growth is not too large, we can use the elliptic structure of the equation to recover the spatial regularity and estimates. 
We remark that such a structure is a consequence of the monoatomic gas viscous coefficients \eqref{weakviscosity}.

On the other hand, when studying the thermodynamic model \eqref{eq:RGS-ptb-ln-tm-vb}, besides the similar structure as in the isentropic case for \subeqref{eq:RGS-ptb-ln-tm-vb}{1} , the equation \subeqref{eq:RGS-ptb-ln-tm-vb}{2} has the form
\begin{equation*}
	g_\tau - e^{2k\tau} g_{xx} = e^{k\tau}(f_\tau^2 + f_{x\tau}^2) + e^{2k\tau} (f_{xx} + f_x) + \cdots.
\end{equation*}
Then similar $ L^2 $-estimate yields
\begin{equation*}
	\dfrac{1}{2}\norm{g}{\stnorm{\infty}{2}}^2 +  \int \dfrac{e^{2k\tau}}{2} \norm{g_x}{L_x^2}^2 \,d\tau \leq \int \dfrac{e^{2k\tau}}{2} \norm{f}{L_x^2}^2 \,d\tau + \cdots,
\end{equation*}
where we have no estimate on the first integral on the right. However, we know, similarly as in the isentropic case, $ \norm{f}{\stnorm{\infty}{2}} $ will be bounded. Therefore by choosing the temporal weight $ e^{-\kappa \tau} $ with $ \kappa > 0 $ large enough, we can obtain the temporal weighted estimate. Such structures exist for all the higher order estimates. We carefully track the temporal weights through out our analysis and eventually perform elliptic estimates to establish the spatial regularity.  

Unless stated otherwise, we adopt the following notations in the above and the rest of this work:
\begin{equation*}
	\int \cdot \,dx = \int_0^{R_0} \cdot \,dx,~ \int \cdot \,dt = \int_0^T \cdot \,dt,
\end{equation*}
where $ 0 < R_0, T < \infty $ are fixed positive constants and our solutions live in the space time domain $ (0,R_0) \times (0,T) $. Also, $ \norm{\cdot}{\stnorm{p}{q}} $ denotes the standard space-time Sobolev norm in the space time variable $ (x,t) \in (0, R_0) \times (0,T) $. For any fixed time $ t \in (0,T) $, $ \norm{\cdot}{L_x^{q}} $ is the standard Sobolev norm in the space variable $ x\in (0,R_0) $. Here $ p,q \in (1,\infty] $. $ \sup_t $ will be used to represent $ \sup_{0< t < T} $. Similar notations are also adopted for the rescaled temporal variable $ s, \tau $ defined above. $ A \lesssim B $ is used to denote that there exists some positive constant $ C $ such that $ A \leq C\cdot B $. $ A \simeq B $ will mean $ A \lesssim B $ and $ B \lesssim A $. $ C = C(\cdot) $ is a constant which is different from line to line and depends on the inputs. The following form of the Hardy's inequality will be employed in this work.

\begin{lm}[Hardy's inequality, \cite{Jang2014}]
	Let $ k $ be a given real number, and let $ g $ be a function satisfying $ \int_0^1 s^k (g^2 + g'^2)\,ds < \infty $.
	\begin{enumerate}
		\item If $ k > 1$, then we have
		\begin{equation*}
			\int_0^1 s^{k-2}g^2 \,ds \leq C \int_0^1 s^k (g^2+g'^2 )\,ds.
		\end{equation*}
		\item If $ k < 1 $, then $ g $ has a trace at $ x=0 $ and moreover
		\begin{equation*}
			\int_0^1 s^{k-2}(g-g(0))^2\,ds\leq C\int_0^1 s^k g'^2 \,ds.
		\end{equation*}
	\end{enumerate}
\end{lm}

In this work, inspired by \cite{LuoXinZeng2016} let us define the relative entropy functional as 
\begin{equation}\label{relative-entropy-functional}
\mathfrak H(h) : = \log(1+h)^2(1+h+xh_x),
\end{equation}
where $ h: (0,R_0)\times(0,T) \mapsto \mathbb R $ is any smooth function. 
We will have the following estimates on the relative entropy.
\begin{lm}\label{lm:estimates-of-relative-entropy}
For $ h $ satisfies, 
\begin{equation}\label{lm:s003} 
\max \lbrace \norm{h}{\stnorm{\infty}{\infty}}, \norm{xh_x}{\stnorm{\infty}{\infty}}, \norm{h_t}{\stnorm{\infty}{\infty}},  \norm{xh_{xt}}{\stnorm{\infty}{\infty}} 
\rbrace < \varepsilon,
\end{equation} 
with some $ 0 < \varepsilon  < 1 $ small enough, the following estimates of the function $ \mathfrak H (h) $ hold,
\begin{gather}
\int ( h_x^2 + x^2 h_{xx}^2 ) \,dx  \lesssim  \int \mathfrak H(h)_x^2 \,dx ,  \label{lm:s001}  \\
\int (h_{xt}^2 + x^2 h_{xxt}^2) \,dx \lesssim \int \mathfrak H(h)_{xt}^2 \,dx + \varepsilon \int (h_x^2 + x^2 h_{xx}^2 )\,dx . \label{lm:s002}
\end{gather}
\end{lm}

\begin{pf}
Notice, \begin{equation*}
\begin{aligned}
\mathfrak H(h)_x &  = \dfrac{1}{(1+h)(1+h+xh_x)} \bigl\lbrace 2(1+h+xh_x) h_x + (1+h)(2h_x + xh_{xx}) \bigr\rbrace\\
& = \dfrac{4h_x + xh_{xx}}{(1+h)(1+h+xh_x)} + \dfrac{2(h+xh_x)h_x + h (2h_x + xh_{xx})}{(1+h)(1+h+xh_x)}.
\end{aligned}
\end{equation*}
Under the assumption \eqref{lm:s003}, we have
\begin{equation*}
\int \mathfrak H(h)_x^2 \,dx \gtrsim \dfrac{1}{1+\varepsilon} \underbrace{ \int ( 4h_x + xh_{xx} )^2\,dx }_{\mathfrak A} - \varepsilon \int h_x^2 + x^2 h_{xx}^2 \,dx.
\end{equation*}
$ \mathfrak A $ can be calculated as follows,
\begin{align*}
\mathfrak  A & = 16\int h_x^2 \,dx + \int x^2 h_{xx}^2 \,dx + 8\int x  h_x h_{xx}\,dx = 16 \int h_x^2 \,dx + \int x^2 h_{xx}^2 \,dx - 4 \int h_x^2 \,dx + 4 xh_{x}^2 \bigr|_{x=R_0} \\
& \geq 12 \int h_x^2 \,dx + \int x^2 h_{xx}^2 \,dx. 
\end{align*}
Therefore, \eqref{lm:s001} follows easily. Similarly, 
\begin{equation*}
\begin{aligned}
\mathfrak H(h)_{xt} & = \dfrac{4h_{xt} + xh_{xxt}}{(1+h)(1+h+xh_x)} + (4h_x+xh_{xx})\bigl\lbrace\dfrac{1}{(1+h)(1+h+xh_x)} \bigr\rbrace_t \\
& ~~~~~~ + \bigl\lbrace \dfrac{2(h+xh_x)h_x + h (2h_x + xh_{xx})}{(1+h)(1+h+xh_x)} \bigr\rbrace_t,
\end{aligned}
\end{equation*}
from which we have
\begin{equation*}
\begin{aligned}
\int \mathfrak H(h)_{xt}^2 \,dx & \gtrsim \dfrac{1}{1+\varepsilon} \int (4 h_{xt} + xh_{xxt})^2 \,dx - \varepsilon \int ( h_x^2 + x^2 h_{xx}^2 + h_{xt}^2 + x^2h_{xxt}^2 ) \,dx\\
& \gtrsim \dfrac{1}{1+\varepsilon} \int (12 h_{xt}^2 + x^2 h_{xxt} )\,dx  - \varepsilon \int ( h_x^2 + x^2 h_{xx}^2 + h_{xt}^2 + x^2h_{xxt}^2 ) \,dx.
\end{aligned}
\end{equation*}
This finishes the proof.
\end{pf}


\section{Instability of Self-similarly Expanding Homogeneous Solutions for the Isentropic Model}\label{Section:instability}

We start with the property of the self-similarly expanding homogeneous solutions of the isentropic model. It is easy to verify from \eqref{first-cosmology-speed} and \eqref{energy-inviscid} that the energy for such a solution vanishes.

From \eqref{energy-eulercoordinate} and \eqref{weakviscosity}, $ D = 0 $ if and only if $ u \propto r $ in $ (0, R(t)) $. Therefore, unless we adopt a perturbation $ \phi $ with $ \phi \propto r $, the physical energy $ E $ defined in \eqref{energy-eulercoordinate} decreases over time. For the instability theory, we will show that, in contrast to \cite{Hadzic2016}, the vanishing of the initial energy is not enough to guarantee the stability. To do so, we shall write down the energy $ E $ in term of the perturbation variable $ \phi $ in the Lagrangian coordinates. We will use the notation $ E(\phi,\phi_s) = E(t) = E(s) $ and $ (\phi_0, \phi_1) = (\phi(\cdot,0), \phi_s(\cdot,0)) $ to denote the physical energy and the initial perturbation. Similarly $ D(\phi,\phi_s) = D(t) = D(s) $ denotes the dissipation. 
In particular, $ E(\phi_0,\phi_1) $ is the initial energy corresponding to the initial perturbation $ (\phi_0, \phi_1) $. Here $ \phi $ is defined in \eqref{ptbv:isentropic} and \eqref{idi:ss=growingrate}, $(x,t)$ are the Lagrangian variables defined through \eqref{Lagrangian001} and \eqref{Lagrangian002}, and $ s $ is the self-similar temporal variable defined in \eqref{def:self-similar-time-isentropic}. Then $ E $ can be written as, using \eqref{HomogeneousSol02} and \eqref{HomogeneousSol03},
\begin{equation}\label{energy-ss-perturbation}
\begin{aligned}
	& E(t) = \dfrac{1}{2}\int x^2 \bar\rho r_t^2 \,dx + \delta\int \dfrac{x^5}{r} \bar\rho \,dx + \int x^2 \bar\rho^{4/3} \bigl\lbrace 3 \dfrac{x^{2/3}}{(r^2r_x)^{1/3}} - 4 \dfrac{x}{r} + \dfrac{x^2r_x}{r^2} \bigr\rbrace \,dx\\
	& ~~ = \bar\alpha^{-1}\int x^4 \bar\rho \bigl( \dfrac{1}{2} \phi_s^2 + \sqrt{2|\delta|} (1+\phi)\phi_s - \delta (1+\phi)^2 + \dfrac{\delta}{1+\phi} \bigr) \,dx \\
	& ~~~~~~ + \bar\alpha^{-1}\int x^2 \bar\rho^{4/3} \bigl\lbrace \dfrac{3}{((1+\phi)^2(1+\phi+x\phi_x))^{1/3}} -\dfrac{3}{1+\phi} + \dfrac{x\phi_x}{(1+\phi)^2} \bigr\rbrace \,dx  = E(s).
\end{aligned}
\end{equation}
Similarly,
\begin{equation}\label{dissipation-ss-perturbation}
	\begin{aligned}
		& D(t)  = \int \bigl\lbrack 2\mu \bigl( r^2 \dfrac{r_{xt}^2}{r_x^2} + 2 r_t^2 \bigr) - \dfrac{2\mu}{3} \bigl(r \dfrac{r_{xt}}{r_x} + 2 r_t \bigr)^2 \bigr\rbrack r_x \,dx \\
		& ~~ = \dfrac{4\mu}{3} \int x^2 \dfrac{((1+\phi)x\phi_{xs} - x\phi_x \phi_s)^2}{1+\phi+x\phi_x} \,dx = D(s),
	\end{aligned}
\end{equation}
and 
\begin{equation}\label{energy-ss-conservation}
	\dfrac{d}{ds} E + \bar\alpha^{3/2} D = 0.
\end{equation}

Consider the perturbation given by $ \phi = \phi(s) $ as a function independent of the spatial variable, which satisfies
\begin{equation}\label{eq:ss=001}
	(\phi_s + \sqrt{2|\delta|}(1+\phi))^2 + \dfrac{2\delta}{1+\phi} = 0.
\end{equation}
It is easy to verify, for such perturbations, $ E \equiv 0 $ and $ D \equiv 0 $. 
Direct calculation yields the zero of \eqref{eq:ss=001} on the $ (\phi,\phi_s)$-plane which contains $(0,0)$ (corresponding to the self-similarly expanding homogeneous solution) is given by
\begin{equation}\label{eq:ss=002}
	\phi_s = -\sqrt{2|\delta|} (1+\phi) + \sqrt{2|\delta|} (1+\phi)^{-1/2},
\end{equation}
and it is directly verified that $ (0,0) $ is the steady point of \eqref{eq:ss=001}. 
Meanwhile, it is easy to calculate the zeros of \eqref{eq:ss=001} is given by \eqref{eq:ss=002} and $ \phi_s = -\sqrt{2|\delta|} (1+\phi) - \sqrt{2|\delta|} (1+\phi)^{-1/2} $. 
Also, on the half plane $ \lbrace \phi > -1 \rbrace \subset \lbrace (\phi,\phi_s) \rbrace $, we have, 
\begin{equation*}
	\begin{aligned}
		E(\phi,\phi_s) = 0 & ~~~~ \text{if} ~ (\phi,\phi_s) ~ \text{on the curve determined by \eqref{eq:ss=001}},\\
		E(\phi,\phi_s) < 0 & ~~~~ \text{if} ~ -\sqrt{2|\delta|} (1+\phi) - \sqrt{2|\delta|}(1+\phi)^{-1/2} < \phi_s < -\sqrt{2|\delta|} (1+\phi) + \sqrt{2|\delta|}(1+\phi)^{-1/2},\\
		E(\phi,\phi_s) > 0 & ~~~~ \text{otherwise}.\\
	\end{aligned}
\end{equation*} 

\subsection*{Homogeneous perturbations}
Inspired by the above analysis, we shall consider the perturbations of the form $ \phi = \phi(s) $ first. For the sake of convenience, such a kind of perturbations will be referred to as the homogeneous perturbations. Now, from the conservation of energy which is due to the fact $ D \equiv 0 $ for the homogeneous perturbations, together with \eqref{eq:perturbation=self-similar}, $ \phi $ satisfies the following equations,
\begin{gather}
	\phi_{ss} + \dfrac{b}{2}\phi_s + |\delta| \biggl( \dfrac{1}{(1+\phi)^2} - (1+\phi) \biggr) = 0, \label{eq:homogeneousptb-1} \\
	\dfrac{1}{2} (\phi_s + \sqrt{2|\delta|}(1+\phi))^2 + \dfrac{\delta}{1+\phi} = \bar\alpha \dfrac{E(\phi,\phi_s)}{\int x^4 \bar\rho \,dx } = \bar\alpha \dfrac{E(\phi_0,\phi_1)}{\int x^4 \bar\rho\,dx }, \label{eq:homogeneousptb-2}
\end{gather}
where $ (\phi_0,\phi_1) = (\phi(0), \phi_s(0)) $ and $ \phi_0 > -1 $. Indeed, $ \phi $ is determined by the equation \eqref{eq:homogeneousptb-1}, and \eqref{eq:homogeneousptb-2} is the consequence of \eqref{eq:homogeneousptb-1}. The constraint $ \phi_0 > -1 $ is to say that the perturbation is well-defined (see \eqref{ptbv:isentropic} and \eqref{idi:ss=growingrate}).

Here one can perform the phase space analysis of \eqref{eq:homogeneousptb-1} on the $(\phi,\phi_s) $ plane. Indeed, \eqref{eq:homogeneousptb-1} can be written as
\begin{equation}
	\dfrac{d}{ds}\biggl(\begin{array}{c}
		\phi \\ \phi_s
	\end{array}\biggr) = \biggl(\begin{array}{c}
		\phi_s \\ - \dfrac{b}{2} \bigl\lbrack\phi_s + \sqrt{2|\delta|}\bigl(\dfrac{1}{(1+\phi)^2} - (1+\phi)\bigr)\bigr\rbrack
	\end{array} \biggr).
\end{equation}
Notice first, it is easy to verified,
\begin{equation}\label{ss:002}
	(\phi_s + \sqrt{2|\delta|}((1+\phi)-(1+\phi)^{-1/2}))_s = \dfrac{\sqrt{2|\delta|}}{2}(1+(1+\phi)^{-3/2})(\phi_s + \sqrt{2|\delta|}((1+\phi)-(1+\phi)^{-1/2})).
\end{equation}

\noindent {\bf Case 1} $ \phi_1 > -\sqrt{2|\delta|}((1+\phi_0) - (1+\phi_0)^{-1/2}) $.  
From \eqref{ss:002}
\begin{equation}\label{ss:001}
	\phi_s > - \sqrt{2|\delta|}((1+\phi)-(1+\phi)^{-1/2}) + e^{\frac{\sqrt{2|\delta|}}{2}\cdot s}(\phi_1 + \sqrt{2|\delta|} ((1+\phi_0) - (1+\phi_0)^{-1/2})) .
\end{equation}

\begin{itemize}
\item  If $ \phi_s \leq 0 $ in $ s \in [0, s_1) $, then $ \phi \leq  \phi(0) = \phi_0 $. From \eqref{ss:001}, $ \phi_s > - \sqrt{2|\delta|}((1+\phi_0)-(1+\phi_0)^{-1/2}) + e^{\frac{\sqrt{2|\delta|}}{2}\cdot s}(\phi_1 + \sqrt{2|\delta|} ((1+\phi_0) - (1+\phi_0)^{-1/2})) $. Therefore, $ \exists \, s_1 < \infty $ such that $ \phi_s (s_1) > 0 $.

\item If $ \phi_s > 0 $ and $\phi \leq 0 $ in $ s\in [s_1,s_2) $, from \eqref{ss:001}, $ \phi_{s} > \phi_1 + \sqrt{2|\delta|} ((1+\phi_0) - (1+\phi_0)^{-1/2}):=c > 0   $. Therefore, $ \exists \, s_1 < s_2 < \infty $ such that $ \phi(s_2) >c (s_2-s_1) + \phi(s_1) > 0 $.

\item If $\phi_s > 0 $, $ \phi > 0$,  $\phi_s > \sqrt{2|\delta|}\bigl((1+\phi) - \dfrac{1}{(1+\phi)^2}\bigr)$ in $ s \in [s_2,s_3) $, then $ \phi > \phi(s_2) > 0 $ and $ \phi_{ss} < 0 $. Therefore $ \phi_s(s_2) > \phi_s > \sqrt{2|\delta|}\bigl((1+\phi(s_2)) - \dfrac{1}{(1+\phi(s_2))^2} \bigr) :=c > 0 $. Consequently, $ \phi > c (s-s_2) + \phi(s_2) $. Suppose $ s_3 = \infty $ then $ \exists \, s_3' <\infty $ large enough such that $ \phi_s(s_3') > \sqrt{2|\delta|}\bigl((1+\phi(s_3')) - \dfrac{1}{(1+\phi(s_3'))^2}\bigr) > \phi_s(s_2)  $ which is a contradiction. Hence $ s_3 < \infty $.

\item If $\phi_s > 0 $, $ \phi > 0$,  $\phi_s \leq \sqrt{2|\delta|}\bigl((1+\phi) - \dfrac{1}{(1+\phi)^2}\bigr)$ in $ s \in [s_3,s_4) $, then $ \phi_{ss} \geq 0 $. We claim $ s_4 = \infty $. Otherwise, $ \exists \, s_4'< s_4'' < \infty $ such that $ \phi_s(s_4') = \sqrt{2|\delta|}\bigl((1+\phi(s_4') - \dfrac{1}{(1+\phi(s_4')^2}\bigr) > 0 $ but $ \phi_s > \sqrt{2|\delta|}\bigl((1+\phi) - \dfrac{1}{(1+\phi)^2}\bigr)>0 $, $ \phi_{ss} < 0 $ in $ s \in (s_4',s_4''] $. Consequently, $ \phi_{s}(s_4'') < \phi_s(s_4') = \sqrt{2|\delta|}\bigl((1+\phi(s_4')) - \dfrac{1}{(1+\phi(s_4'))^2}\bigr) < \sqrt{2|\delta|} \bigl((1+\phi(s_4'')) - \dfrac{1}{(1+\phi(s_4''))^2}\bigr) < \phi_s(s_4'') $, which is a contradiction. Therefore $ s_4 = \infty $. Moreover, it is easy to derive, $ \phi_s > \phi_s(s_3) > 0 $ and $ \phi(s) \rightarrow \infty $ as $ s \rightarrow \infty $. In the meantime, \eqref{eq:homogeneousptb-1} implies $ \phi_{ss} + \dfrac{b}{2}\phi_s \rightarrow \infty $ as $ s \rightarrow \infty $. Therefore $ \phi_s \rightarrow \infty $ as $ s \rightarrow \infty $. 

\end{itemize}

\todo{9 Sep 2017}

\noindent {\bf Case 2} $ \phi_1 < -\sqrt{2|\delta|}((1+\phi_0) - (1+\phi_0)^{-1/2}) $. Similar as before, from \eqref{ss:002},
\begin{equation}\label{ss:003}
	\phi_s < - \sqrt{2|\delta|}((1+\phi)-(1+\phi)^{-1/2}) + e^{\frac{\sqrt{2|\delta|}}{2} \cdot s} (\phi_1 + \sqrt{2|\delta|} ((1+\phi_0) - (1+\phi_0)^{-1/2})) .
\end{equation}
\begin{itemize}
	\item If $ \phi_s \geq 0 $ in $ s \in [0, s_1) $, then $ \phi \geq \phi(0) = \phi_0 $. Therefore from \eqref{ss:003}, $ \phi_s < -\sqrt{2|\delta|} ((1+\phi(0)) - (1+\phi(0))^{-1/2}) + e^{\frac{\sqrt{2|\delta|}}{2}\cdot s }(\phi_1 + \sqrt{2|\delta|}((1+\phi_0)-(1+\phi_0)^{-1/2})) $. Then $ \exists \, s_1 <\infty $ such that $ \phi_s(s_1) < 0 $. 
	\item If $ \phi_s < 0 $,  $ \phi \geq 0 $ in $ s\in [s_1,s_2) $, from \eqref{ss:003}, $ \phi_s < \phi_1+\sqrt{2|\delta|}((1+\phi_0)-(1+\phi_0)^{-1/2}): = -c < 0 $ Therefore, $ \exists s_1 < s_2 < \infty $ such that $ \phi(s_2)< -c(s_2-s_1) + \phi(s_1) <0 $.
	\item If $ \phi_s < 0 $, $ \phi < 0 $, $ \phi_s \leq \sqrt{2|\delta|}\bigl( (1+\phi) - \dfrac{1}{(1+\phi)^2}\bigr) $ in $ s \in [s_2, s_3) $, then $ \phi < \phi(s_2) < 0 $, $ \phi_{ss} \geq 0 $ and $ \phi_s \geq \phi_s(s_2) $. Therefore $ \phi_s < \sqrt{2|\delta|} \bigl( (1+\phi(s_2)) - \dfrac{1}{(1+\phi(s_2))^2}\bigr) :=-c < 0 $ and $ \phi < - c(s-s_2) + \phi(s_2) $. In particular, $ \phi_s(s_2) \leq \phi_s(s_3) \leq \sqrt{2|\delta|}\bigl( (1+\phi(s_3)) - \dfrac{1}{(1+\phi(s_3))^2}\bigr) < \sqrt{2|\delta|}\bigl( (1-c(s_3-s_2) + \phi(s_2)) - \dfrac{1}{(1-c(s_3-s_2) + \phi(s_2))^2}\bigr)  $ which yields $ s_3 < \infty $.
	\item If $ \phi_s < 0 $, $\phi< 0$, $ \phi_s > \sqrt{2|\delta|}\bigl( (1+\phi) - \dfrac{1}{(1+\phi)^2}\bigr) $ in $ s \in [s_3,s_4) $, then from \eqref{eq:homogeneousptb-1} $ \phi_{ss} < 0 $. $ \phi_s < \phi_s(s_3) < 0 $. We claim $ s_4 < \infty $. Otherwise, since, similar as before, $\phi < \phi_s(s_3)(s-s_3) + \phi(s_3) $, we will have $ \exists \, s_4' < \infty $ such that $ \phi \rightarrow -1 $ as $ s \rightarrow s_4' $. This is a contradiction.
\end{itemize}

\begin{rmk}
	Indeed, the homogeneous perturbations in Case 1 and Case 2 discussed above correspond to the solutions with velocity larger or smaller than the escape velocity given
	in Section \ref{section:introduction_of_expanding_sol}.
\end{rmk}

\subsection* {Inhomogeneous perturbations} 
We will briefly demonstrate in the following, for any inhomogeneous initial perturbation (non-constant) $ (\phi_0,\phi_1) = (\phi_0(x),\phi_1(x))  $ with $ E(\phi_0,\phi_1) \leq 0 $ and $ D(\phi_0, \phi_1) \neq 0 $, there is no $ 0.5 > \bar\varepsilon_0 > 0 $ (no matter how small the perturbation is) such that $ \max \lbrace \norm{\phi}{L_x^\infty}, \norm{x\phi_x}{L_x^\infty} \rbrace < \bar\varepsilon_0  $ for all $ s > 0 $. 
From the energy equality \eqref{energy-ss-conservation}, there exists $ \bar s $ such that $ E(\bar s) < 0 $.
 Hence, without loss of generality, it is assumed $ E_0 = E(\phi_0, \phi_1) < 0 $. Then from the definition of $ E $ and \eqref{energy-ss-conservation}, the following holds
\begin{equation}\label{ss-contradict-001}
	- \int x^4\bar\rho \dfrac{|\delta|}{1+\phi}\,dx - \int x^2 \bar\rho^{4/3} \bigl( \dfrac{3}{1+\phi} - \dfrac{x\phi_x}{(1+\phi)^2} \bigr) \,dx  \leq \bar\alpha E_0 < 0.
\end{equation}
Suppose $ \exists \, 0 < \bar\varepsilon_0 < 0.5 $ such that $ \max \lbrace \norm{\phi}{\infty}, \norm{x\phi_x}{\infty} \rbrace < \bar\varepsilon_0 $ for all $ s > 0 $. Then from \eqref{ss-contradict-001}, 
\begin{equation*}
	\bar\alpha(s) |E_0| \leq \abs{ - \int x^4\bar\rho \dfrac{|\delta|}{1+\phi}\,dx - \int x^2 \bar\rho^{4/3} \bigl(  \dfrac{3}{1+\phi} - \dfrac{x\phi_x}{(1+\phi)^2} \bigr) \,dx  }{} \lesssim \bar\varepsilon_0 + 1,
\end{equation*}
which contradicts \eqref{idi:ss=growingrate} for $ s $ large enough. 

Combining the above analysis, we have shown the following theorem.

\begin{thm}[Instability of self-similar expanding homogeneous solutions]
	Considering an initial perturbation $ (\phi_0,\phi_1) $ of the self-similarly expanding homogeneous solutions for the isentropic model of radiation gaseous stars, if it satisfies one of the following conditions,
	\begin{enumerate}
		\item the perturbation is homogeneous and $ (\phi_0,\phi_1 ) $ is not on the curve \eqref{eq:ss=002};
		\item $ E(\phi_0,\phi_1) \leq 0 $ and the perturbation is inhomogeneous,
	\end{enumerate}
	then for any $ 0 < \bar\varepsilon_0 < 0.5 $, there exists a $ \bar T > 0 $ such that $ \max \lbrace \norm{\phi(s)}{L_x^\infty}, \norm{x\phi_x(s)}{L_x^\infty} \rbrace > \bar \varepsilon_0 $ for $ t > \bar T $, no matter how small the initial perturbation is.
\end{thm}

\begin{rmk}
	At this moment, we don't know how to deal with the case when the initial perturbation $ (\phi_0,\phi_1) $ is inhomogeneous and $ E(\phi_0,\phi_1) > 0 $, in which case there might be various asymptotic behaviors. However, our analysis has already demonstrate that the self-similarly expanding solution of the isentropic model for radiation gaseous stars is unstable, even the perturbation is on the manifold of vanishing energy solutions, which behaviors differently to the result in \cite{Hadzic2016}. In fact, from the homogeneous perturbations investigated above, a perturbation of the self-similarly expanding solution might make the star configuration collapse or linearly expand as time grows up.
\end{rmk}


\section{Stability of the Linearly Expanding Homogeneous Solution for the Isentropic Model}\label{section:sb-linearly-isentropic}

To study the stability of the linearly expanding solution of the isentropic model for radiation gaseous stars, we denote the following point-wise bounds of the perturbation $ \vartheta $. 
\begin{equation}\label{le:supnorm}
	\tilde\omega := \sup_{\tau}\lbrace \norm{x\vartheta_x(\tau)}{L_x^\infty}, \norm{\vartheta(\tau)}{L_x^\infty}, \norm{x\vartheta_{x\tau}(\tau)}{L_x^\infty}, \norm{\vartheta_\tau(\tau)}{L_x^\infty} \rbrace,
\end{equation}
where $ \vartheta $ is the solution to the equation \eqref{eq:perturbationrbation=linearly} with the corresponding initial data and boundary conditions, $ (\vartheta(\cdot,0),\vartheta_\tau(\cdot,0)) = (\vartheta_0(\cdot),\vartheta_1(\cdot)) $ and \eqref{BC:isentropic,linearlyexpanding} respectively. 
Our goal is to show that for the perturbation small enough initially, the perturbation amplitude $ \tilde\omega $ will keep being small over time. In order to do so, we split the following analysis into three parts. 

First, by assuming that $ \tilde \omega $ is bounded by some small constant $ \varepsilon_0 $, we perform weighted estimates of the perturbation in Section \ref{sec:energy-est-isentropic}. The lower and higher order energy functionals are calculated here. Some temporal weighted estimates will be employed. Such estimates will capture the growth of the energy functionals. 

Next in Section  \ref{sec:interior-est-isentropic}, some interior estimates will be established. In particular, we explore the specific growth of the energy functionals which captures the regularity of the perturbation near the coordinate center $ x = 0 $. 

Finally in Section \ref{sec:pointwise-bounde-isentropic}, we collect all the estimates and show that the quantity $ \tilde\omega $ is bounded by the total energy functional defined in \eqref{isentropic-linearly-total-energy} and therefore bounded by the initial energy \eqref{isentropic-linearly-total-initial-energy}. To this end, we can employ continuous arguments to show the asymptotic stability of the linearly expanding homogeneous solution of the isentropic model.

In the following, it is also denoted the polynomial of $ \tilde \omega $ as $ \tilde \omega $. Also, the total energy and dissipation functionals are given by, for $ a \in (0,1) $,
\begin{align}
	& \mathcal{\tilde E}(t) : = ({\tilde\alpha} + {\tilde\alpha}^{1+a} ) \int x^4 \bar\rho \vartheta_\tau^2 \,dx + \int x^4 \bar\rho^{4/3} \vartheta_x^2 \,dx  + \int x^4 \vartheta_x^2 \,dx + \int x^4 \bar\rho \vartheta^2 \,dx + \tilde\alpha^{a-3}\int x^4 \bar\rho \vartheta_{\tau\tau}^2 \,dx {\nonumber} \\
	& ~~~~ + \tilde\alpha^{1+a} \int x^2 \lbrack (1+\vartheta) x\vartheta_{x\tau} - x\vartheta_x \vartheta_\tau\rbrack^2 \,dx + \tilde\alpha^{1+a} \int \chi(\vartheta_\tau^2 + x^2 \vartheta_{x\tau}^2) \,dx + \int \chi (\vartheta^2 + x^2 \vartheta_x^2)\,dx {\nonumber} \\
	& ~~~~ + \tilde\alpha^{a-5}\int\chi x^2 \bar\rho \vartheta_{\tau\tau}^2 \,dx + \int \mathcal G_x^2 \,dx + \tilde\alpha^{a-1} \int \mathcal{G}_{x\tau}^2 \,dx + \int \vartheta_x^2 \,dx + \int x^2 \vartheta_{xx}^2\,dx {\nonumber}\\
	& ~~~~ + \tilde\alpha^{a-1} \int \vartheta_{x\tau}^2\,dx + \tilde\alpha^{a-1} \int x^2 \vartheta_{xx\tau}^2\,dx ,  \label{isentropic-linearly-total-energy}   \\
	& \mathcal{\tilde D}(t) : = \int ( \tilde\alpha + \tilde\alpha^{1+a}) \int x^4 \bar\rho \vartheta_\tau^2 \,dx\,dt + \int ( \tilde\alpha^3 + \tilde\alpha^{3+a} ) \int x^2\lbrack (1+\vartheta)x\vartheta_{x\tau} - x \vartheta_x \vartheta_\tau \rbrack^2 \,dx \,dt {\nonumber}\\
	& ~~~~ + \int \tilde\alpha^{a-3} \int x^4 \bar\rho\vartheta_{\tau\tau}^2 \,dx\,d\tau + \int \tilde\alpha^{a-1}\int x^2 \lbrack(1+\vartheta)x\vartheta_{x\tau\tau} - x\vartheta_x\vartheta_{\tau\tau}\rbrack^2 \,dx\,d\tau {\nonumber}\\
	& ~~~~ + \int \tilde\alpha^{1+a} \int \chi ( \vartheta_\tau^2 + x^2 \vartheta_{x\tau}^2 ) \,dx + \int \tilde\alpha^{a-5}\int \chi x^2 \bar\rho\vartheta_{\tau\tau}^2 \,dx \,d\tau + \int \tilde\alpha^{a-3}\int\chi(\vartheta_{\tau\tau}^2 + x^2 \vartheta_{x\tau\tau}^2  ) \,dx\,d\tau {\nonumber}\\
	& ~~~~ + \int \tilde\alpha^{-3} \int \bar\rho^{4/3} \mathcal{G}_x\,dx \,d\tau. \label{isentropic-linearly-total-dissipation}
\end{align}
The initial energy is given by
\begin{equation}\label{isentropic-linearly-total-initial-energy}
	\begin{aligned}
	& \mathcal{\tilde E}_0 = \mathcal{\tilde E}_0(\vartheta_0,\vartheta_1) : = \int x^4 \bar\rho \vartheta_1^2 \,dx + \int x^4 \bar\rho \vartheta_0^2 \,dx + \int x^4 \vartheta_{0,x}^2 \,dx + \int x^4 \bar\rho^{4/3} \vartheta_{0,x}^2 \,dx \\
	& ~~~~ + \int x^4 \bar\rho\vartheta_{2}^2 \,dx + \int \chi ( \vartheta_0^2 + x^2 \vartheta_{0,x}^2 ) \,dx + \int \chi x^2 \bar\rho \vartheta_{2}^2 \,dx + \int ( \vartheta_{0,x}^2 + x^2 \vartheta_{0,xx}^2 )\,dx,
	 \end{aligned}
\end{equation}
where $ \vartheta_2 $ is the initial value of $ \vartheta_{\tau\tau} $ defined by \eqref{eq:perturbationrbation=linearly}. That is, $ \vartheta_2 $ is defined by the following identity
\begin{equation*}
\begin{aligned}
	& \dfrac{1}{(1+\vartheta_0)^2}\bigl( a_0 x\bar\rho \vartheta_{2} + a_0a_1 x\bar\rho \vartheta_1 \bigr) + \bigl( \dfrac{1}{1+\vartheta_0} - \dfrac{1}{(1+\vartheta_0)^4} \bigr)\delta x\bar\rho \\
		& ~~~~ + \bigl\lbrack \bigl(\dfrac{\bar\rho}{(1+\vartheta_0)^2(1+\vartheta_0+x\vartheta_{0,x})}\bigr)^{4/3}\bigr\rbrack_x - \dfrac{(\bar\rho^{4/3})_x}{(1+\vartheta_0)^4} = \dfrac{4\mu}{3}a_0^3 \bigl( \dfrac{\vartheta_1+x\vartheta_{1,x}}{1+\vartheta_0+x\vartheta_{0,x}} +  \dfrac{2\vartheta_1}{1+\vartheta_0} \bigr)_x. 
\end{aligned}	
\end{equation*}

Here $ \mathcal G,\chi $ are the relative entropy and the cut-off function defined as 
\begin{align}
\mathcal{G} & := \mathfrak H(\vartheta) = \log{(1+\vartheta)^2(1+\vartheta+x\vartheta_x)} {\label{isen:RelativeEntropy}} ,\\
\chi(x) & :=  \begin{cases}
		1 & 0 \leq x \leq  R_0/2, \\
		0 & 3 R_0 / 4 \leq x \leq R_0,
	\end{cases} {\label{interiorCutOffFunction}} 
\end{align}
and $ -4 \leq \chi'(x) \leq 0 $. Here $ \mathfrak H $ is the relative entropy functional defined in \eqref{relative-entropy-functional}.

\subsection{Energy Estimates}\label{sec:energy-est-isentropic}

We start our estimates with the basic energy estimate. The following lemma is a direct consequence of the $ L^2 $-estimate of \eqref{eq:perturbationrbation=linearly}. 

\begin{lm}\label{lm:isentropic-linearly-basic-energy}
	Considering a smooth solution $ \vartheta $ to \eqref{eq:perturbationrbation=linearly} with the corresponding boundary condition \eqref{BC:isentropic,linearlyexpanding}, if the constant $ \delta $ and $ a_0, a_1 > 0 $ for $ \tilde \alpha = \tilde\alpha(\tau) > 0 $ defined by \eqref{def:growing-rate-isentropic}  
	satisfy
	\begin{equation}\label{constraint-002}
	\delta > - \dfrac{a_0a_1^2}{8},
	\end{equation}
there exists a constant $ 0 < \varepsilon_0 <1 $ such that $ \tilde \omega < \varepsilon_0 $ will imply the  following basic energy estimate,
\begin{equation}\label{lm:001}
\begin{aligned}
& \tilde\alpha\int x^4 \bar\rho \vartheta_\tau^2 \,dx + \int x^4 \bar\rho^{4/3} \vartheta_x^2 \,dx + \int \tilde\alpha \int x^4 \bar\rho \vartheta_\tau^2 \,dx \,d\tau \\
& ~~~~~~ + \int \tilde\alpha^3 \int x^2 \lbrack (1+\vartheta) x \vartheta_{x\tau} - x\vartheta_x\vartheta_\tau \rbrack^2 \,dx\,d\tau \lesssim \mathcal{\tilde E}_0.
\end{aligned}
\end{equation}
Moreover, such a basic energy estimate yields the bound
\begin{equation}\label{estimates-002}
\begin{aligned}
\int x^4 \vartheta_x^2 \,dx + \int x^4 \bar\rho \vartheta^2 \,dx \lesssim \mathcal{\tilde E}_0.
\end{aligned}
\end{equation}
Additionally, for any positive constant $ 0 < a < 1 $, which will be fixed in following paragraph, a faster decay estimate holds. That is, for $ \varepsilon_0 $ small enough, 
\begin{equation}\label{lm:002}
\begin{aligned}
& \tilde\alpha^{1+a} \int x^4 \bar\rho\vartheta_\tau^2 \,dx + \int \tilde\alpha^{1+a}\int x^4 \bar\rho\vartheta_\tau^2 \,dx\,d\tau \\
& ~~~~~~ + \int \tilde\alpha^{3+a}\int x^2 \lbrack (1+\vartheta)x\vartheta_{x\tau} - x\vartheta_{x}\vartheta_{\tau} \rbrack^2 \,dx\,d\tau \lesssim \mathcal{\tilde E}_0.
\end{aligned}
\end{equation}
\end{lm}

\begin{pf}
Multiply \eqref{eq:perturbationrbation=linearly} with $ x^3 (1+\vartheta)^2 \vartheta_\tau $ and integrate the resulting equation in the spatial variable. After applying integration by parts, it follows,
\begin{equation}\label{ene:le=001}
		\dfrac{d}{dt} \tilde E_1 + \tilde D_1 = 0,
\end{equation}
where
\begin{equation*}
\begin{aligned}
	& \tilde E_1 : = \tilde E_{11} + \tilde E_{12} + \tilde E_{13} = \dfrac{\tilde\alpha}{2} \int x^4 \bar\rho \vartheta_\tau^2 \,dx + \delta \int x^4 \bar\rho \biggl( \dfrac{(1+\vartheta)^2}{2} + \dfrac{1}{1+\vartheta}- \dfrac{3}{2}\biggr) \,dx  \\
	& ~~~~ + \int x^2 \bar\rho^{4/3} \biggl\lbrace \dfrac{3}{((1+\vartheta)^2(1+\vartheta+x\vartheta_x))^{1/3}} - \dfrac{3}{1+\vartheta} + \dfrac{x\vartheta_x}{(1+\vartheta)^2} \biggr\rbrace\,dx, \\
	& \tilde D_1 : = \tilde D_{11} + \tilde D_{12}  = \dfrac{\tilde\alpha_\tau}{2} \int x^4 \bar\rho\vartheta_\tau^2\,dx + \dfrac{4\mu}{3} \tilde\alpha^3 \int \dfrac{x^2 \lbrack(1+\vartheta)x\vartheta_{x\tau} - x\vartheta_x\vartheta_\tau\rbrack^2}{1+\vartheta + x \vartheta_x}\,dx.
\end{aligned}	
\end{equation*}
Multiply \eqref{ene:le=001} with $ \tilde\alpha^a $, $ a> 0 $. It follows
\begin{equation}\label{ene:le=002}
\begin{aligned}
	& \dfrac{d}{dt} \dfrac{\tilde\alpha^{1+a}}{2} \int x^4 \bar\rho \vartheta_\tau^2 \,dx + \dfrac{(1-a)\tilde\alpha^{a}\tilde\alpha_\tau}{2} \int x^4 \bar\rho \vartheta_\tau^2\,dx + \tilde\alpha^{a} \tilde D_{12} = - \tilde\alpha^{a}\dfrac{d}{dt} \bigl( \tilde E_{12} + \tilde E_{13} \bigr).
\end{aligned}	
\end{equation}
Integrating \eqref{ene:le=001} and \eqref{ene:le=002} over the temporal variable yields,
\begin{align}
	& \tilde E_1  + \int \tilde D_1 \, d\tau \lesssim \mathcal{\tilde E}_0, {\label{ine:001}} \\
	& \dfrac{\tilde \alpha^{1+a}}{2} \int x^4\bar\rho\vartheta_\tau^2 \,dx +  \int \dfrac{(1-a)\tilde\alpha^a\tilde\alpha_\tau}{2} \int x^4 \bar\rho \vartheta_\tau^2\,dx \,d\tau + \int\tilde\alpha^{a}\tilde D_{12} \,d\tau {\nonumber}\\
	& ~~~~ \lesssim (1+|\delta|) (1+\tilde\omega) \int \tilde\alpha^{a}\biggl( \int x^4 \bar\rho \vartheta^2 \,dx + \int x^4 \vartheta_{x}^2 \,dx \biggr)^{1/2} {\nonumber}\\ & ~~~~~~~\times \biggl( \int x^4 \bar\rho \vartheta_\tau^2 \,dx + \int x^4 \vartheta_{x\tau}^2 \,dx \biggr)^{1/2} \,d\tau + \mathcal{\tilde E}_0 . {\label{ine:002}}  
\end{align}

Here we have applied the Hardy's inequality and the fact $ \bar\rho^{1/3}(x) \simeq (R_0-x) $ to get the following inequalities,
\begin{gather*}
	\int x^2 \vartheta^2 \,dx \lesssim \int x^4 \bar\rho^{1/3} \vartheta^2 \,dx + \int x^4 \vartheta_x^2 \,dx 
	 \lesssim \int x^4 \bar\rho \vartheta^2 \,dx + \int x^4 \vartheta_x^2 \,dx, \\
	\int x^2 \vartheta_\tau^2 \,dx \lesssim \int x^4 \bar\rho\vartheta_\tau^2 \,dx + \int x^4\vartheta_{x\tau}^2 \,dx.
\end{gather*}
We will use this type of estimates in the following without mentioning it again. The following analysis of the energy and dissipation functionals defined above holds after applying the H\"older inequality and the Hardy's inequality. 
\todo{13 Sep 2017}
\begin{align*}
	\tilde E_{11} & \gtrsim \tilde \alpha \int x^4 \bar\rho \vartheta_\tau^2\,dx  \gtrsim a_0 e^{\beta_1\tau} \int x^4 \bar\rho \vartheta_\tau^2 \,dx , ~~~~~~ 
	\tilde E_{12} \geq \min \biggl\lbrace (3/2+\tilde \omega) \delta \int x^4 \bar\rho \vartheta^2 \,dx, 0 \biggr\rbrace, \\
	\tilde E_{13} & \gtrsim \dfrac{2}{3} \int x^4 \bar\rho^{4/3} \vartheta_x^2 \,dx - \tilde \omega \int x^2 \bar\rho^{4/3} (\vartheta^2 + x^2 \vartheta_{x}^2) \,dx {\nonumber} \\
	& \gtrsim (1-\tilde\omega)\int x^4 \bar\rho^{4/3} \vartheta_x^2 \,dx - \tilde \omega \int x^4 \bar\rho \vartheta^2 \,dx ,\\
	\tilde D_{11} & = \dfrac{\tilde\alpha_\tau}{2} \int x^4\bar\rho \vartheta_\tau^2\,dx \geq \dfrac{\beta_1}{2} \tilde \alpha \int x^4 \bar\rho \vartheta_\tau^2 \,dx \geq \dfrac{a_0\beta_1}{2} e^{\beta_1\tau} \int x^4 \bar\rho \vartheta_\tau^2\,dx , \\
	\tilde D_{12} & \gtrsim (1-\tilde\omega) \tilde\alpha^3 \int x^2 \biggl\lbrack \bigl( \dfrac{x\vartheta_x}{1+\vartheta} \bigr)_\tau \biggr\rbrack^2 \,dx {\nonumber}  = (1-\tilde\omega)\tilde\alpha^3 \int x^4 \biggl\lbrack \bigl(\dfrac{\vartheta_\tau}{1+\vartheta}\bigr)_x\biggr\rbrack^2\,dx {\nonumber}\\ 
	& \gtrsim (1-\tilde\omega) \tilde\alpha^3 \int x^2 \lbrack(1+\vartheta)x\vartheta_{x\tau} - x\vartheta_x\vartheta_\tau\rbrack^2 \,dx.
\end{align*}
In addition, 
\begin{align*}
	& \bigl(\int x^4 \bar\rho \vartheta^2 \,dx \bigr)^{1/2} \dfrac{d}{d\tau} \bigl(\int x^4 \bar\rho \vartheta^2 \,dx \bigr)^{1/2} = \dfrac{1}{2} \dfrac{d}{d\tau} \int x^4\bar\rho \vartheta^2 \,dx {\nonumber}\\
	& = \int x^4 \bar\rho \vartheta\vartheta_\tau\,dx \leq \bigl(\int x^4 \bar\rho \vartheta^2 \,dx \bigr)^{1/2} \bigl(\int x^4 \bar\rho \vartheta_\tau^2 \,dx \bigr)^{1/2}, \\
	& \bigl( \int x^2 \bigl(\dfrac{x\vartheta_x}{1+\vartheta} \bigr)^2\,dx \bigr)^{1/2} \dfrac{d}{d\tau}\bigl( \int x^2 \bigl(\dfrac{x\vartheta_x}{1+\vartheta} \bigr)^2\,dx \bigr)^{1/2} = \dfrac{1}{2}\dfrac{d}{d\tau} \int x^2 \bigl(\dfrac{x\vartheta_x}{1+\vartheta} \bigr)^2\,dx {\nonumber}\\
	& \leq \bigl( \int x^2 \bigl(\dfrac{x\vartheta_x}{1+\vartheta} \bigr)^2\,dx \bigr)^{1/2}  \bigr(\int x^2 \bigl\lbrack\bigl(\dfrac{x\vartheta_x}{1+\vartheta} \bigr)_\tau\bigr\rbrack^2\,dx \bigr)^{1/2},
\end{align*}
from which one can derive
\begin{align*}
	& \dfrac{d}{d\tau} \bigl(\int x^4 \bar\rho \vartheta^2 \,dx \bigr)^{1/2} \leq \bigl( \int x^4 \bar\rho \vartheta_\tau^2\,dx \bigr)^{1/2}, \\
	& \dfrac{d}{d\tau} \bigl( \int x^2 \bigl( \dfrac{x\vartheta_x}{1+\vartheta}\bigr)^2\,dx \bigr)^{1/2} \leq \bigr(\int x^2 \bigl\lbrack\bigl(\dfrac{x\vartheta_x}{1+\vartheta} \bigr)_\tau\bigr\rbrack^2\,dx \bigr)^{1/2},
\end{align*}
and therefore,
\begin{align}
& \bigl( \int x^4 \bar\rho \vartheta^2 \,dx \bigr)^{1/2} \leq \mathcal{\tilde E}_0^{1/2} + \int \bigl( \int x^4 \bar\rho \vartheta_\tau^2 \,dx \bigr)^{1/2} \,d\tau, {\label{ine:003}}  \\
& \bigl( \int x^4 \vartheta_x^2 \,dx \bigr)^{1/2} \leq (1+\tilde\omega) \bigl( \int x^2 \bigl(\dfrac{x\vartheta_x}{1+\vartheta}\bigr)^2\,dx \bigr)^{1/2} \leq \mathcal{\tilde E}_0^{1/2} {\nonumber}\\
& ~~~~ + \int \bigr(\int x^2 \bigl\lbrack\bigl(\dfrac{x\vartheta_x}{1+\vartheta} \bigr)_\tau\bigr\rbrack^2\,dx \bigr)^{1/2} \,d\tau. {\label{ine:006}}
\end{align}

Then \eqref{ine:001} implies,
\begin{equation}\label{estimates-003}
	\begin{aligned}
		& \tilde\alpha \int x^4 \bar\rho \vartheta_\tau^2 \,dx + (1-\tilde\omega)\int x^4 \bar\rho^{4/3} \vartheta_x^2 \,dx \\
		& ~~~~~~ + \biggl( 1 + \dfrac{(6+\tilde\omega)\min\lbrace \delta, 0 \rbrace - \tilde\omega}{a_0\beta_1^2} \biggr)\ \int  \dfrac{\tilde\alpha_\tau}{2} \int x^4\bar\rho\vartheta_\tau^2\,dx \,d\tau  \\
		& ~~~~~~ + (1-\tilde\omega) \int \tilde\alpha^3 \int x^2 \bigl\lbrack \bigl( \dfrac{x\vartheta_x}{1+\vartheta} \bigr)_\tau  \bigr\rbrack^2\,dx \,d\tau  \\
		& ~~~~ \lesssim (1 + |\delta|)(1 + \tilde\omega) \mathcal{\tilde E}_0 \lesssim \mathcal{\tilde E}_0 ,
	\end{aligned}
\end{equation}
where it has been use the fact from \eqref{ine:003},  
\begin{equation}\label{ine:007}
\begin{aligned}
	&  \int x^4 \bar\rho\vartheta^2 \,dx \leq \mathcal{\tilde E}_0 + 2 \bigl( \int \bigl( \int x^4 \bar\rho\vartheta_\tau^2 \,dx \bigr)^{1/2} \,d\tau\bigr)^2 \\
	& ~~~~ \leq \mathcal{\tilde E}_0 + 4 \bigl( \int \tilde\alpha_\tau^{-1}\,d\tau\bigr) \bigl( \int \dfrac{\tilde\alpha_\tau}{2} \int x^4 \bar\rho\vartheta_\tau^2 \,dx\,d\tau  \bigr) \\
	& ~~~~ \leq \mathcal{\tilde E}_0 + \dfrac{4}{a_0\beta_1^2} \cdot \int \dfrac{\tilde\alpha_\tau}{2}  \int x^4 \bar\rho\vartheta_\tau^2 \,dx\,d\tau .
\end{aligned}
\end{equation}
Notice, we need the following constraint on $ |\delta| $, $ a_0 $ and $ \beta_1 $ to ensure \eqref{estimates-003} to give an appropriate bound, recalling \eqref{le:growingspeed}
\begin{equation*}
1 + 6 \dfrac{\min\lbrace \delta, 0 \rbrace}{a_0 \beta_1^2} > 0 ~~ \text{or equivalently} ~ \delta > - \dfrac{a_0\beta_1^2}{6}.
\end{equation*}
We only have to consider the case when $ \delta < 0 $. Then $ \beta_1^2 = a_1^2 + \dfrac{2\delta}{a_0} $. 
Therefore, we impose the condition
\begin{equation*}\tag{\ref{constraint-002}}
	\delta > - \dfrac{a_0a_1^2}{8},
\end{equation*}
and \eqref{lm:001} follows from \eqref{estimates-003}.
Additionally, from  \eqref{lm:001}, \eqref{ine:006} and \eqref{ine:007},
\begin{equation*}\tag{\ref{estimates-002}}
\begin{aligned}
	& \int x^4 \vartheta_x^2 \,dx 
	\lesssim \mathcal{\tilde E}_0 + \bigl(\int \tilde\alpha^{-3} \,d\tau \bigr) \bigl( \int \tilde\alpha^3 \int x^2 \bigl\lbrack\bigl(\dfrac{x\vartheta_x}{1+\vartheta} \bigr)_\tau\bigr\rbrack^2\,dx\,d\tau \bigr)\\
	& ~~~~ \lesssim \mathcal{\tilde E}_0 + \dfrac{1}{a_0^3 \beta_1} \int \tilde\alpha^3 \int x^2 \bigl\lbrack\bigl(\dfrac{x\vartheta_x}{1+\vartheta} \bigr)_\tau\bigr\rbrack^2\,dx\,d\tau \lesssim \mathcal{\tilde E}_0, \\
	& \int x^4 \bar\rho \vartheta^2 \,dx \lesssim \mathcal{\tilde E}_0.
\end{aligned}
\end{equation*}

Similarly, \eqref{ine:002} implies,
\begin{equation}\label{ine:005}
	\begin{aligned}
		& \tilde\alpha^{1+a} \int x^4 \bar\rho \vartheta_\tau^2\,dx + (1-a) \int \tilde\alpha^a\tilde\alpha_\tau \int x^4 \bar\rho \vartheta_\tau^2\,dx \,d\tau\\
		& ~~~~~~ + (1-\tilde\omega) \int \tilde\alpha^{3+a}  \int x^2 \bigl\lbrack \bigl(\dfrac{x\vartheta_x}{1+\vartheta}\bigr)_\tau\bigr\rbrack^2 \,dx \,d\tau \\
		& ~~ \lesssim (1+|\delta|)(1+\tilde\omega)  \int \tilde\alpha^a \biggl( \int x^4\bar\rho \vartheta_\tau^2\,dx + \int x^2 \bigl\lbrack \bigl(\dfrac{x\vartheta_x}{1+\vartheta}\bigr)_\tau\bigr\rbrack^2 \,dx \biggr)^{1/2} \,d\tau \\ & ~~~~~~ \times {\mathcal{\tilde E}_0}^{1/2}  + \mathcal{\tilde E}_0 \\
		& ~~ \lesssim (1+|\delta|)(1+\tilde\omega) \biggl\lbrace  \bigl(  \int \tilde\alpha^{a} \tilde\alpha_\tau \int x^4\bar\rho \vartheta_\tau^2\,dx \,d\tau\bigr)^{1/2}\bigl( \int \tilde\alpha^a\tilde\alpha_\tau^{-1} \,d\tau \bigr)^{1/2} \\
		& ~~~~ + \biggl( \int \tilde\alpha^{3+a} \int x^2 \bigl\lbrack \bigl(\dfrac{x\vartheta_x}{1+\vartheta}\bigr)_\tau\bigr\rbrack^2 \,dx \,d\tau \biggr)^{1/2}\bigl(  \int \tilde\alpha^{a-3} \,d\tau \bigr)^{1/2}  \biggr\rbrace \times {\mathcal{\tilde E}_0}^{1/2}  + \mathcal{\tilde E}_0,
		\end{aligned}
\end{equation}
where the following inequality is made used of on the right hand side together with \eqref{estimates-002},
\begin{align*}
	& \int x^4 \vartheta_{x\tau}^2 \,dx = \int x^2 \dfrac{\lbrack(1+\vartheta)x\vartheta_{x\tau}\rbrack^2}{(1+\vartheta)^2}\\
	& ~~~~ \lesssim \int x^2 \dfrac{\lbrack (1+\vartheta)x\vartheta_{x\tau} - x\vartheta_x \vartheta_\tau \rbrack^2}{(1+\vartheta)^2}\,dx + \int x^2 \dfrac{x^2\vartheta_x^2 \vartheta_\tau^2}{(1+\vartheta)^2}\,dx\\
	& ~~~~ \lesssim (1+\tilde\omega)\int x^2 \bigl\lbrack\bigl(\dfrac{x\vartheta_x}{1+\vartheta} \bigr)_\tau\bigr\rbrack^2\,dx + \tilde\omega \int x^4 \bigl( \bar\rho \vartheta_\tau^2 + \bigl\lbrack\bigl(\dfrac{\vartheta_\tau}{1+\vartheta}\bigr)_x \bigr\rbrack^2 \bigr) \,dx \\
	& ~~~~ \lesssim (1+\tilde\omega) \int x^2 \bigl\lbrack\bigl(\dfrac{x\vartheta_x}{1+\vartheta} \bigr)_\tau\bigr\rbrack^2\,dx + \tilde \omega \int x^4 \bar\rho \vartheta_\tau^2 \,dx.
\end{align*}
By applying the Cauchy's inequality on the right of \eqref{ine:005} together with \eqref{idi:le=growingrate}, it concludes
\begin{equation}\label{estimates-001}
	\begin{aligned}
		& \tilde\alpha^{1+a} \int x^4 \bar\rho \vartheta_\tau^2\,dx + (1- a -\tilde\omega) \int \tilde \alpha^{1+a} \int x^4 \bar\rho \vartheta_\tau^2\,dx \,d\tau\\
		& ~~~~~~ + (1-\tilde\omega)  \int \tilde\alpha^{3+a} \int x^2 \bigl\lbrack \bigl(\dfrac{x\vartheta_x}{1+\vartheta}\bigr)_\tau\bigr\rbrack^2 \,dx \,d\tau \lesssim \mathcal{\tilde E}_0,
	\end{aligned}
\end{equation}
provided
\begin{equation*}
	\int \tilde\alpha^a\tilde\alpha_{\tau}^{-1} \,d\tau \simeq \int \tilde\alpha^{a-1} \,d\tau < \infty, ~ \int \tilde\alpha^{a-3}\,d\tau < \infty,
\end{equation*}
or equivalently,
\begin{equation}\label{constraint-001}
	0 < a < 1.
\end{equation}

\end{pf}


Next, we shall study the energy estimates on the temporal derivative of \eqref{eq:perturbationrbation=linearly}. 
Multiply \eqref{eq:perturbationrbation=linearly} with $ \tilde\alpha^{-3}( 1 + \vartheta )^2 $ and apply the temporal derivative to the resulting equation. It reads,
\begin{equation}\label{eq:temporalderivative=linearlyexpanding}
	\begin{aligned}
		& \tilde \alpha^{-2} x \bar\rho \vartheta_{\tau\tau\tau} - \tilde\alpha^{-3} \tilde \alpha_\tau x\bar\rho \vartheta_{\tau\tau} + ( \tilde\alpha^{-3} \tilde\alpha_{\tau\tau} - 3 \tilde\alpha^{-4}\tilde\alpha_\tau^2 ) x \bar\rho \vartheta_{\tau} \\
		& ~~ + \biggl\lbrace  \tilde\alpha^{-3} \bigl( \vartheta_\tau + 2\dfrac{\vartheta_{\tau}}{(1+\vartheta)^3} \bigr) - 3 \tilde\alpha^{-4} \tilde\alpha_\tau \bigl( (1+\vartheta) - \dfrac{1}{(1+\vartheta)^2}\bigr) \biggr\rbrace \delta x\bar\rho \\
		& ~~ + \biggl\lbrace \tilde\alpha^{-3} (1+\vartheta)^2 \biggl\lbrace \bigl\lbrack \bigl(\dfrac{\bar\rho}{(1+\vartheta)^2(1+\vartheta+x\vartheta_x)}\bigr)^{4/3}\bigr\rbrack_x - \dfrac{(\bar\rho^{4/3})_x}{(1+\vartheta)^4} \biggr\rbrace  \biggr\rbrace_\tau \\
		& = (1+\vartheta)^2\bigl( \mathfrak{\tilde B}_{x\tau} + 4\mu \bigl(\dfrac{\vartheta_\tau}{1+\vartheta}\bigr)_{x\tau} \bigr) + 2 (1+\vartheta)\vartheta_\tau\bigl(\mathfrak{\tilde B}_x + 4\mu \bigl(\dfrac{\vartheta_\tau}{1+\vartheta}\bigr)_x \bigr).
	\end{aligned}
\end{equation}

We will perform corresponding energy estimates on \eqref{eq:temporalderivative=linearlyexpanding} in the next lemma.
\begin{lm}\label{lm:isentropic-linearly-temporal-derivative}
Under the same assumption as in Lemma \ref{lm:isentropic-linearly-basic-energy}, for
 $ \varepsilon_0 $ small enough, we have the following higher order energy estimate
\begin{equation}\label{lm:101}
	\begin{aligned}
	& \tilde\alpha^{a-3} \int x^4 \bar\rho \vartheta_{\tau\tau}^2 \,dx + \int \tilde\alpha^{a-3} \int x^4 \bar\rho \vartheta_{\tau\tau}^2 \,dx \,d\tau \\
	& ~~~~~~ + \int \tilde\alpha^{a-1} \int x^2\lbrack (1+\vartheta)x\vartheta_{x\tau\tau} - x\vartheta_{x}\vartheta_{\tau\tau} \rbrack^2\,dx\,d\tau \lesssim \mathcal{\tilde E}_0.
	\end{aligned}
	\end{equation}
\end{lm}

\begin{pf}
Multiply the above with $ \tilde\alpha^{b}x^3 \vartheta_{\tau\tau} $ and integrate the resulting equation in the spatial variable where $ b < 0 $. After integration by parts, it follows, 
\begin{equation}\label{ene:le=003}
	\dfrac{d}{d\tau} \tilde E_2 + \tilde D_2 = \tilde L_1 + \tilde L_2 + \tilde L_3 + \tilde L_4 + \tilde L_5, 
\end{equation}
where
\begin{align*}
	& \tilde E_2 
	:= \dfrac{\tilde\alpha^{b-2}}{2} \int x^4 \bar\rho \vartheta_{\tau\tau}^2 \,dx , \\
	& \tilde D_2 
	:= - \dfrac{b}{2} \tilde\alpha^{b-3} \tilde \alpha_\tau \int x^4 \bar\rho \vartheta_{\tau\tau}^2 \,dx {\nonumber}\\ & ~~~~~~ + \dfrac{4\mu}{3} \tilde \alpha^b \int \dfrac{x^2\lbrack (1+\vartheta)x\vartheta_{x\tau\tau} - x\vartheta_x\vartheta_{\tau\tau}\rbrack^2}{1+\vartheta+x\vartheta_x}\,dx , \\ 
	& \tilde L_1 : = - (\tilde\alpha^{b-3}\tilde\alpha_{\tau\tau} - 3 \tilde\alpha^{b-4}\tilde\alpha_\tau^2)\int x^4 \bar\rho \vartheta_\tau\vartheta_{\tau\tau}\,dx , \\
	& \tilde L_2 : = - \int \biggl\lbrace  \tilde\alpha^{b-3} \bigl( \vartheta_\tau + 2\dfrac{\vartheta_{\tau}}{(1+\vartheta)^3} \bigr) - 3 \tilde\alpha^{b-4} \tilde\alpha_\tau \bigl( (1+\vartheta) - \dfrac{1}{(1+\vartheta)^2}\bigr) \biggr\rbrace {\nonumber} \\
	& ~~~~~~ \times \delta x^4\bar\rho \vartheta_{\tau\tau} \,dx , \\
	& \tilde L_3 : = - \dfrac{4}{3}  \mu \tilde\alpha^b \int x^2 \biggl\lbrace -  \dfrac{(1+\vartheta)^2(\vartheta_\tau+x\vartheta_{x\tau})^2(\vartheta_{\tau\tau} + x\vartheta_{x\tau\tau})}{(1+\vartheta+x\vartheta_x)^2} {\nonumber} \\
	& ~~~~~~ + 2 \dfrac{(1+\vartheta)\vartheta_\tau(\vartheta_\tau+x\vartheta_{x\tau})(\vartheta_{\tau\tau}+x\vartheta_{x\tau\tau})}{1+\vartheta+x\vartheta_x} 
	- \vartheta_\tau^2(\vartheta_{\tau\tau} + x\vartheta_{x\tau\tau}) \biggr\rbrace \,dx , \\
	& \tilde L_4 : = -3 \tilde\alpha^{b-4}\tilde\alpha_\tau\int x^2 \bar\rho^{4/3} \biggl\lbrace \dfrac{2 x\vartheta_x \vartheta_{\tau\tau} }{(1+\vartheta)^{5/3}(1+\vartheta+x\vartheta_x)^{4/3}} + 2 \dfrac{x\vartheta_x \vartheta_{\tau\tau}}{(1+\vartheta)^3} {\nonumber}\\ 
	& ~~~~~~ + \dfrac{3 \vartheta_{\tau\tau} + x \vartheta_{x\tau\tau}}{(1+\vartheta)^{2/3}(1+\vartheta+x\vartheta_x)^{4/3}} - \dfrac{3 \vartheta_{\tau\tau} + x \vartheta_{x\tau\tau}}{(1+\vartheta)^2} \biggr\rbrace\,dx , \\
	& \tilde L_5 : = \tilde\alpha^{b-3} \int x^2 \bar\rho^{4/3} \biggl\lbrace (3\vartheta_{\tau\tau} + x \vartheta_{x\tau\tau}) \bigl\lbrack \dfrac{1}{(1+\vartheta)^{2/3}(1+\vartheta+x\vartheta_x)^{4/3}} - \dfrac{1}{(1+\vartheta)^2} \bigr\rbrack_\tau  {\nonumber}\\ 
	& ~~~~~~ + \theta_{\tau\tau} \bigl\lbrack \dfrac{2x\vartheta_x}{(1+\vartheta)^{5/3}(1+\vartheta+x\vartheta_x)^{4/3}} + 2 \dfrac{x\vartheta_x}{(1+\vartheta)^3} \bigr\rbrack_\tau \biggr\rbrace \,dx .
\end{align*}

\noindent Integrating \eqref{ene:le=003} over the temporal variable yields,
\begin{equation}\label{ine:101}
 \tilde E_2 + \int \tilde D_2 \,d\tau = \sum_{k=1}^5 \int \tilde L_k\,d\tau.	
\end{equation}
The left of \eqref{ine:101} admits the following bound
\begin{align*}
	& \tilde E_2 + \int \tilde D_2 \,d\tau \gtrsim \tilde\alpha^{b-2} \int x^4 \bar\rho 
	\vartheta_{\tau\tau}^2 \,dx  + |b| \beta_1 \int \tilde\alpha^{b-2} \int x^4\bar\rho \vartheta_{\tau\tau}^2 \,dx \, d\tau \\
	& ~~ + (1-\tilde\omega) \int \tilde\alpha^b \int x^2 ((1+\vartheta)x\vartheta_{x\tau\tau} - x\vartheta_x \vartheta_{\tau\tau})^2\,dx \,d\tau.
\end{align*}
In the meantime, the right of \eqref{ine:101} can be estimated as in the following. 
\begin{align*}
	& \int \tilde L_1 \,d\tau \lesssim \sup_\tau \tilde\alpha^{(b-a-3)/2}\cdot \bigl( \int \tilde a ^{b-2} \int x^4 \bar\rho \vartheta_{\tau\tau}^2 \,dx \,d\tau \bigr)^{1/2} \bigl( \int \tilde \alpha^{1+a} \int x^4 \bar\rho \vartheta_\tau^2 \,dx \,d\tau \bigr)^{1/2} {\nonumber} \\
	& ~~~~ \lesssim \epsilon \int \tilde\alpha^{b-2} \int x^4 \bar\rho \vartheta_{\tau\tau}^2 \,dx \,d\tau + C_\epsilon \sup_\tau \tilde\alpha^{b-a-3} \cdot \int \tilde\alpha^{1+a} \int x^4 \bar\rho \vartheta_\tau^2 \,dx \,d\tau, \\
	& \int \tilde L_2 \,d\tau 
	\lesssim \epsilon \int \tilde\alpha^{b-2} \int x^4 \bar\rho \vartheta_{\tau\tau}^2 \,dx \,d\tau  + C_\epsilon \tilde\alpha^{b-a-5} \int \tilde\alpha^{1+a} \int x^4 \bar\rho \vartheta_\tau^2\,dx \,d\tau {\nonumber}\\
	& ~~~~ + C_\epsilon \int \tilde\alpha^{b-4} \,d\tau \cdot \sup_\tau \int x^4 \bar\rho \vartheta^2\,dx, \\
	& \int \tilde L_3 \,d\tau 
	\lesssim \tilde \omega \int \tilde\alpha^{b} \bigl( \int x^2 (\vartheta_\tau^2 + x^2 \vartheta_{x\tau}^2) \,dx \bigr)^{1/2} \bigl( \int x^2 (\vartheta_{\tau\tau}^2 + x^2 \vartheta_{x\tau\tau}^2) \,dx \bigr)^{1/2} \,d\tau {\nonumber} \\
	& ~~ \lesssim \tilde\omega \biggl\lbrace \sup_\tau \tilde\alpha^{b-a+1} \cdot \int \tilde\alpha^{1+a} \int x^4 \bar\rho \vartheta_\tau^2 \,dx \,d\tau + \int \tilde\alpha^{b-2} \int x^4\bar\rho\vartheta_{\tau\tau}^2\,dx \,d\tau {\nonumber}\\
	& ~~~~ + \sup_\tau \tilde\alpha^{b-a-1} \cdot \int \tilde\alpha^{1+a}\int x^4 \bar\rho\vartheta_\tau^2 \,dx \,d\tau + \int \tilde\alpha^{b}  \int x^2 ((1+\vartheta)x\vartheta_{x\tau\tau} - x\vartheta_x \vartheta_{\tau\tau})^2 \,dx \,d\tau {\nonumber}\\
	& ~~~~ + \sup_\tau \tilde\alpha^{b-a-1}\cdot \int \tilde\alpha^{3+a} \int x^2 ((1+\vartheta)x\vartheta_{x\tau}-x\vartheta_x\vartheta_\tau)^2 \,dx \,d\tau + \int \tilde\alpha^{b-2}  \int x^4 \bar\rho \vartheta_{\tau\tau}^2 \,dx\,d\tau {\nonumber}\\
	& ~~~~ + \sup_\tau \tilde\alpha^{b-a-3} \cdot \int \tilde\alpha^{3+a} \int x^2 ((1+\vartheta)x\vartheta_{x\tau} - x\vartheta_x\vartheta_\tau)^2 \,dx \, d\tau {\nonumber}\\
	& ~~~~ + \int \tilde\alpha^{b} \int x^2 ((1+\vartheta)x\vartheta_{x\tau\tau} - x\vartheta_x\vartheta_{\tau\tau})^2 \,dx \,d\tau\biggr\rbrace, \\
	& \int \tilde L_4\, d\tau \lesssim (1+\tilde\omega) \int \tilde\alpha^{b-3} \int x^2 \bar\rho^{4/3} (|\vartheta| + |x\vartheta_x|) (|\vartheta_{\tau\tau}| + |x\vartheta_{x\tau\tau}|) \,dx\,d\tau {\nonumber} \\
	& ~~~~ \lesssim \epsilon \biggl\lbrace \int \tilde\alpha^{b-2} \int x^4 \bar\rho\vartheta_{\tau\tau}^2 \,dx \,d\tau + \int \tilde\alpha^{b} \int x^2 ((1+\vartheta)x\vartheta_{x\tau\tau} - x\vartheta_x\vartheta_{\tau\tau})^2 \,dx \,d\tau \biggr\rbrace {\nonumber}\\
	& ~~~~  + C_\epsilon \sup_\tau \biggl(\int x^4 \bar\rho \vartheta^2 \,dx + \int x^4\vartheta_x^2 \,dx \biggr) \biggl\lbrace \int \tilde\alpha^{b-4} \,d\tau + \int \tilde\alpha^{b-6} \,d\tau \biggr\rbrace, \\
	& \int \tilde L_5 \,d\tau \lesssim (1+\tilde\omega) \int \tilde\alpha^{b-3} \int x^2 \bar\rho^{4/3} (|\vartheta_\tau| + |x\vartheta_{x\tau}|) (|\vartheta_{\tau\tau}| + |x\vartheta_{x\tau\tau}|) \,dx\,d\tau {\nonumber} \\
	& ~~~~ \lesssim \epsilon \biggl\lbrace \int \tilde\alpha^{b-2} \int x^4 \bar\rho\vartheta_{\tau\tau}^2 \,dx \,d\tau + \int \tilde\alpha^{b} \int x^2 ((1+\vartheta)x\vartheta_{x\tau\tau} - x\vartheta_x\vartheta_{\tau\tau})^2 \,dx \,d\tau \biggr\rbrace {\nonumber}\\
	& ~~~~ + C_\epsilon \biggl\lbrace \sup_\tau( \tilde\alpha^{b-a-5} + \tilde\alpha^{b-a-7} )  \cdot \int \tilde\alpha^{1+a} \int x^4 \bar\rho \vartheta_\tau^2 \,dx \,d\tau {\nonumber}\\
	& ~~~~  + \sup_\tau(\tilde\alpha^{b-a-7} + \tilde\alpha^{b-a-9} )\cdot \int \tilde\alpha^{3+a} \int x^2 ((1+\vartheta)x\vartheta_{x\tau} - x\vartheta_x\vartheta_\tau)^2\,dx\,d\tau \biggr\rbrace,
\end{align*}
where the following estimates have been applied, by applying the Hardy's inequality repeatedly
\begin{align}
	& \int x^2 (\vartheta_\tau^2 + x^2 \vartheta_{x\tau}^2)\,dx \lesssim (1+\tilde\omega) \biggl\lbrace \int x^2 \vartheta_\tau^2 \,dx + \int x^2 ((1+\vartheta)x\vartheta_{x\tau} - x\vartheta_x\vartheta_\tau)^2\,dx \biggr\rbrace {\nonumber}\\
	& ~~ \lesssim (1+\tilde\omega)\biggl\lbrace \int x^4 \bar\rho^{1/3} \vartheta_\tau^2 \,dx + \int x^4 \bar\rho^{1/3} \vartheta_{x\tau}^2 \,dx +\int x^2 ((1+\vartheta)x\vartheta_{x\tau} - x\vartheta_x\vartheta_\tau)^2\,dx \biggr\rbrace {\nonumber}\\
	& ~~ \lesssim (1+\tilde\omega) \biggl\lbrace \int x^4 \bar\rho^{1/3} \vartheta_\tau^2 \,dx + \int x^2 ((1+\vartheta)x\vartheta_{x\tau} - x\vartheta_x\vartheta_\tau)^2\,dx  \biggr\rbrace {\nonumber}\\
	& ~~  \lesssim (1+\tilde\omega) \biggl\lbrace \int x^4 \bar\rho \vartheta_\tau^2\,dx + \int x^2 ((1+\vartheta)x\vartheta_{x\tau} - x\vartheta_x\vartheta_\tau)^2\,dx\biggr\rbrace,  {\label{ine:102}}\\
	& \int x^2 (\vartheta_{\tau\tau}^2 + x^2 \vartheta_{x\tau\tau}^2)\,dx \lesssim (1+\tilde\omega) \biggl\lbrace\int x^4 \bar\rho \vartheta_{\tau\tau}^2\,dx + \int x^2 ((1+\vartheta)x\vartheta_{x\tau\tau} - x\vartheta_x\vartheta_{\tau\tau})^2\,dx \biggr\rbrace. {\label{ine:103}}
\end{align}
From the above inequalities, it is easy to see, the right of \eqref{ine:101} is integrable and bounded by $ \initial $ provided $ b $ satisfies the conditions,
\begin{equation}\label{constraint-003}
b \leq a-1.	
\end{equation}
After choosing $ b = a- 1 < 0 $ and $ \epsilon $ small enough, \eqref{ine:101} yields the estimate, together with Lemma \ref{lm:isentropic-linearly-basic-energy}, 
\begin{equation}\label{estimates-006}
\begin{aligned}
& \tilde\alpha^{b-2} \int x^4\bar\rho\vartheta_{\tau\tau}^2 \,dx + (1-\tilde\omega)\int \tilde\alpha^{b-2} \int x^4 \bar\rho \vartheta_{\tau\tau}^2 \,dx \,d\tau \\
& ~~~~~~ + (1-\tilde\omega)\int \tilde\alpha^{b} \int x^2((1+\vartheta)x\vartheta_{x\tau\tau}-x\vartheta_x\vartheta_{\tau\tau})^2 \,dx \,d\tau \lesssim \mathcal{\tilde E}_0,
\end{aligned}
\end{equation}
and consequently  \eqref{lm:101} with $ \tilde\omega < \varepsilon_0 $ small enough.
\end{pf}


\subsection{Interior Estimates}\label{sec:interior-est-isentropic}

In the following, we will perform the interior estimates. To do so, recall that the interior cut-off function is defined as
\begin{equation*}\tag{\ref{interiorCutOffFunction}}
	\chi(x) : = \begin{cases}
		1 & 0 \leq x \leq 1/2, \\
		0 & 3/4 \leq x \leq 1,
	\end{cases} 
\end{equation*}
and $ -4 \leq \chi'(x) \leq 0 $. The following lemma explores some interior estimates which can be bounded by the energy estimates obtained in the previous section. 

\begin{lm}\label{lm:isentropic-linearly-interior-estimates-1}
Under the same assumption as in Lemma \ref{lm:isentropic-linearly-temporal-derivative}, for $\varepsilon_0 $ small enough, we have the following estimates concerning the regularity of $ \theta_\tau $,
\begin{gather}
	\int \tilde\alpha^{1+a} \int \chi (  \vartheta_\tau^2 + x^2 \vartheta_{x\tau}^2 ) \,dx \,d\tau \lesssim \mathcal{\tilde E}_0  ,  \label{lm:301}\\
	\int x^2 \lbrack (1+\vartheta)x\vartheta_{x\tau} - x\vartheta_{x}\vartheta_\tau \rbrack^2 \,dx \lesssim \tilde\alpha^{-a-1} \mathcal{\tilde E}_0, \label{lm:302}\\
	\int \chi (\vartheta_\tau^2 + x^2 \vartheta_{x\tau}^2 )\,dx \lesssim \tilde \alpha^{-a-1} \mathcal{\tilde E}_0, \label{lm:303} \\
	\int \chi (\vartheta^2 + x^2 \vartheta_x^2 )\,dx \lesssim \mathcal{\tilde E}_0. \label{lm:304}
\end{gather}

\end{lm}

\begin{pf}
Multiply \eqref{eq:perturbationrbation=linearly} with $ \chi x \vartheta_\tau $ and integrate the resulting equations in the spatial variable. After integration by parts, it follows,
\begin{equation}\label{ene:le=004}
   \dfrac{4\mu}{3} \tilde\alpha^3 \int \chi \biggl\lbrace \dfrac{(\vartheta_\tau+x\vartheta_{x\tau})^2}{1+\vartheta+x\vartheta_x} + \dfrac{(1+\vartheta+x\vartheta_x)\vartheta_\tau^2}{(1+\vartheta)^2} \biggr\rbrace \,dx =  \tilde I_1 + \tilde I_2 + \tilde I_3 + \tilde I_4,  
\end{equation}
where
\begin{align*}
	& \tilde I_1 := - \dfrac{8\mu}{3}\tilde\alpha^3 \int\chi' \dfrac{x\vartheta_\tau^2}{1+\vartheta} \,dx - \tilde\alpha^3 \int \chi' x \mathfrak{\tilde B} \theta_\tau\,dx {\nonumber} \\
   & ~~~~~~~~ + \int \chi' \bar\rho^{4/3} \bigl(\dfrac{x\vartheta_\tau}{(1+\vartheta)^{8/3}(1+\vartheta+x\vartheta_x)^{4/3}}  - \dfrac{x\vartheta_\tau}{(1+\vartheta)^4}\bigr)\,dx, \\
	& \tilde I_2 := \int \chi\bar\rho^{4/3} \bigl( \dfrac{\vartheta_\tau+x\vartheta_{x\tau}}{(1+\vartheta)^{8/3}(1+\vartheta+x\vartheta_x)^{4/3}} - \dfrac{\vartheta_\tau+x\vartheta_{x\tau}}{(1+\vartheta)^4} + \dfrac{4x \vartheta_x \vartheta_\tau }{(1+\vartheta)^5} \bigr)\,dx, \\
	& \tilde I_3 := - \tilde\alpha \int \chi \dfrac{x^2 \bar\rho \vartheta_{\tau\tau}\vartheta_\tau}{(1+\vartheta)^2}  \,dx - \tilde\alpha_\tau  \int \chi \dfrac{x^2\bar\rho\vartheta_\tau^2}{(1+\vartheta)^2} \,dx, \\
	& \tilde I_4 := - \delta\int \chi \bigl(\dfrac{1}{1+\vartheta} - \dfrac{1}{(1+\vartheta)^4} \bigr) x^2 \bar\rho\vartheta_\tau\,dx .
\end{align*}

Before employing temporal-weighted estimates on \eqref{ene:le=004}, we shall start with analysing the right. Indeed, the following estimates hold for $ \tilde I_k $'s.
\begin{align*}
	& \tilde I_1 \lesssim \tilde\alpha^{3} \int x^4 \bar\rho \vartheta_\tau^2 \,dx + \tilde\alpha^{3} \int x^2 ((1+\vartheta)x\vartheta_{x\tau}-x\vartheta_x\vartheta_\tau)^2\,dx \\
	& ~~~~ + \bigr(\int x^4\bar\rho\vartheta_\tau^2 \,dx \bigl)^{1/2} \bigr( \int x^4 \bar\rho\vartheta^2 \,dx + \int x^4 \vartheta_x^2 \,dx \bigr)^{1/2},\\
	& \tilde I_2 \lesssim \bigl( \int \chi ((\vartheta + x\vartheta_{x})^2 + \vartheta^2) \,dx\bigr)^{1/2}\bigl( \int \chi ((\vartheta_\tau + x\vartheta_{x\tau})^2 + \vartheta_\tau^2) \,dx\bigr)^{1/2},  \\
	& \tilde I_3 \lesssim \tilde\alpha \biggl\lbrace \bigl( \int x^4 \bar\rho \vartheta_{\tau\tau}^2 \,dx \bigr)^{1/2} + \bigl( \int x^4 \bar\rho \vartheta_{\tau}^2 \,dx \bigr)^{1/2} \biggr\rbrace \times  \bigl( \int \chi \vartheta_\tau^2 \,dx \bigr)^{1/2}, \\
	& \tilde I_4 \lesssim \bigl( \int \chi\vartheta^2\,dx  \bigr)^{1/2}\bigl( \int \chi \vartheta_\tau^2\,dx \bigr)^{1/2}.
\end{align*}
Now multiply \eqref{ene:le=004} with $ \tilde\alpha^{c} $ and apply the Cauchy's inequality to the resulting equation. It holds,
\begin{equation}\label{ene:le=005}
	\begin{aligned}
		& (1-\epsilon)\tilde\alpha^{3+c} \int\chi ((\vartheta_\tau+x\vartheta_{x\tau})^2 + \vartheta_\tau^2)\,dx \lesssim C_\epsilon \biggl\lbrace \tilde\alpha^{c-3}\int \chi ((\vartheta + x\vartheta_{x})^2 + \vartheta^2) \,dx \\
		& ~~~~ + \tilde\alpha^{c-1} \int x^4 \bar\rho \vartheta_{\tau\tau}^2 \,dx + \tilde\alpha^{c-1} \int x^4 \bar\rho \vartheta_\tau^2\,dx  \biggr\rbrace\\
		& ~~~~ + \tilde\alpha^{3+c} \int x^4\bar\rho \vartheta_\tau^2 \,dx + \tilde\alpha^{3+c} \int x^2 ((1+\vartheta)x\vartheta_{x\tau} - x\vartheta_x\vartheta_\tau)^2 \,dx \\
		& ~~~~ + \biggl\lbrace \int x^4\bar\rho \vartheta^2\,dx + \int x^4 \vartheta_x^2 \,dx \biggr\rbrace \times \tilde\alpha^{c-3}.
	\end{aligned}
\end{equation}
In addition, 
\begin{align*}
	& \bigl( \int \chi ((\vartheta+x\vartheta_x)^2 + \vartheta^2 )\,dx \bigr)^{1/2} \dfrac{d}{d\tau} \bigl( \int \chi ((\vartheta+x\vartheta_x)^2 + \vartheta^2 )\,dx\bigr)^{1/2} = \dfrac{1}{2} \dfrac{d}{d\tau} \int \chi ((\vartheta+x\vartheta_x)^2 + \vartheta^2 )\,dx \\
	& \leq \bigl( \int \chi ((\vartheta+x\vartheta_x)^2 + \vartheta^2 )\,dx \bigr)^{1/2} \bigl(\int \chi ((\vartheta_\tau+x\vartheta_{x\tau})^2 + \vartheta_\tau^2 )\,dx \bigr)^{1/2},
\end{align*}
from which it follows,
\begin{equation}\label{ine:105}
	\bigl( \int \chi ((\vartheta+x\vartheta_x)^2 + \vartheta^2 )\,dx \bigr)^{1/2}\lesssim \mathcal{\tilde E}_{0}^{1/2} + \int\bigl( \int \chi ((\vartheta_\tau+x\vartheta_{x\tau})^2 + \vartheta_\tau^2 )\,dx \bigr)^{1/2} \,d\tau.
\end{equation}
Then after choosing $ \epsilon = 1/2 $ in \eqref{ene:le=005}, the following estimate holds,
\begin{equation}\label{ene:le=006}
	\begin{aligned}
		& \tilde\alpha^{3+c} \int \chi ((\vartheta_\tau + x\vartheta_{x\tau})^2 + \vartheta_\tau^2 )\,dx \lesssim \tilde\alpha^{c-3} \cdot \int \tilde\alpha^{-(3+c)}\,d\tau \\
		& ~~~~ \times \int \tilde\alpha^{3+c} \int \chi ((\vartheta_\tau + x \vartheta_{x\tau})^2 + \vartheta_\tau^2 )\,dx\,d\tau + \tilde R_1 , 
	\end{aligned}
\end{equation}
where 
\begin{align*}
	& \tilde R_1 : = \mathcal{\tilde E}_0 \cdot \tilde\alpha^{c-3} +  \tilde\alpha^{c-b+1} \cdot \tilde\alpha^{b-2} \int x^4 \bar\rho \vartheta_{\tau\tau}^2 \,dx \\
	& ~~~~ + \tilde\alpha^{c-a-2} \cdot \tilde\alpha^{1+a} \int x^4 \bar\rho \vartheta_\tau^2\,dx  + \tilde\alpha^{c-a+2} \cdot \tilde\alpha^{1+a} \int x^4\bar\rho \vartheta_\tau^2 \,dx \\
	& ~~~~ +\tilde\alpha^{c-a} \cdot \tilde\alpha^{3+a}\int x^2 ((1+\vartheta)x\vartheta_{x\tau} - x\vartheta_x\vartheta_\tau)^2 \,dx \\
	& ~~~~ + \tilde\alpha^{c-3} \cdot  \sup_\tau \bigl( \int x^4\bar\rho \vartheta^2\,dx + \int x^4 \vartheta_x^2 \,dx \bigr).
\end{align*}
Then by applying the Gr\"onwall's inequality, it admits
\begin{equation*}\tag{\ref{lm:301}} \label{ine:104}
	\begin{aligned}
		& \int \tilde\alpha^{3+c} \int \chi ((\vartheta_\tau+x\vartheta_{x\tau})^2 + \vartheta_\tau^2 )\,dx \,d\tau \lesssim e^{C \int \tilde\alpha^{c-3} \,d\tau \cdot \int \tilde\alpha^{-(3+c)}\,d\tau} \\
		& ~~ \times \int \tilde R_1 \,d\tau \lesssim \mathcal{\tilde E}_0,
	\end{aligned}
\end{equation*}
provided
\begin{equation}\label{constraint-004}
	- 3 < c \leq a-2.
\end{equation}
From now on, we choose $ c = a-2 $.

In the meantime, from \eqref{ene:le=001}, together with \eqref{ine:102},
\begin{align*}
	& \tilde\alpha^3 \int x^2 ((1+\vartheta)x\vartheta_{x\tau}-x\vartheta_x\vartheta_\tau)^2 \,dx 
	\lesssim \epsilon \tilde\alpha^3 \int x^2((1+\vartheta)x\vartheta_{x\tau} - x \vartheta_x\vartheta_\tau)^2 \,dx\\
	& ~~~~ + C_\epsilon \tilde\alpha^{-3} \bigl( \int x^4 \bar\rho\vartheta^2 \,dx + \int x^4 \vartheta_x^2 \,dx \bigr) \\
	& ~~~~ + \tilde\alpha^{-(1+a)/2 }\bigl( \int x^4\bar\rho \vartheta^2\,dx \bigr)^{1/2}\bigl( \tilde\alpha^{1+a} \int x^4\bar\rho \vartheta_\tau^2 \,dx  \bigr)^{1/2}\\
	& ~~~~ + \tilde\alpha^{(3-a-b)/2} \bigl( \tilde\alpha^{1+a} \int x^4 \bar\rho \vartheta_\tau^2 \,dx \bigr)^{1/2} \bigl( \tilde\alpha^{b-2} \int x^4 \bar\rho \vartheta_{\tau\tau}^2 \,dx \bigr)^{1/2}.
\end{align*}
Therefore, 
\begin{equation*}\tag{\ref{lm:302}}\label{estimates-005}
	\int x^2 ((1+\vartheta)x\vartheta_{x\tau} - x\vartheta_x\vartheta_\tau)^2 \,dx \lesssim \tilde\alpha^{d_1} \cdot \mathcal{\tilde E}_0 ,
\end{equation*}
where $$ d_1 = \max\lbrace -6, -(7+a)/2, -(3+a+b)/2 \rbrace = -a - 1. $$
Meanwhile, from and \eqref{ene:le=006}, 
\begin{align*}
	& \int \chi ((\vartheta_\tau + x\vartheta_{x\tau})^2 + \vartheta_\tau^2 )\,dx \lesssim \tilde\alpha^{-6}  \cdot \int \tilde\alpha^{-(3+c)}\,d\tau \\
	& ~~~~ \times \int \tilde\alpha^{3+c} \int \chi ((\vartheta_\tau+x\vartheta_{x\tau})^2+\vartheta_\tau^2 )\,dx\,d\tau + \tilde R_2,
\end{align*}
where\begin{align*}
	& \tilde R_2 := \tilde\alpha^{-(3+c)} \tilde R_1 = \mathcal{\tilde E}_0 \cdot \tilde\alpha^{-6}+ \tilde\alpha^{-b-2} \cdot \sup_\tau \tilde\alpha^{b-2} \int x^4 \bar\rho\vartheta_{\tau\tau}^2 \,dx \\
	& ~~~~ + ( \tilde\alpha^{-a-5} + \tilde\alpha^{-a-1} ) \cdot \sup_\tau \tilde\alpha^{1+a} \int x^4\bar\rho \vartheta_\tau^2 \,dx \\
	& ~~~~ + \int x^2 ((1+\vartheta)x\vartheta_{x\tau} - x\vartheta_x\vartheta_\tau)^2\,dx + \tilde\alpha^{-6} \cdot \sup_\tau \bigl( \int x^4\bar\rho \vartheta^2\,dx + \int x^4 \vartheta_x^2 \,dx \bigr)  .
\end{align*}
Consequently, 
\begin{equation*}\tag{\ref{lm:303}}\label{estimates-004}
	\int \chi ((\vartheta_\tau + x\vartheta_{x\tau})^2 + \vartheta_\tau^2 )\,dx\lesssim \tilde\alpha^{d_2} \cdot \mathcal{\tilde E}_0 ,	
\end{equation*}
where
$$ d_2 = \max\lbrace - 6, -b-2, -a-5, -a-1, d_1 \rbrace = - a - 1. $$
As a corollary, from \eqref{ine:105}, \eqref{ine:104},
\begin{equation*}\tag{\ref{lm:304}}\label{estimates-007}
\int \chi (x^2 \vartheta_x^2 + \vartheta^2) \,dx \lesssim \mathcal{\tilde E}_0 + \int\tilde\alpha^{-(3+c)} \,d\tau \cdot \int \tilde\alpha^{3+c} \int \chi ((\vartheta_\tau + x\vartheta_{x\tau})^2 + \vartheta_\tau^2) \,dx\,d\tau \lesssim \mathcal{\tilde E}_0.	
\end{equation*}
\end{pf}



In addition to the regularity obtained above, in the following two lemmas, we will derive some more energy estimates precisely concerning the interior regularity. These will be the last block to close the stability analysis. 

\begin{lm}\label{lm:isentropic-linearly-interior-estimates-2}
Under the same assumption as in Lemma \ref{lm:isentropic-linearly-interior-estimates-1}, for $ \varepsilon_0 $ small enough, we have the following estimate, 
\begin{equation}\label{lm:401}
\tilde\alpha^{a-5} \int \chi x^2 \bar\rho \vartheta_{\tau\tau}^2 \,dx + \int \tilde\alpha^{a-5} \int \chi x^2 \bar\rho \vartheta_{\tau\tau}^2 \,dx \,d\tau + \int \tilde\alpha^{a-3} \int \chi (\vartheta_{\tau\tau}^2 + x^2 \vartheta_{x\tau\tau}^2 ) \,dx\,d\tau \lesssim \mathcal{\tilde E}_0. 
\end{equation}

\end{lm}

\begin{pf}
Multiply \eqref{eq:temporalderivative=linearlyexpanding} with $ \tilde\alpha^{\mathfrak{a}} \chi x \vartheta_{\tau\tau} $, $ \mathfrak a < 0 $ and integrate the resulting equation in the spatial variable. After integration by parts, it follows,

\begin{equation}\label{ene:le=007}
	\dfrac{d}{d\tau} \tilde E_{3} + \tilde D_{3} = \tilde K_1 + \tilde K_2 + \tilde K_3 + \tilde K_4 + \tilde K_5 + \tilde K_6 + \tilde K_7,
\end{equation}
where
\begin{align*}
	& \tilde E_{3} : = \dfrac{\tilde\alpha^{\mathfrak a-2}}{2} \int \chi x^2 \bar \rho \vartheta_{\tau\tau}^2 \,dx, \\
	& \tilde D_{3} : = - \dfrac{\mathfrak a}{2} \tilde\alpha^{\mathfrak a -3} \tilde\alpha_\tau \int \chi x^2 \bar\rho \vartheta_{\tau\tau}^2 \,dx {\nonumber} \\
	& ~~~~~~ + \dfrac{4}{3}\mu \tilde\alpha^{\mathfrak a} \int \chi \biggl\lbrace \dfrac{(1+\vartheta)^2x^2 \vartheta_{x\tau\tau}^2}{1+\vartheta+x\vartheta_x} + \dfrac{\lbrack2(1+\vartheta+x\vartheta_x)^2 - (1+\vartheta)^2\rbrack\vartheta_{\tau\tau}^2}{1+\vartheta+x\vartheta_x} \biggr\rbrace \,dx ,  \\
	& \tilde K_1 : = - (\tilde\alpha^{\mathfrak a-3} \tilde\alpha_{\tau\tau} - 3 \tilde\alpha^{\mathfrak a -4}\tilde\alpha_\tau^2 ) \int \chi x^2 \bar\rho \vartheta_\tau\vartheta_{\tau\tau} \,dx, \\
	& \tilde K_2 : =- \int \biggl\lbrace \tilde\alpha^{\mathfrak a -3} \bigl(\vartheta_\tau + 2 \dfrac{\vartheta_\tau}{(1+\vartheta)^3}\bigr) - 3\tilde\alpha^{\mathfrak a-4}\tilde\alpha_\tau\bigl((1+\vartheta) - \dfrac{1}{(1+\vartheta)^2}\bigr) \biggr\rbrace {\nonumber} \\
	& ~~~~~~ \times \delta\chi x^2 \bar\rho \vartheta_{\tau\tau}\,dx, \\
	& \tilde K_3 : = \dfrac{4}{3}\mu\tilde\alpha^{\mathfrak a} \int \chi \biggl\lbrace (x(1+\vartheta)^2\vartheta_{\tau\tau}) _x\bigl( \dfrac{(\vartheta_\tau+x\vartheta_{x\tau})^2}{(1+\vartheta+x\vartheta_x)^2} - \dfrac{\vartheta_\tau^2}{(1+\vartheta)^2} \bigr) {\nonumber} \\
	& ~~~~~~ + 3(1+\vartheta)^2\vartheta_{\tau\tau} \bigl( - \dfrac{2 x\vartheta_{x\tau} \vartheta_\tau}{(1+\vartheta)^2} + \dfrac{2x\vartheta_x \vartheta_{\tau}^2}{(1+\vartheta)^3} \bigr) + 6 x(1+\vartheta)\vartheta_\tau\vartheta_{\tau\tau} \bigl(\dfrac{\vartheta_\tau}{1+\vartheta}\bigr)_x {\nonumber}\\
	& ~~~~~~ - 2 (x(1+\vartheta)\vartheta_\tau\vartheta_{\tau\tau})_x \bigl( \dfrac{\vartheta_\tau + x \vartheta_{x\tau}}{1+\vartheta+x\vartheta_x} - \dfrac{\vartheta_\tau}{1+\vartheta} \bigr) \biggr\rbrace \,dx  ,   \\
	& \tilde K_4 : = - \tilde\alpha^{\mathfrak a} \int \chi' \biggl\lbrace x(1+\vartheta)^2\vartheta_{\tau\tau}\mathfrak{\tilde B}_\tau + 2 x(1+\vartheta)\vartheta_{\tau}\vartheta_{\tau\tau}\mathfrak{\tilde B} + \dfrac{4}{3}\mu x(1+\vartheta)\vartheta_{\tau\tau}^2  \biggr\rbrace \,dx, \\
	& \tilde K_5 : = - 3 \tilde\alpha^{\mathfrak a-4}\tilde\alpha_\tau \int \chi \bar\rho^{4/3} \biggl\lbrace \dfrac{(x(1+\vartheta)^2\vartheta_{\tau\tau})_x}{(1+\vartheta)^{8/3}(1+\vartheta+x\vartheta_x)^{4/3}} - \biggl( \dfrac{x\vartheta_{\tau\tau}}{(1+\vartheta)^2} \biggr)_x \biggr\rbrace \, dx,       \\
	& \tilde K_6 : =  \tilde\alpha^{\mathfrak a-3} \int \chi \bar\rho^{4/3} \biggl\lbrace (\vartheta_{\tau\tau} + x\vartheta_{x\tau\tau}) \bigl\lbrack \dfrac{1}{(1+\vartheta)^{2/3}(1+\vartheta+x\vartheta_x)^{4/3}} - \dfrac{1}{(1+\vartheta)^2} \bigr\rbrack_\tau {\nonumber} \\
	& ~~~~~~ + \vartheta_{\tau\tau} \bigl\lbrack \dfrac{2 x\vartheta_x}{(1+\vartheta)^{5/3}(1+\vartheta+x\vartheta_x)^{4/3}} + \dfrac{2x\vartheta_x}{(1+\vartheta)^3} \bigr\rbrack_\tau \biggr\rbrace \,dx,   \\
	& \tilde K_7 : = - 3 \tilde\alpha^{\mathfrak a-4}\tilde\alpha_\tau \int \chi'  x \bar\rho^{4/3} \biggl\lbrace \dfrac{\vartheta_{\tau\tau}}{(1+\vartheta)^{2/3}(1+\vartheta+x\vartheta_x)^{4/3}} - \dfrac{\vartheta_{\tau\tau}}{(1+\vartheta)^2} \biggr\rbrace \,dx {\nonumber} \\
	& ~~~~~~ + \tilde\alpha^{\mathfrak a-3} \int \chi' x \bar\rho^{4/3} \vartheta_{\tau\tau} \biggl\lbrace \dfrac{1}{(1+\vartheta)^{2/3}(1+\vartheta+x\vartheta_x)^{4/3}}  - \dfrac{1}{(1+\vartheta)^2} \biggr\rbrace_\tau \,dx.
\end{align*}
Integrating \eqref{ene:le=007} over the temporal variable yields, 
\begin{equation}\label{ine:201}
\tilde E_3 + \int \tilde D_3 \,d\tau \lesssim \sum_{k=1}^{7} \int \tilde K_k + \mathcal{\tilde E}_0.	
\end{equation}
Similarly as before, the left of \eqref{ine:201} admits,
\begin{equation*}
	\begin{aligned}
		& \tilde E_3 + \int \tilde D_3 \,d\tau \gtrsim \tilde\alpha^{\mathfrak a-2} \int \chi x^2 \bar\rho \vartheta_{\tau\tau}^2 \,dx + |\mathfrak a| \int \tilde\alpha^{\mathfrak a -2} \int \chi x^2 \bar\rho \vartheta_{\tau\tau}^2 \,dx\,d\tau \\
		& ~~~~ + (1-\tilde \omega)  \int \tilde\alpha^{\mathfrak a } \int \chi ( x^2 \vartheta_{x\tau\tau}^2 + \vartheta_{\tau\tau}^2) \,dx \,d\tau. 
	\end{aligned}
\end{equation*}
In the following, the right of \eqref{ine:201} would be estimated.
\begin{align*}
	& \int \tilde K_1 \,d\tau \lesssim (1+\tilde\omega)\int  \tilde\alpha^{\mathfrak a - 2} \bigl( \int \chi \vartheta_{\tau\tau}^2 \,dx \bigr)^{1/2} \bigl( \int \chi x^4 \bar\rho\vartheta_\tau^2 \,dx  \bigr)^{1/2} \,d\tau {\nonumber} \\
	& ~~ \lesssim  \epsilon \int  \tilde\alpha^{\mathfrak a} \int \chi \vartheta_{\tau\tau}^2 \,dx \,d\tau + C_\epsilon \sup_\tau \tilde\alpha^{\mathfrak a - a - 5} \cdot \int\tilde\alpha^{1+a} \int x^4 \bar\rho \vartheta_\tau^2 \,dx \,d\tau ,   \\
	& \int \tilde K_2 \,d\tau \lesssim (1+\tilde\omega) \int \tilde\alpha^{\mathfrak a - 3} \bigl( \int \chi\vartheta_{\tau\tau}^2 \,dx \bigr)^{1/2} \bigl( \int x^4 \bar\rho \vartheta^2 \,dx  + \int x^4 \bar\rho \vartheta_\tau^2 \,dx \bigr)^{1/2} \,d\tau {\nonumber} \\
	& ~~ \lesssim \epsilon \int \tilde\alpha^{\mathfrak a } \int \chi \vartheta_{\tau\tau}^2 \,dx\,d\tau + C_\epsilon \biggl\lbrace \int \tilde\alpha^{\mathfrak a - 6} \,d\tau \cdot \sup_\tau \int x^4 \bar\rho \vartheta^2 \,dx  {\nonumber} \\ & ~~~~~~ + \sup_\tau \tilde\alpha^{\mathfrak a - a- 7} \cdot \int \tilde\alpha^{1+a}\int x^4 \bar\rho \vartheta_\tau^2\,dx  \,d\tau \biggr\rbrace,   \\
	& \int \tilde K_3 \,d\tau \lesssim \tilde\omega \int \tilde\alpha^{\mathfrak a} \int \chi (|x\vartheta_{x\tau\tau}|+|\vartheta_{\tau\tau}|)(|x\vartheta_{x\tau}|+|\vartheta_\tau|) \,dx \,d\tau {\nonumber} \lesssim \epsilon \int \tilde\alpha^{\mathfrak a} \int \chi (x^2 \vartheta_{x\tau\tau}^2 + \vartheta_{\tau\tau}^2 ) \,dx {\nonumber} \\
	& ~~~~~~ + \tilde \omega C_\epsilon \sup_\tau \tilde\alpha^{\mathfrak a - c - 3} \cdot \int \tilde\alpha^{3+c} \int \chi ((\vartheta_\tau+x\vartheta_{x\tau})^2 + \vartheta_\tau^2 )\,dx\,d\tau,  \\
	& \int \tilde K_4 \,d\tau \lesssim \int \tilde\alpha^{\mathfrak a} \int x^3 \bar\rho^{1/2} |\vartheta_{\tau\tau}| \bigl( |(1+\vartheta)x \vartheta_{x\tau\tau} - x\vartheta_x \vartheta_\tau| + x \bar\rho^{1/2}|\vartheta_{\tau\tau}| {\nonumber}\\
	& ~~~~~~ + |(1+\vartheta)x\vartheta_{x\tau} - x\vartheta_x\vartheta_\tau| + x \bar\rho^{1/2} |\vartheta_\tau|\bigr) \,dx\,d\tau {\nonumber}\\
	& ~~ \lesssim \sup_\tau \lbrace \tilde\alpha^{2\mathfrak a - 2b + 2} + \tilde\alpha^{\mathfrak a -b + 2} + \tilde\alpha^{2\mathfrak a - a -b -1} + \tilde\alpha^{2\mathfrak a -a -b +1} \rbrace \cdot \int \tilde\alpha^{b-2} \int x^4\bar\rho \vartheta_{\tau\tau}^2 \,dx \,d\tau {\nonumber} \\
	& ~~~~~~  + \int \tilde\alpha^{b} \int x^2 ((1+\vartheta)x\vartheta_{x\tau\tau} - x\vartheta_x\vartheta_{\tau\tau})^2 \,dx\,d\tau {\nonumber}\\
	& ~~~~~~ + \int \tilde\alpha^{3+a} \int x^2 ((1+\vartheta)x\vartheta_{x\tau}-x\vartheta_x\vartheta_\tau)^2 \,dx\,d\tau + \int \tilde\alpha^{1+a} \int x^4 \bar\rho\vartheta_\tau^2\,dx\,d\tau,    \\
	& \int \tilde K_5 \,d\tau \lesssim \int \tilde\alpha^{\mathfrak a - 3}\int \chi (|x\vartheta_{x\tau\tau}|+|\vartheta_{\tau\tau}|)(|x\vartheta_{x}|+ |\vartheta|) \,dx \,d\tau  {\nonumber} \\
	& ~~ \lesssim \epsilon \int \tilde\alpha^{\mathfrak a} \int \chi (x^2 \vartheta_{x\tau\tau}^2 + \vartheta_{\tau\tau}^2) \,dx \, d\tau + C_\epsilon \int \tilde\alpha^{\mathfrak a - 6} \,d\tau \times \sup_\tau \int \chi (x^2 \vartheta_x^2 + \vartheta^2) \,dx,\\
	& \int \tilde K_6 \,d\tau \lesssim \int \tilde\alpha^{\mathfrak a - 3} \int \chi (|x\vartheta_{x\tau\tau}|+|\vartheta_{\tau\tau}|)(|x\vartheta_{x\tau}|+ |\vartheta_\tau|)  \,dx {\nonumber}\\
	& ~~ \lesssim \epsilon \int \tilde\alpha^{\mathfrak a } \int \chi (x^2 \vartheta_{x\tau\tau}^2 + \vartheta_{\tau\tau}^2)\,dx\,d\tau + C_\epsilon \sup_\tau \tilde\alpha^{\mathfrak a - c - 9} \cdot \int \tilde\alpha^{3+c}  \int \chi (x^2 \vartheta_{x\tau}^2 + \vartheta_\tau^2 )\,dx \,d\tau,   \\
	& \int \tilde K_7 \,d\tau \lesssim \int \tilde\alpha^{\mathfrak a - 3} \int x^4 \bar\rho |\vartheta_{\tau\tau}| ( |x\vartheta_x| + |\vartheta| ) \,dx \,d\tau + \int \tilde\alpha^{\mathfrak a - 3} \int x^4 \bar\rho |\vartheta_{\tau\tau}| (|(1+\vartheta)x\vartheta_{x\tau} \\
	& ~~~~~~ - x\vartheta_x\vartheta_\tau| + |\vartheta_\tau|) \,dx\,d\tau {\nonumber} \lesssim \sup_\tau \bigl\lbrace  \int \tilde\alpha^{2\mathfrak a - b - 4} \,d\tau  +  \tilde\alpha^{2\mathfrak a - a -b -5} + \tilde\alpha^{2\mathfrak a - a- b -7} \bigr\rbrace\\
	& ~~~~~~ \times \int \tilde\alpha^{b-2} \int x^4 \bar\rho\vartheta_{\tau\tau}^2\,dx \,d\tau {\nonumber} + \sup_\tau \bigl\lbrace \int x^4 \bar\rho \vartheta^2 \,dx + \int x^4 \vartheta_x^2 \,dx \bigr\rbrace\\
	& ~~~~~~  + \int \tilde\alpha^{1+a} \int x^4 \bar\rho \vartheta_\tau^2 \,dx\,d\tau {\nonumber} + \int \tilde\alpha^{3+a}\int x^2 ((1+\vartheta)x\vartheta_{x\vartheta} - x\vartheta_x\vartheta_\tau)^2 \,dx \,d\tau, 
\end{align*}
where \eqref{ine:102} and \eqref{ine:103} are applied.
Similarly as before, it shall be imposed the following constraints,
\begin{equation}\label{constraint-101}
	\mathfrak a \leq 
	a - 3 < 0 .
\end{equation}
Therefore, after choosing $ \mathfrak a = a-3 $ and $ \epsilon $ small enough, depending on $ |\mathfrak a| $, \eqref{ine:201} implies,
\begin{equation*}\tag{\ref{lm:401}}\label{estimates=102}
	\begin{aligned}
		& \tilde\alpha^{\mathfrak a -2} \int \chi x^2 \bar\rho\vartheta_{\tau\tau}^2 \,dx + (1-\tilde\omega) \int \tilde\alpha^{\mathfrak a -2 } \int \chi x^2 \bar\rho \vartheta_{\tau\tau}^2 \,dx \,d\tau \\
		& ~~~~~~ + (1-\tilde\omega) \int \tilde\alpha^{\mathfrak a } \int \chi (x^2 \vartheta_{x\tau\tau}^2 + \vartheta_{\tau\tau}^2)\,dx\,d\tau \lesssim \mathcal{\tilde E}_0.
	\end{aligned}
\end{equation*}
\end{pf}

In the following, I shall present the elliptic structure of \eqref{eq:perturbationrbation=linearly}. To start with, recalling that the relative entropy is defined as
\begin{equation*}\tag{\ref{isen:RelativeEntropy}}
	\mathcal G = \log{(1+\vartheta)^2(1+\vartheta+x\vartheta_x)}.
\end{equation*}
Then \eqref{eq:perturbationrbation=linearly} can be written as
\begin{equation}\label{eq:alternative=linearly}
	\begin{aligned}
		& \dfrac{4\mu}{3} \tilde \alpha^3 \mathcal G_{x\tau} + \dfrac{4}{3} \dfrac{\bar\rho^{4/3}}{(1+\vartheta)^{8/3}(1+\vartheta+x\vartheta_x)^{4/3}} \mathcal G_x \\
		&  ~~ = (\bar\rho^{4/3})_x\biggl\lbrace \dfrac{1}{(1+\vartheta)^{8/3}(1+\vartheta+x\vartheta_x)^{4/3}} - \dfrac{1}{(1+\vartheta)^4} \biggr\rbrace \\
		& ~~~~~~ + \dfrac{1}{(1+\vartheta)^2} \bigl( \tilde\alpha x \bar\rho \vartheta_{\tau\tau} + \tilde\alpha_\tau x\bar\rho \vartheta_\tau\bigr) + \bigl( \dfrac{1}{1+\vartheta} - \dfrac{1}{(1+\vartheta)^4}\bigr) \delta x\bar\rho.
	\end{aligned}
\end{equation}

\begin{lm}\label{lm:isentropic-linearly-interior-estimates-3}
Under the same assumption as in Lemma \ref{lm:isentropic-linearly-interior-estimates-2}, for $ \varepsilon_0 $ small enough, the following estimates of the relative entropy $ \mathcal G $ hold,
\begin{gather}
\int \mathcal G_x^2 \,dx + \int \tilde \alpha^{-3} \int \bar\rho^{4/3} \mathcal G_x^2 \,dx \,d\tau \lesssim \mathcal{\tilde E}_0,  \label{lm:501}\\
\int \mathcal G_{x\tau}^2 \,dx \lesssim \tilde\alpha^{1-a} \mathcal{\tilde E}_0. \label{lm:502}
\end{gather}
As a consequence, applying Lemma \ref{lm:estimates-of-relative-entropy} implies
\begin{gather}
	\norm{\vartheta_x}{L_x^2}^2 + \norm{x\vartheta_{xx}}{L_x^2}^2  \lesssim \mathcal{\tilde E}_0, \label{lm:503}\\
	\norm{\vartheta_{x\tau}}{L_x^2}^2 + \norm{x\vartheta_{xx\tau}}{L_x^2}^2  \lesssim (1+\tilde\alpha^{1-a})\mathcal{\tilde E}_0. \label{lm:504}
\end{gather}

\end{lm}

\begin{pf}

Multiply \eqref{eq:alternative=linearly} with $ \tilde\alpha^{-3} \mathcal G_x $ and integrate the resulting equation. It holds,
\begin{equation}\label{ene:le=301}
	\begin{aligned}
		& \dfrac{d}{d\tau} \dfrac{2\mu}{3} \int \mathcal G_{x}^2 \,dx 
		+ \dfrac{4}{3} \tilde\alpha^{-3} \int \dfrac{\bar\rho^{4/3} \mathcal G_{x}^2}{(1+\vartheta)^{8/3}(1+\vartheta+x\vartheta_x)^{4/3}}\,dx \\
		& ~~ = \tilde\alpha^{-3} \int  \biggl\lbrace  \dfrac{(\bar\rho^{4/3})_x}{(1+\vartheta)^{8/3}(1+\vartheta+x\vartheta_x)^{4/3}} - \dfrac{(\bar\rho^{4/3})_x}{(1+\vartheta)^4} \biggr\rbrace \mathcal G_x\,dx \\
		& ~~~~ + \tilde\alpha^{-2} \int \dfrac{ x \bar\rho \vartheta_{\tau\tau}  }{(1+\vartheta)^2} \mathcal G_x \,dx + \tilde\alpha^{-3} \tilde\alpha_\tau \int \dfrac{x \bar\rho \vartheta_\tau}{(1+\vartheta)^2} \mathcal G_x \,dx \\
		& ~~~~ + \tilde\alpha^{-3} \int \biggl( \dfrac{1}{1+\vartheta} - \dfrac{1}{(1+\vartheta)^4}\biggr) \delta x\bar\rho \mathcal G_x \,dx : = \sum_{k=1}^{4}\mathfrak{\tilde L}_k.
	\end{aligned}
\end{equation}
Now it is time to estimate the right of \eqref{ene:le=301}. Indeed, by noticing 
$$ \abs{(\bar\rho^{4/3})_x}{} \lesssim x \bar\rho, $$
the following estimates hold.
\begin{align*}
& \mathfrak{\tilde L}_1, \mathfrak{\tilde L}_4 \lesssim \tilde\alpha^{- 3} \int x\bar\rho (|\vartheta|+|x\vartheta_x|) \mathcal G_x \,dx \lesssim \epsilon \tilde \alpha^{-3} \int \dfrac{\bar\rho^{4/3} \mathcal G_{x}^2}{(1+\vartheta)^{8/3}(1+\vartheta+x\vartheta_x)^{4/3}} \,dx {\nonumber}\\
& ~~~~~~ + C_\epsilon  \tilde\alpha^{- 3} \int x^2 \bar\rho^{2/3}(\vartheta^2 + x^2 \vartheta_x^2 )\,dx \lesssim \epsilon \tilde \alpha^{-3} \int \dfrac{\bar\rho^{4/3} \mathcal G_{x}^2}{(1+\vartheta)^{8/3}(1+\vartheta+x\vartheta_x)^{4/3}} \,dx {\nonumber} \\
& ~~~~~~ + C_\epsilon  \tilde\alpha^{ - 3} \bigl( \int x^4 \bar\rho \vartheta^2 \,dx + \int x^4 \vartheta_x^2 \,dx\bigr), \\
& \mathfrak{\tilde L}_2 \lesssim \tilde\alpha^{-2} \int x\bar\rho |\vartheta_{\tau\tau}| |\mathcal{G}_x| \,dx 
\lesssim \epsilon \tilde\alpha^{-2-b} \int \mathcal G_x^2 \,dx + C_\epsilon \tilde\alpha^{b-2} \int x^2 \bar\rho \vartheta_{\tau\tau}^2 \,dx {\nonumber}\\
& ~~ \lesssim \epsilon \tilde\alpha^{-2-b} \cdot\sup_\tau \int \mathcal G_x^2\,dx + C_\epsilon \tilde\alpha^{b-2} \int x^4 \bar\rho\vartheta_{\tau\tau}^2 \,dx  {\nonumber} \\
& ~~~~~~ + C_\epsilon \tilde\alpha^{-2}\cdot \tilde\alpha^{b} \int x^2 ((1+\vartheta)x\vartheta_{x\tau\tau} - x\vartheta_x\vartheta_{\tau\tau})^2 \,dx, \\
& \mathfrak{\tilde L}_3 \lesssim \tilde\alpha^{-2} \int x\bar\rho |\vartheta_\tau||\mathcal G_x|\,dx \lesssim \epsilon \tilde\alpha^{-3} \int \dfrac{\bar\rho^{4/3}\mathcal G_x^2}{(1+\vartheta)^{8/3}(1+\vartheta+x\vartheta_x)^{4/3}} \,dx {\nonumber} \\
& ~~~~~~ + C_\epsilon \tilde\alpha^{-a-2} \cdot \tilde\alpha^{1+a}\int x^4 \bar\rho\vartheta_\tau^2 \,dx + C_\epsilon \tilde\alpha^{-a-4} \cdot\tilde\alpha^{3+a} \int x^2 ((1+\vartheta)x\vartheta_{x\tau} - x\vartheta_x\vartheta_\tau)^2 \,dx,  
\end{align*}
where \eqref{ine:102} and \eqref{ine:103} are applied. Therefore, after integrating \eqref{ene:le=301} in the temporal variable and choosing $ \epsilon $ small enough, it admits,
\begin{equation}\label{estimates-101}
	\begin{aligned}
		& \int \mathcal G_x^2\,dx + \int \tilde\alpha^{-3} \int \bar\rho^{4/3} \mathcal G_x^2 \,dx \,d\tau \lesssim \mathcal{\tilde E}_0.
	\end{aligned}
\end{equation}
Notice, $ b = a - 1 \in (-1,0) $.
On the other hand, from \eqref{eq:alternative=linearly}, 
\begin{align*}
	& \tilde \alpha^6 \int  \mathcal G_{x\tau}^2 \,dx \lesssim \int  \mathcal G_x^2 \,dx + \int x^2 \bar\rho (\vartheta^2 + x^2\vartheta_x^2)\,dx + \tilde \alpha^2 \biggl\lbrace \int x^2 \bar\rho \vartheta_{\tau\tau}^2 \,dx + \int x^2 \bar\rho \vartheta_{\tau}^2 \,dx  \biggr\rbrace \\
	& ~~ \lesssim \int  \mathcal G_x^2 \,dx + \int (x^4 \bar\rho \vartheta^2 + x^4\vartheta_x^2)\,dx  + \tilde\alpha^{4 - \mathfrak a } \cdot \tilde\alpha^{\mathfrak a - 2} \int \chi x^2 \bar\rho \vartheta_{\tau\tau}^2 \,dx + \tilde\alpha^{4-b} \cdot \tilde\alpha^{b-2} \int x^4 \bar\rho \vartheta_{\tau\tau}^2 \,dx \\
	& ~~~~~~ + \tilde\alpha^{1-a}\cdot \tilde\alpha^{1+a} \int x^4 \bar\rho\vartheta_\tau^2 \,dx + \tilde\alpha^{1-a} \cdot \tilde\alpha^{1+a} \int \chi \vartheta_\tau^2 \,dx.
\end{align*}
Such a fact, together with \eqref{estimates-004} , \eqref{estimates-101} , yields
\begin{equation}
	\int \mathcal G_{x\tau}^2 \,dx \lesssim \tilde\alpha^{d_3}\cdot \mathcal{\tilde E}_0,
\end{equation}
where
\begin{equation*}
	d_3 = \max \lbrace - 6, -2-\mathfrak a, -2-b, -5-a, - 4 + d_2 \rbrace = 1 - a .
\end{equation*}

\end{pf}

\subsection{Pointwise Bounds and Asymptotic Stability}\label{sec:pointwise-bounde-isentropic}

To this end, we have prepared enough estimates to close the a priori assumption {\eqref{le:supnorm}}. We shall employ the standard embedding theory to obtain the bound of the a prior bounded quantity $ \tilde \omega $. Indeed, it will be shown
\begin{equation}\label{estimates-201}
	\tilde \omega^2 \leq C \mathcal{\tilde E}(\tau) + C\mathcal{\tilde E}_0 \leq C \mathcal{\tilde E}_0,
\end{equation} 
with some constant $ C >0 $ when $ \tilde\omega < \varepsilon_0 < 1 $ small enough. Therefore,  when $ \mathcal{\tilde E}_0 $ small enough, \eqref{estimates-201} shows shat $ \tilde \omega \leq \varepsilon_0/2 $. Then by employing the local well-posedness theory of \eqref{eq:perturbationrbation=linearly} and continuous arguments, the asymptotic stability theory for the linearly expanding homogeneous solutions in the isentropic case is justified. This will finish the proof of Theorem \ref{theorem2}.

Now we will show \eqref{estimates-201}. After applying the Hardy's inequality and the Sobolev embedding theory, it holds,
\begin{align}
& \norm{x\vartheta_x}{L_x^\infty}^2 \lesssim \norm{x\vartheta_x}{L_x^2}^2 + \norm{\vartheta_x + x\vartheta_{xx}}{L_x^2}^2 \lesssim \norm{\vartheta_x}{L_x^2}^2 + \norm{x\vartheta_{xx}}{L_x^2}^2, \label{lm:602} 
\\
& \norm{\vartheta}{L_x^\infty}^2 \lesssim \norm{\vartheta}{L_x^2}^2 + \norm{\vartheta_x}{L_x^2}^2 \lesssim \norm{x\bar\rho^{1/3}\vartheta}{L_x^2}^2 + \norm{\vartheta_x}{L_x^2}^2 {\nonumber} \\
& ~~~~~~ \lesssim \norm{x^2\bar\rho^{1/2}\vartheta}{L_x^2}^2 + \norm{\vartheta_x}{L_x^2}^2,  \label{lm:603} 
\end{align}
where we have used the fact $ \bar\rho^{1/3}(x) \simeq (R_0-x) $ from \eqref{densitydistance}. On the other hand, applying the fundamental theory of calculus and the H\"older inequality, 
\begin{align*}
	& \norm{x\vartheta_{x\tau}}{L_x^\infty}^2 \lesssim \norm{x\vartheta_{x\tau}}{L_x^2(R_0/4,R_0/2)}^2 + \int_0^{R_0} \abs{(x^2\vartheta_{x\tau}^2)_x}{} \,dx \\
	& \lesssim \norm{x\vartheta_{x\tau}}{L_x^2(R_0/4,R_0/2)}^2 + \norm{x\vartheta_{x\tau}}{L_x^2} \bigl( \norm{\vartheta_{x\tau}}{L_x^2} + \norm{x\vartheta_{xx\tau}}{L_x^2} \bigr)\\
	& \lesssim (1+\tilde\omega) \bigl(\norm{x((1+\vartheta)x\vartheta_{x\tau}-x\vartheta_x\vartheta_\tau)}{L_x^2}^2 + \norm{x^2\bar\rho^{1/2}\vartheta_\tau}{L_x^2}^2 \bigr) \\
	& ~~~~ + \bigl(\norm{\chi^{1/2}x\vartheta_{x\tau}}{L_x^2} + \norm{x^2
	\vartheta_{x\tau}}{L_x^2} \bigr) \bigl( \norm{\vartheta_{x\tau}}{L_x^2} + \norm{x\vartheta_{xx\tau}}{L_x^2} \bigr).
\end{align*}
Here we have separated the norm of $ \norm{x\vartheta_{x\tau}}{L_x^2} $ in the interior and boundary subdomains and used the following form of the mean value theory
\begin{equation*} \exists x_0 \in (R_0/4, R_0/2),~ \text{such that} ~ \abs{x\vartheta_{x\tau}(x_0,\tau)}{} \leq \dfrac{4}{R_0} \norm{x\vartheta_{x\tau}}{L_x^2(R_0/4,R_0/2)}. \end{equation*}
Then employing \eqref{ine:102} yields
\begin{equation}\label{lm:604}
	\begin{aligned}
		& \norm{x\vartheta_{x\tau}}{L_x^\infty}^2 \lesssim (1+\tilde\omega)\bigl( \norm{x((1+\vartheta)x\vartheta_{x\tau}-x\vartheta_x\vartheta_\tau)}{L_x^2}^2 + \norm{x^2\bar\rho^{1/2}\vartheta_\tau}{L_x^2}^2 \bigr) \\
		& ~~ + (1+\tilde\omega) \bigl(\norm{\chi^{1/2}x\vartheta_{x\tau}}{L_x^2} + \norm{x^2\bar\rho^{1/2}\vartheta_\tau}{L_x^2} + \norm{x((1+\vartheta)x\vartheta_{x\tau}-x\vartheta_{x}\vartheta_\tau)}{L_x^2} \bigr) \\ & ~~~~ \times \bigl( \norm{\vartheta_{x\tau}}{L_x^2} + \norm{x\vartheta_{xx\tau}}{L_x^2} \bigr).
	\end{aligned}
\end{equation}
Similarly,
\begin{equation}\label{lm:605}
\begin{aligned}
& \norm{\vartheta_\tau}{L_x^\infty}^2 \lesssim 	\norm{x^2 \bar\rho^{1/2}\vartheta_\tau}{L_x^2}^2 + (1+\tilde\omega)\bigl( \norm{\chi^{1/2}\vartheta_\tau}{L_x^2} + \norm{x^2\bar\rho^{1/2}\vartheta_\tau}{L_x^2} \\
& ~~ + \norm{x((1+\vartheta)x\vartheta_{x\tau}-x\vartheta_{x}\vartheta_\tau)}{L_x^2}  \bigr)\norm{\vartheta_{x\tau}}{L_x^2}.
\end{aligned}
\end{equation}

Then we have the following lemma.
\begin{lm}\label{lm:isentropic-linearly-embedding}
Under the same assumption as in Lemma \ref{lm:isentropic-linearly-interior-estimates-3}, for $ \varepsilon_0 $ small enough in any fixed time $ \tau > 0 $, we have
\begin{equation}\label{lm:601}
\norm{x\vartheta_x}{L_x^\infty}^2, \norm{\vartheta}{L_x^\infty}^2, \norm{x\vartheta_{x\tau}}{L_x^\infty}^2,\norm{\vartheta_\tau}{L_x^\infty}^2 \leq C \mathcal{\tilde E}(\tau) + C \mathcal{\tilde E}_0.
\end{equation}
In particular, this will show the first inequality in \eqref{estimates-201}.
\end{lm}

\begin{pf}
We only need to show the right of \eqref{lm:602}, \eqref{lm:603}, \eqref{lm:604} and \eqref{lm:605} can be bounded by $ \mathcal{\tilde E}(\tau) $. From the definition of  the total energy functional in \eqref{isentropic-linearly-total-energy}, directly we have the following
\begin{align*}
	& ~~~~~~~~~~~~~ \norm{x\vartheta_x(\tau)}{L_x^\infty}^2 , \norm{\vartheta(\tau)}{L_x^\infty}^2\lesssim  \mathcal{\tilde E}(\tau) + \mathcal{\tilde E}_0,\\
	& \norm{x\vartheta_{x\tau}(\tau)}{L_x^\infty}^2, \norm{\vartheta_\tau}{L_x^\infty}^2 \lesssim (1+\tilde\omega) \tilde\alpha^{-1-a} \mathcal{\tilde E} \\
	& ~~~~~~~~ + (1+\tilde\omega) \tilde\alpha^{(-1-a)/2} \mathcal{\tilde E}^{1/2} \times \tilde\alpha^{(1-a)/2} \mathcal{\tilde E}_0^{1/2} \lesssim \tilde\alpha^{-a}(\mathcal{\tilde E}(\tau) +\mathcal{\tilde E}_0 ).
\end{align*}
This finishes the proof.
\end{pf}

The second part of the inequality in \eqref{estimates-201} is a direct consequence of Lemma \ref{lm:isentropic-linearly-basic-energy} , Lemma \ref{lm:isentropic-linearly-temporal-derivative}, Lemma \ref{lm:isentropic-linearly-interior-estimates-1}, Lemma \ref{lm:isentropic-linearly-interior-estimates-2} and Lemma \ref{lm:isentropic-linearly-interior-estimates-3}.
\begin{lm}
Under the same assumption as in Lemma \ref{lm:isentropic-linearly-embedding}, we have $ \mathcal{\tilde E}(\tau) + \mathcal{\tilde D}(\tau) \leq C \mathcal{\tilde E}_0 $ for some positive constant $ C $. 
\end{lm}


\section{Stability of the Linearly Expanding Homogeneous Solution for the Thermodynamic Model}\label{section:sb-linearly-thermodynamic}

In this section, we will study the linearly expanding homogeneous solution of the thermodynamic model, which is given in Section \ref{section:introduction_of_expanding_sol}. More precisely, we will study the stability of such solutions. Similarly as before, we work with the perturbation variables $ (\xi,\zeta) $ (defined in \eqref{def:pt-variable-thermo}) in the Lagrangian coordinates. That is, we will work in the space-time $ (x,\tau) $, where $ \tau $ is the linearly temporal variable. The equations satisfied by $ (\xi, \zeta) $ are given in \eqref{eq:RGS-ptb-ln-tm-vb} together with the boundary condition \eqref{BC:thermodynamic}. The initial condition is given as $ (\xi(\cdot,0),\xi_\tau(\cdot,0),\zeta(\cdot,0) )=( \xi_0(\cdot),\xi_1(\cdot),\zeta_0(\cdot) ) $. To start with, we denote the point-wise bounds of the perturbation variables as
\begin{equation}\label{lsthermo-supnorm}
	\hat \omega : = \sup_\tau \lbrace \norm{x\xi_x(\tau)}{L_x^\infty} , \norm{\xi(\tau)}{L_x^\infty}, \norm{x\xi_{x\tau}(\tau)}{L_x^\infty}, \norm{\xi_{\tau}(\tau)}{L_x^\infty}, \norm{(\zeta/\sigma)(\tau)}{L_x^\infty} \rbrace,
\end{equation}
where $ \sigma = R_0 - x $ denotes the distance to the boundary. As in Section \ref{section:sb-linearly-isentropic}, we analysis the perturbation variables $ (\xi,\zeta) $ in three steps. 

First, by assuming that $ \hat \omega $ is bounded by some small constant $ \hat\varepsilon_0 $, we perform weighted estimates of the perturbation in Section \ref{sec:energy-est-thermodynamic}. Again, we will establish the lower and higher order energy estimates with temporal weights. The temporal weight will track the growth of the energy functionals.

In Section \ref{sec:interior-est-thermodynamic}, we will explore some interior estimates. In particular, these will recover the regularity of the solutions near the coordinate center $ x = 0 $. Moreover, the estimates will establish the interchange of temporal weights between the temporal derivatives and the spatial derivatives. 

Notice, in contrast to the isentropic model, we will make use of the large enough growing rate of the expanding solution to manipulate the extra forcing terms. Therefore, while establishing the temporal weighted estimates, we shall track carefully the constant $ a_1 $ during the proof.

To summing up the analysis, in Section \ref{sec:pointwise-bounde-thermodynamic}, we show that $ \hat \omega $ can be indeed bounded by the total energy functional defined in \eqref{thermodynamic-linearly-total-energy} and therefore bounded by the initial energy \eqref{thermodynamic-initial-energy}.
Then by choosing the initial energy small enough, we will end up with a bounded $ \hat\omega $. 
Then a continuity argument will demonstrate the asymptotic stability of the linearly expanding homogeneous solution for the thermodynamic model. 

In the following, the polynomial of $ \hat\omega $ will be also denoted as $ \hat\omega $. The total energy and dissipation functionals are given by
\begin{align}
	& \mathcal{\hat E}(t) :=  e^{(1+r_1)a_1\tau} \int x^4 \bar\rho \xi_\tau^2 \,dx + e^{l_1a_1\tau}\int x^2\bar\rho\zeta^2 \,dx + \int x^4 \bar\rho \xi^2\,dx + \int x^4 \xi_x^2 \,dx + \int x^2 \xi^2 \,dx {\nonumber}\\
	& ~~~~ + e^{(r_2-2)a_1\tau} \int x^4 \bar\rho\xi_{\tau\tau}^2 \,dx + e^{(l_2-2)a_1\tau}\int x^2 \bar\rho \zeta_\tau^2 \,dx + e^{((l_1+l_2)/2 + 1)a_1\tau} \int x^2\zeta_x^2 \,dx  {\nonumber}\\
	& ~~~~  + e^{((r_1+r_2)/2 + 3/2)a_1\tau} \int x^2\lbrack(1+\xi)x\xi_{x\tau} - x\xi_x\xi_\tau\rbrack^2 \,dx + \int \chi(\xi^2 +x^2 \xi_x^2 )\,dx {\nonumber}\\
	& ~~~~ + e^{(r_2+2)a_1\tau} \int \chi (x^2\xi_{x\tau}^2 +\xi_\tau^2 )\,dx + e^{(r_3-2)a_1\tau}\int \chi x^2 \bar\rho \xi_{\tau\tau}^2 \,dx + \int \xi_x^2 \,dx + \int x^2 \xi_{xx}^2 \,dx {\nonumber} \\
	& ~~~~ + \int \zeta_x^2 \,dx + e^{(r_3+2)a_1\tau}\int (\xi_{x\tau}^2 + x^2 \xi_{xx\tau}^2) \,dx + \int x^2 \zeta_{xx}^2\,dx  , \label{thermodynamic-linearly-total-energy} \\
	& \mathcal{\hat D}(t) := \int a_1 e^{(1+r_1)a_1 \tau}\int x^4 \bar\rho \xi_\tau^2 \,dx\,d\tau + \int e^{(3+r_1)a_1\tau} \int x^2 \lbrack (1+\xi)x\xi_{x\tau} - x\xi_x\xi_\tau\rbrack^2\,dx\,d\tau {\nonumber} \\
	& ~~~~ + \int a_1 e^{l_1a_1\tau}\int x^2 \bar\rho \zeta^2 \,dx\,d\tau + \int e^{(2+l_1)a_1\tau}\int x^2 \zeta_x^2 \,dx\,d\tau + \int a_1 e^{(r_2-2)a_1\tau}\int x^4 \bar\rho \xi_{\tau\tau}^2 \,dx\,d\tau {\nonumber}\\
	& ~~~~ + \int e^{r_2a_1\tau}\int x^2 \lbrack(1+\xi)x\xi_{x\tau\tau} - x\xi_x\xi_{\tau\tau}\rbrack^2 \,dx\,d\tau + \int e^{l_2a_1\tau}\int x^2 \zeta_{x\tau}^2 \,dx\,d\tau {\nonumber}\\
	& ~~~~ + \int e^{(3+\mathfrak r)a_1\tau} \int \chi \bigl( x^2 \xi_{x\tau}^2 + \xi_\tau^2 \bigr)  \,dx\,d\tau + \int e^{r_3a_1\tau}\int \chi(x^2\xi_{x\tau\tau}^2 + \xi_{\tau\tau}^2 )\,dx\,d\tau , \label{thermodynamic-linearly-dissipation}
\end{align}
where $ r_1,l_1 $ satisfies the following
\begin{equation}\label{total-constraints}
	\begin{gathered}
		-1 < r_1 < 1, ~ r_1 - 3 \leq  l_1  < -2, \\
		r_2 \leq r_1 - 1, ~ l_2 + 2\leq 0, ~ 0 \leq r_2 - l_2 \leq 2,\\
		-3 < \mathfrak r \leq r_2 - 1, ~ r_3 \leq r_2 - 2, ~ 
		l_2 + 2 \geq 0. 
	\end{gathered}
\end{equation}
The initial energy is given by
\begin{equation}\label{thermodynamic-initial-energy}
	\begin{aligned}
		& \mathcal{\hat E}_0 = \mathcal{\hat E}_0(\xi,\zeta) : =  \int x^4 \bar\rho \xi_1^2 \, dx + \int x^2 \bar\rho \zeta_0^2\,dx + \int x^4 \bar\rho \xi_0^2\,dx + \int x^4 \xi_{0,x}^2 \,dx  + \int x^4 \bar\rho \xi_2^2 \,dx \\
		& ~~~~+ \int x^2 \bar\rho \zeta_1^2 \,dx + \int \chi (\xi_0^2 + x^2 \xi_{0,x}^2 )\,dx + \int \chi x^2 \bar\rho \xi_2^2 \,dx + \int (\xi_{0,x}^2 + x^2 \xi_{0,xx}^2) \,dx,
	\end{aligned}
\end{equation}
where $ \chi $ is the cut-off function defined before as
\begin{equation*}
	\chi(x) = \begin{cases}
		1 & 0 \leq x \leq  R_0/2, \\
		0 & 3 R_0 / 4 \leq x \leq R_0,
	\end{cases}
\end{equation*}
and $ -4 \leq \chi'(x) \leq 0 $ and $ \xi_2, \zeta_1 $ are the initial date corresponding to $ \xi_{\tau\tau}, \zeta_\tau$ defined by the equations \eqref{eq:RGS-ptb-ln-tm-vb}. That is
\begin{equation*}
	\begin{aligned}
		& \dfrac{1}{(1+\xi_0)^2} \bigl(a_0 x\bar\rho \xi_{2} + a_0a_1 x\bar\rho\xi_1\bigr) + \bigl\lbrack \dfrac{K\bar\rho(\zeta_0 + \bar\theta)}{(1+\xi_0)^2(1+\xi_0+x\xi_{0,x})} \bigr\rbrack_x - \dfrac{(K\bar\rho\bar\theta)_x}{(1+\xi_0)^4} \\
		& ~~~~~~ = \dfrac{4\mu}{3} a_0^3 \bigl( \dfrac{\xi_1 + x\xi_{1,x}}{1+\xi_0+x\xi_{0,x}} +  \dfrac{2\xi_1}{1+\xi_0} \bigr)_x, \\
		& 3K x^2 \bar\rho \zeta_1 + K \dfrac{\bar\rho(\zeta_0+\bar\theta)(x^3(1+\xi_0)^2\xi_1)_x}{(1+\xi_0)^2(1+\xi_0+x\xi_{0,x})} - a_0^2 \bigl\lbrack \dfrac{(1+\xi_0)^2}{1+\xi_0+x\xi_{0,x}}x^2\zeta_{0,x} \\
		& ~~~~~~~~~~ + \bigl(\dfrac{(1+\xi_0)^2}{1+\xi_0+x\xi_{0,x}} - 1\bigr) x^2 \bar\theta_x \bigr\rbrack_x = \dfrac{4\mu}{3} a_0 x^2(1+\xi_0)^2(1+\xi_0+x\xi_{0,x}) \\
		& ~~~~~~~~~~~~~~\times \lbrack \dfrac{\xi_1+x\xi_{1,x}}{1+\xi_0+x\xi_{0,x}} - \dfrac{\xi_1}{1+\xi_0}  \rbrack^2.
	\end{aligned}
\end{equation*}

\begin{rmk}
	Here we write down an example for all the parameters, which satisfies the constraints \eqref{total-constraints},
	\begin{equation}\label{example_of_indices}
		r_1 = 1/2, ~ r_2 = - 1/2, ~ l_1 = -5/2, ~ l_2 = -2, ~ \mathfrak r = -3/2, r_3 = -5/2.
	\end{equation}
\end{rmk}

\subsection{Energy Estimates}\label{sec:energy-est-thermodynamic}
We start with the $ L^2 $-estimate of \eqref{eq:RGS-ptb-ln-tm-vb}.

\begin{lm}\label{lemma:thermodynamic-basic-energy-est} Considering a smooth solution $ (\xi,\zeta) $ to \eqref{eq:RGS-ptb-ln-tm-vb} with the corresponding boundary conditions \eqref{BC:thermodynamic}, define the functionals
	\begin{equation}\label{lm:901}
		\begin{aligned}
			& \mathcal{\hat E}_{\xi,1} : = e^{(1+r_1)a_1\tau} \int x^4 \bar\rho\xi_\tau^2 \,dx, ~~ \mathcal{\hat E}_{\zeta,1} : = e^{l_1 a_1\tau} \int x^2 \bar\rho \zeta^2 \,dx,  \\
			& \mathcal{\hat D}_{\xi,1} = \mathcal{\hat D}_{\xi,11} + \mathcal{\hat D}_{\xi,12} := \int a_1 e^{(1+r_1)a_1\tau} \int x^4 \bar\rho \xi_\tau^2 \,dx\,d\tau \\
			& ~~~~~~ + \int e^{(3+r_1)a_1\tau} \int x^2 \lbrack (1+\xi)x\xi_{x\tau}-x\xi_x\xi_\tau\rbrack^2 \,dx\,d\tau, \\
			& \mathcal{\hat D}_{\zeta,1}=\mathcal{\hat D}_{\zeta,11} + \mathcal{\hat D}_{\zeta,12}: = \int a_1 e^{l_1 a_1\tau} \int x^2 \bar\rho \zeta^2\,dx\,d\tau + \int e^{(2+l_1)a_1\tau}\int x^2 \zeta_x^2 \,dx \,d\tau,
		\end{aligned}
	\end{equation}
	where $ r_1, l_1 $ are some constants satisfying the following constraints
	\begin{equation}\label{constrain=001}
		-1 < r_1 < 1, ~  l_1 - r_1 \geq - 3 , ~ l_1 + 2 < 0.
	\end{equation}
	Suppose $ a_1 > 0 $ is large enough and $ \hat\omega < \hat \varepsilon_0 $ is small enough. Then we have
	\begin{equation}\label{lm:902}
		\mathcal{\hat E}_{\xi,1} + \mathcal{\hat E}_{\zeta,1} + \mathcal{\hat D}_{\xi,1}  + \mathcal{\hat D}_{\zeta,1} \lesssim C(\hat\varepsilon_0,a_1,r_1,l_1) \mathcal{\hat E}_0.
	\end{equation}
	Also,
	\begin{equation}\label{Rene:ub=001}
		\norm{x^2 {\bar\rho}^{1/2} \xi }{\stnorm{\infty}{2}}^2 + \norm{x^2 \xi_x}{\stnorm{\infty}{2}}^2 + \norm{x\xi}{\stnorm{\infty}{2}}^2 \lesssim C(\hat\varepsilon_0, a_1,r_1,l_1) \mathcal{\hat E}_0.
	\end{equation}
	
\end{lm}
\begin{pf}
	Multiply \subeqref{eq:RGS-ptb-ln-tm-vb}{1} with $ \tilde\alpha^{r_1} x^3 (1+\xi)^2 \xi_\tau $ and integrate the resulting equation over the spatial variable. It yields the following identity.
	\begin{equation}\label{Rene:id=001}
		\begin{aligned}
			& \dfrac{d}{d\tau} \hat E_{\xi,1} + \hat D_{\xi,1} = \tilde\alpha^{r_1} \int K x^2 \bar\rho \zeta \cdot \bigl( \dfrac{2 \xi_\tau}{1+\xi} + \dfrac{\xi_\tau+x\xi_{x\tau}}{1+\xi+x\xi_x} \bigr) \,dx \\
			& ~~~~~~~~ + \tilde\alpha^{r_1} \int K x^2 \bar\rho\bar\theta \bigl(\dfrac{1}{1+\xi+x\xi_x}-\dfrac{1}{(1+\xi)^2}\bigr)\cdot(3\xi_\tau+x\xi_{x\tau}) \,dx \\
			& ~~~~~~~~ + \tilde\alpha^{r_1} \int 2K x^2 \bar\rho\bar\theta \bigl( \dfrac{x\xi_x}{(1+\xi)^3} + \dfrac{x\xi_x}{(1+\xi)(1+\xi+x\xi_x)} \bigr) \cdot \xi_\tau\,dx =: \hat I_1 + \hat I_2 + \hat I_3, 
		\end{aligned}
	\end{equation}
	where
	\begin{equation*}
		\begin{aligned}
			\hat E_{\xi,1} & := \dfrac{\tilde\alpha^{1+r_1}}{2}\int x^4 \bar\rho \xi_\tau^2\,dx ,\\
			\hat D_{\xi,1} & := \dfrac{(1-r_1)\tilde\alpha^{r_1}\tilde\alpha_\tau}{2}\int x^4 \bar\rho\xi_\tau^2\,dx + \dfrac{4\mu}{3} \tilde\alpha^{3+r_1} \int \dfrac{x^2\lbrack(1+\xi)x\xi_{x\tau} - x\xi_x \xi_\tau\rbrack^2}{1+\xi+x\xi_x}\,dx.
		\end{aligned}
	\end{equation*}
	\todo{17:04}
	In the meantime, multiply \subeqref{eq:RGS-ptb-ln-tm-vb}{2} with $ \tilde\alpha^{l_1} \zeta $ and integrate the resulting in the spatial variable. 
	\begin{equation}\label{Rene:id=002}
		\begin{aligned}
			& \dfrac{d}{d\tau} \hat E_{\zeta,1} + \hat D_{\zeta,1} = - \tilde\alpha^{2+l_1}\int x^2 \bar\theta_x \bigl( \dfrac{(1+\xi)^2}{1+\xi+x\xi_x} - 1\bigr) \cdot \zeta_x \,dx \\
			& ~~~~ - \tilde\alpha^{l_1} \int \dfrac{K\bar\rho(\zeta+\bar\theta)\lbrack x^3 (1+\xi)^2\xi_\tau\rbrack_x}{(1+\xi)^2 (1+\xi+x\xi_x)} \cdot \zeta \,dx  + \tilde\alpha^{1+l_1} \int x^2 (1+\xi)^2 (1+\xi+x\xi_x) \mathfrak{\hat F}(\xi) \cdot\zeta \,dx\\
			& ~~~~~~~~ =: \hat I_4 + \hat I_5 + \hat I_6,
		\end{aligned}
	\end{equation}
	where
	\begin{equation*}
		\begin{aligned}
			\hat E_{\zeta,1} & : = \dfrac{3K \tilde\alpha^{l_1}}{2} \int x^2 \bar\rho \zeta^2 \,dx , \\
			\hat D_{\zeta,1} & : = - l_1 \dfrac{3K \tilde\alpha^{l_1-1} \tilde\alpha_\tau }{2} \int x^2 \bar\rho \zeta^2 \,dx + \tilde\alpha^{2+l_1} \int \dfrac{(1+\xi)^2}{1+\xi+x\xi_x} \cdot x^2 \zeta_x^2\,dx.
		\end{aligned}
	\end{equation*}
	
	\todo{17:07}
	By choosing $ 1-r_1 > 0, - l_1 > 0  $, the following estimates on the energy and dissipation functionals $ \hat E_{\xi,1}, \hat E_{\zeta,1},  \hat D_{\xi,1}, \hat D_{\zeta,1} $ hold.
	\begin{align*}
		& \hat E_{\xi,1} \gtrsim \tilde\alpha^{1+r_1} \int x^4 \bar\rho \xi_\tau^2 \,dx \gtrsim e^{(1+r_1)a_1\tau} \int x^4 \bar\rho \xi_\tau^2 \,dx ,\\
		& \hat E_{\zeta,1} \gtrsim \tilde\alpha^{l_1} \int x^2 \bar\rho \zeta^2\,dx \gtrsim e^{l_1 a_1 \tau} \int x^2\bar\rho \zeta^2\,dx, \\
		& \hat D_{\xi,1} \gtrsim (1-r_1)\tilde\alpha^{r_1}\tilde\alpha_\tau \int x^4 \bar\rho \xi_\tau^2 \,dx + (1-\hat \omega) \tilde\alpha^{3+r_1} \int x^2 \lbrack (1+\xi)x\xi_{x\tau} - x\xi_x\xi_\tau\rbrack^2 \,dx {\nonumber}\\
		& ~~~~~~ \gtrsim (1-r_1) a_1 e^{(r_1+1)a_1\tau} \int x^4 \bar\rho \xi_\tau^2 \,dx + (1-\hat\omega) e^{(3+r_1)a_1\tau} \int x^2 \lbrack (1+\xi)x\xi_{x\tau} - x\xi_x\xi_\tau\rbrack^2 \,dx ,\\
		& \hat D_{\zeta,1} \gtrsim - l_1 \tilde\alpha^{l_1-1} \tilde\alpha_\tau \int x^2 \bar\rho \zeta^2\,dx + (1-\hat\omega) \tilde\alpha^{2+l_1}\int x^2 \zeta_x^2 \,dx \\
		& ~~~~~~ \gtrsim - l_1 a_1 e^{l_1 a_1 \tau}\int x^2 \bar\rho \zeta^2 \,dx +(1-\hat \omega) e^{(2+l_1) a_1 \tau} \int x^2 \zeta_x^2 \,dx.
	\end{align*}
	Integration in the temporal variable of \eqref{Rene:id=001} and \eqref{Rene:id=002} yields the following.
	\begin{equation}\label{Rene:ineq=001}
		\begin{gathered}
			\mathcal{\hat E}_{\xi,1} + (1-r_1)\mathcal{\hat D}_{\xi, 11} + \mathcal{\hat D}_{\xi,12}
			\lesssim \int \hat I_1 \,d\tau + \int \hat I_2 \,d\tau + \int \hat I_3 \,d\tau + \mathcal{\hat E}_0 , \\
			\mathcal{\hat E}_{\zeta,1} + (-l_1) \mathcal{\hat D}_{\zeta,11} + \mathcal{\hat D}_{\zeta,12}
			\lesssim \int \hat I_4\,d\tau + \int \hat I_5 \,d\tau + \int \hat I_6 \,d\tau + \mathcal{\hat E}_0 . 
		\end{gathered}
	\end{equation}
	What is left is to perform estimates on the right of \eqref{Rene:ineq=001}. In order to do so, we shall need the following estimates similar to \eqref{ine:003} and \eqref{ine:006}. 
	\begin{equation}\label{Rene:ineq=002}
		\begin{aligned}
			& \int x^4 \bar\rho \xi^2 \,dx \lesssim \int x^4 \bar\rho \xi_0^2\,dx + \biggl\lbrack \int \bigl( \int x^4 \bar\rho \xi_\tau^2 \,dx \bigr)^{1/2} \,d\tau \biggr\rbrack^2  \\
			& ~~~~~~~~~~ \lesssim \int x^4 \bar\rho \xi_0^2 \,dx + \int a_1^{-1} e^{-(1+r_1)a_1\tau} \,d\tau \cdot \mathcal{\hat D}_{\xi,11}, \\ 
			& \int x^4 \xi_x^2 \,dx \lesssim \int x^4 \xi_{0,x}^2 \,dx + \biggl\lbrack \int \bigl( \int x^2\lbrack(1+\xi)x\xi_{x\tau} - x\xi_x \xi_\tau\rbrack^2\,dx  \bigr)^{1/2} \,d\tau \biggr\rbrack^2   \\
			& ~~~~~~~~~~ \lesssim \int x^4 \xi_{0,x}^2 \,dx + \int e^{-(3+r_1)a_1\tau} \,d\tau \cdot \mathcal{\hat D}_{\xi,12}. 
		\end{aligned}
	\end{equation}
	In the meantime, the right of \eqref{Rene:ineq=001} can be estimated as follows.
	\begin{align}
		& \int \hat I_1 \,d\tau \lesssim (1+\hat\omega)\int e^{r_1 a_1 \tau} \int x^2 \bar\rho |\zeta| ( |\xi_\tau|+|x\xi_{x\tau}|) \,dx \,d\tau \lesssim \epsilon \mathcal{\hat D}_{\xi,1} {\nonumber} \\
		& ~~~~ + C_\epsilon(1+\hat\omega) a_1^{-1} \sup_\tau e^{(r_1-l_1-3)a_1 \tau} \cdot\mathcal{\hat D}_{\zeta,12} + C_\epsilon(1+\hat\omega)a_1^{-1} \sup_\tau e^{(r_1-l_1-3)a_1\tau} \cdot \mathcal{\hat D}_{\zeta,11},
		{\label{Rene:ineq=003}}\\
		& \int \hat I_2 \,d\tau + \int \hat I_3 \,d\tau \lesssim (1+\hat\omega) \int e^{r_1 a_1 \tau} \int x^2 \bar\rho\bar\theta (|\xi| + |x\xi_x|)( |\xi_\tau| + |x \xi_{x\tau}|) \,dx \,d\tau \lesssim \epsilon \mathcal{\hat D}_{\xi,1} {\nonumber} \\
		& ~~~~ + C_\epsilon(1+\hat\omega) \bigl(a_1^{-1}\int e^{(r_1-1)a_1\tau}\,d\tau + \int e^{(r_1-3)a_1\tau}\,d\tau \bigr) \cdot\sup_\tau \biggl\lbrace \int x^4\bar\rho\xi^2 \,dx + \int x^4 \xi_x^2 \,dx \biggr\rbrace,
		\\
		& \int \hat I_4 \,d\tau \lesssim (1+\hat\omega) \int  e^{(2+l_1)a_1\tau} \int x^2 |\bar\theta_x| (|\xi| + |x\xi_x|) |\zeta_x| \,dx \lesssim  \epsilon \mathcal{\hat D}_{\zeta,12} {\nonumber} \\ & ~~~~ + C_\epsilon(1+\hat\omega) \int e^{(2+l_1)a_1\tau} \,d\tau\cdot \sup_\tau \biggl\lbrace \int x^4\bar\rho\xi^2\,dx + \int x^4 \xi_x^2 \,dx \biggr\rbrace, \\
		& \int \hat I_5 \,d\tau \lesssim (1+\hat\omega) \int e^{l_1a_1\tau} \int x^2 \bar\rho (|\zeta| + |\bar\theta|) (|\xi_\tau| + |x\xi_{x\tau}|) \cdot |\zeta| \,dx\,d\tau \lesssim  \epsilon \mathcal{\hat D}_{\zeta, 1} {\nonumber}\\
		& ~~~~ + C_\epsilon(1+\hat\omega)a_1^{-1} \sup_\tau e^{(l_1-r_1 - 3)a_1\tau} \cdot \mathcal{\hat D}_{\xi,11} + C_\epsilon(1+\hat\omega)a_1^{-1} \sup_\tau e^{(l_1-r_1-3)a_1\tau} \cdot \mathcal{\hat D}_{\xi,12}, \\
		& \int \hat I_6 \,d\tau \lesssim (1+\hat\omega) \int e^{(1+l_1)a_1\tau} \int x^2 (\xi_\tau^2 + x^2 \xi_{x\tau}^2) \cdot\zeta\,dx\,d\tau \lesssim \epsilon \mathcal{\hat D}_{\zeta,12} {\nonumber}\\
		& ~~~~ + C_\epsilon \hat\omega a_1^{-1} \sup_\tau e^{(l_1-r_1-1)a_1\tau} \cdot \mathcal{\hat D}_{\xi,11} + C_\epsilon \hat\omega\sup_\tau e^{(l_1-r_1-3)a_1\tau} \cdot\mathcal{\hat D}_{\xi, 12} .
		{\label{Rene:ineq=004}}
	\end{align}
	where the following inequality is made use of,
	\begin{gather}
		\int x^2 \xi_\tau^2 \,dx  + \int x^4 \xi_{x\tau}^2 \,dx \lesssim \int x^4 \bar\rho \xi_\tau^2 \,dx + \int x^2\lbrack (1+\xi)x\xi_{x\tau} - x\xi_x\xi_\tau\rbrack^2 \,dx  , \label{Rene:ineq=005} \\
		\int x^2 \xi^2 \,dx \lesssim \int x^4 \bar\rho \xi^2\,dx + \int x^4 \xi_x^2 \,dx, ~  \int x^2 \zeta^2 \,dx \lesssim \int x^2 \zeta_x^2 \,dx. {\nonumber}
	\end{gather}
	Therefore, for $ r_1, l_1 $ satisfying
	\begin{equation*}\tag{\ref{constrain=001}}
		-1 < r_1 < 1, ~  l_1 - r_1 \geq - 3 , ~ l_1 + 2 < 0,
	\end{equation*}
	\eqref{Rene:ineq=001} and \eqref{Rene:ineq=002} together with the inequalities from \eqref{Rene:ineq=003} to  \eqref{Rene:ineq=004} imply the following energy estimates,
	\begin{align*}
		& \mathcal{\hat E}_{\xi,1} + \mathcal{\hat D}_{\xi,1} \lesssim \epsilon \mathcal{\hat D}_{\xi,1} + C_\epsilon (1+\hat\omega)a_1^{-1} \mathcal{\hat D}_{\zeta,1} {\nonumber} \\
		& ~~~~ + C_\epsilon (1+\hat\omega)(a_1^{-2} + a_1^{-1} ) \times (\mathcal{\hat E}_0 + a_1^{-2} \mathcal{\hat D}_{\xi,11} + a_1^{-1}\mathcal{\hat D}_{\xi,12}) + \mathcal{\hat E}_0, 	\\
		& \mathcal{\hat E}_{\zeta,1} + \mathcal{\hat D}_{\zeta,1} \lesssim \epsilon \mathcal {\hat D}_{\zeta,1} + C_\epsilon(1+\hat\omega)a_1^{-1} \times (\mathcal{\hat E}_0 + a_1^{-2}\mathcal{\hat D}_{\xi,11} + a_1^{-1} \mathcal{\hat D}_{\xi,12} ) {\nonumber}	\\
		& ~~~~ + C_\epsilon(1+\hat\omega) (a_1^{-1} + \hat\omega) \mathcal{\hat D}_{\xi,1} + \mathcal{\hat E}_0. 
	\end{align*}
	Consequently, after choosing $ \epsilon > 0 $ small enough, for sufficiently large $ a_1 > 0 $ and sufficiently small $  \hat\varepsilon_0 > 0 $, the above estimates yield,
	\begin{equation*}\tag{\ref{lm:902}}
		\mathcal{\hat E}_{\xi,1} + \mathcal{\hat E}_{\zeta,1} + \mathcal{\hat D}_{\xi,1}  + \mathcal{\hat D}_{\zeta,1} \lesssim C(\hat\varepsilon_0,a_1,r_l,l_1) \mathcal{\hat E}_0.
	\end{equation*}
	In particular, from \eqref{Rene:ineq=002}, \eqref{Rene:ineq=005},
	\begin{equation*}\tag{\ref{Rene:ub=001}}
		\norm{x^2 \sqrt{\bar\rho} \xi }{\stnorm{\infty}{2}}^2 + \norm{x^2 \xi_x}{\stnorm{\infty}{2}}^2 + \norm{x\xi}{\stnorm{\infty}{2}}^2 \lesssim C(\hat\varepsilon_0, a_1,r_l,l_1) \mathcal{\hat E}_0.
	\end{equation*}
\end{pf}

\todo{19:40, 13 Oct, 2017}

Next, we shall write down the corresponding temporal derivative version of the system \eqref{eq:RGS-ptb-ln-tm-vb}. Direct calculation $ \partial_\tau \lbrack  \tilde\alpha^{-3} (1+\xi)^2 \cdot \eqref{eq:RGS-ptb-ln-tm-vb}_1 \rbrack $ and $ \partial_\tau \lbrack \tilde\alpha^{-2} \cdot \eqref{eq:RGS-ptb-ln-tm-vb}_2  \rbrack $ yields the following system
\begin{equation}\label{eq:1st-RGS-ptb-ln-tm-vb}
	\begin{cases}
		\tilde\alpha^{-2} x \bar\rho \xi_{\tau\tau\tau} - \tilde\alpha^{-3} \tilde\alpha_{\tau} x \bar\rho \xi_{\tau\tau} + (\tilde\alpha^{-3}\tilde\alpha_{\tau\tau} - 3 \tilde\alpha^{-4} \tilde\alpha_\tau^2 ) x \bar\rho \xi_\tau\\
		~~~~~~~~ +\bigl\lbrace \tilde\alpha^{-3}(1+\xi)^2 \bigl\lbrace \bigl\lbrack \dfrac{K\bar\rho(\zeta + \bar\theta)}{(1+\xi)^2(1+\xi+x\xi_x)} \bigr\rbrack_x - \dfrac{(K\bar\rho\bar\theta)_x}{(1+\xi)^4} \bigr\rbrace  \bigr\rbrace_\tau \\
		~~~~~~ = (1+\xi)^2 \bigl( \mathfrak{\hat B}_{x\tau} + 4\mu \bigl( \dfrac{\xi_\tau}{1+\xi}\bigr)_{x\tau} \bigr) + 2 (1+\xi) \xi_\tau \bigl( \mathfrak{\hat B}_{x} + 4\mu \bigl( \dfrac{\xi_{\tau}}{1+\xi}\bigr)_x \bigr),\\
		3K\tilde\alpha^{-2} x^2 \bar\rho \zeta_{\tau\tau} - 6 K \tilde\alpha^{-3} \tilde\alpha_\tau x^2 \bar\rho \zeta_\tau + \bigl\lbrack K\tilde\alpha^{-2} \dfrac{\bar\rho(\zeta+\bar\theta) \lbrack x^3(1+\xi)^2\xi_\tau\rbrack_x}{(1+\xi)^2(1+\xi+x\xi_x)} \bigr\rbrack_\tau \\
		~~~~~~~~ - \bigl\lbrack \dfrac{(1+\xi)^2}{1+\xi+x\xi_x}x^2\zeta_x + \bigl(\dfrac{(1+\xi)^2}{1+\xi+x\xi_x} - 1\bigr) x^2 \bar\theta_x \bigr\rbrack_{x\tau} \\
		~~~~~~ = \lbrack \tilde\alpha^{-1} x^2 (1+\xi)^2(1+\xi+x\xi_x)\cdot \mathfrak{\hat F}(\xi) \rbrack_\tau,
	\end{cases}
\end{equation}
with
\begin{equation*}
	\mathfrak{\hat B}_\tau(R_0,\tau) = 0, ~ \zeta_\tau(R_0,\tau) = 0.
\end{equation*}
Similarly as before, we shall perform the $ L^2 $-based energy estimate in the following.

\begin{lm}\label{lemma:thermodynamic-high-order-est}
	Under the assumption as in Lemma \ref{lemma:thermodynamic-basic-energy-est}, define the functionals
	\begin{equation}\label{lm:1001}
		\begin{aligned}
			& \mathcal{\hat E}_{\xi,2} : = e^{(r_2-2)a_1\tau}\int x^4 \bar\rho \xi_{\tau\tau}^2 \,dx, ~~ \mathcal{\hat E}_{\zeta,2} : = e^{(l_2-2)a_1\tau} \int x^2 \bar\rho\zeta_\tau^2 \,dx , \\
			& \mathcal{\hat D}_{\xi,2} = \mathcal{\hat D}_{\xi,21} + \mathcal{\hat D}_{\xi, 22} : = \int a_1e^{(r_2-2)a_1\tau} \int x^4 \bar\rho\xi_{\tau\tau}^2 \,dx\,d\tau \\
			& ~~~~~~ + \int e^{r_2a_1\tau} \int x^2 \lbrack(1+\xi)x\xi_{x\tau\tau}-x\xi_x\xi_{\tau\tau}\rbrack^2 \,dx\,d\tau ,\\
			& \mathcal{\hat D}_{\zeta,2} : = \int e^{l_2a_1\tau}\int x^2 \zeta_{x\tau}^2 \,dx\,d\tau,
		\end{aligned}	
	\end{equation}
	where $ r_2, l_2 $ satisfy the constraints 
	\begin{equation}\label{constrain=002}
		r_2 \leq r_1 - 1, ~ l_2 + 2 \leq 0, ~ 0 \leq r_2 - l_2 \leq 2.
	\end{equation}
	For $ a_1 $ large enough and $ \hat \varepsilon_0 $ small enough, we have the following estimates,
	\begin{equation}\label{lm:1002}
		\mathcal{\hat E}_{\xi,2} + \mathcal{\hat E}_{\zeta,2} + \mathcal{\hat D}_{\xi,2}  + \mathcal{\hat D}_{\zeta,2} \lesssim C(\hat\varepsilon_0,a_1,r_1,l_1,r_2,l_2)  \mathcal{\hat E}_0.
	\end{equation}
\end{lm}

\begin{pf}
	Multiply \subeqref{eq:1st-RGS-ptb-ln-tm-vb}{1} with $ \tilde\alpha^{r_2} x^3 \xi_{\tau\tau} $ and integrate the resulting in the spatial variable. It holds, 
	\begin{equation}\label{Rene:id=003}
		\begin{aligned}
			& \dfrac{d}{d\tau}\hat E_{\xi,2} + \hat D_{\xi, 2} = \hat L_1 + \hat L_2 + \hat L_3 + \hat L_4,
		\end{aligned}
	\end{equation}
	where
	\begin{align*}
		& \hat E_{\xi,2} : = \dfrac{\tilde\alpha^{r_2-2}}{2} \int x^4 \bar\rho\xi_{\tau\tau}^2\,dx, \\
		& \hat D_{\xi,2} : = - \dfrac{r_2}{2} \tilde\alpha^{r_2-3}\tilde\alpha_{\tau} \int x^4 \bar\rho \xi_{\tau\tau}^2 \,dx + + \dfrac{4}{3}\mu \tilde\alpha^{r_2} \int \dfrac{x^2 \lbrack (1+\xi)x\xi_{x\tau\tau} - x\xi_x \xi_{\tau\tau}\rbrack^2 }{1+\xi+x\xi_x} \,dx, \\
		& \hat L_1 : = - (\tilde\alpha^{r_2-3}\tilde\alpha_{\tau\tau} - 3 \tilde\alpha^{r_2-4}\tilde\alpha_\tau^2 )\int x^4 \bar\rho \xi_\tau \xi_{\tau\tau} \,dx, \\
		& \hat L_2 : = - \dfrac{4}{3}\mu \tilde\alpha^{r_2} \int x^2 \biggl\lbrace - \dfrac{(1+\xi)^2(\xi_\tau+x\xi_{x\tau})^2(\xi_{\tau\tau}+x\xi_{x\tau\tau})}{(1+\xi+x\xi_x)^2}  {\nonumber} \\
		& ~~~~~~ + 2 \dfrac{(1+\xi)\xi_\tau(\xi_\tau+x\xi_{x\tau})(\xi_{\tau\tau}+x\xi_{x\tau\tau})}{1+\xi+x\xi_x} - \xi_\tau^2(\xi_{\tau\tau}+x\xi_{x\tau\tau})
		\biggr\rbrace \,dx, \\
		& \hat L_3 : = - 3 \tilde\alpha^{r_2-4}\tilde\alpha_\tau \int x^2 \biggl\lbrace \dfrac{2K\bar\rho\zeta x \xi_x \xi_{\tau\tau}}{(1+\xi)(1+\xi+x\xi_x)} + \dfrac{2K\bar\rho\bar\theta x\xi_x \xi_{\tau\tau}}{(1+\xi)(1+\xi+x\xi_x)}+ \dfrac{2K\bar\rho\bar\theta x\xi_x \xi_{\tau\tau}}{(1+\xi)^3}   {\nonumber}\\
		& ~~~~~~ + \dfrac{K\bar\rho \zeta (3\xi_{\tau\tau}+x\xi_{x\tau\tau})}{1+\xi+x\xi_x} + \dfrac{K\bar\rho\bar\theta(3\xi_{\tau\tau} + x \xi_{x\tau\tau})}{1+\xi+x\xi_x} - \dfrac{K\bar\rho\bar\theta(3\xi_{\tau\tau} + x \xi_{x\tau\tau})}{(1+\xi)^2} \biggr\rbrace\,dx , \\
		& \hat L_{4}: = \tilde\alpha^{r_2-3} \int x^2 \biggl\lbrace (3 \xi_{\tau\tau} + x\xi_{x\tau\tau}) \bigl\lbrack \dfrac{K\bar\rho \zeta}{1+\xi+x\xi_x} + \dfrac{K\bar\rho\bar\theta}{1+\xi+x\xi_x} - \dfrac{K\bar\rho\bar\theta}{(1+\xi)^2} \bigr\rbrack_\tau {\nonumber} \\
		& ~~~~~~ + \xi_{\tau\tau} \bigl\lbrack \dfrac{2K\bar\rho \zeta x\xi_x}{(1+\xi)(1+\xi+x\xi_x)} + \dfrac{2K\bar\rho\bar\theta x\xi_x }{(1+\xi)(1+\xi+x\xi_x)}  + \dfrac{2K\bar\rho\bar\theta x\xi_x}{(1+\xi)^3} \bigr\rbrack_{\tau} \biggr\rbrace \, dx.
	\end{align*}
	In the meantime, multiply \subeqref{eq:1st-RGS-ptb-ln-tm-vb}{2} with $ \tilde\alpha^{l_2} \zeta_\tau $ and integrate the resulting in the spatial variable. It holds,
	\todo{10:30} 
	\begin{equation}\label{Rene:id=004}
		\dfrac{d}{d\tau} \hat E_{\zeta,2} + \hat D_{\zeta,2} = \hat L_5 + \hat L_6 + \hat L_7 + \hat L_8,
	\end{equation}
	where
	\begin{align*}
		& \hat E_{\zeta,2} := \dfrac{3K}{2} \tilde\alpha^{l_2-2} \int x^2 \bar\rho \zeta_\tau^2 \,dx, \\
		& \hat D_{\zeta,2} := - 3K (l_2/2 + 1) \tilde\alpha^{l_2-3}\tilde\alpha_\tau\int x^2 \bar\rho\zeta_\tau^2 \,dx  + \tilde\alpha^{l_2} \int \dfrac{(1+\xi)^2}{1+\xi+x\xi_x} \cdot x^2 \zeta_{x\tau}^2 \,dx  , \\
		& \hat L_5 := - \tilde\alpha^{l_2} \int x^2 \zeta_{x\tau}\biggl\lbrace  \zeta_x \cdot \bigl\lbrack \dfrac{(1+\xi)^2}{1+\xi+x\xi_x} \bigr\rbrack_\tau + \bigl( \dfrac{(1+\xi)^2}{1+\xi+x\xi_x}-1\bigr)_\tau  \bar\theta_x \biggr\rbrace \,dx, \\
		& \hat L_6 : = 2 K  \tilde\alpha^{l_2-3} \tilde\alpha_\tau  \int \zeta_\tau \cdot \dfrac{\bar\rho (\zeta +\bar\theta ) \lbrack x^3(1+\xi)^2\xi_\tau\rbrack_x}{(1+\xi)^2(1+\xi+x\xi_x)} \,dx, \\
		& \hat L_7 : = - K \tilde\alpha^{l_2-2} \int \zeta_\tau \cdot \biggl\lbrace \dfrac{\bar\rho(\zeta+\bar\theta) \lbrack x^3(1+\xi)^2\xi_\tau\rbrack_x}{(1+\xi)^2(1+\xi+x\xi_x)} \biggr\rbrace_\tau \, dx , \\
		& \hat L_8 := -  \tilde\alpha^{l_2-2}\tilde\alpha_\tau \int x^2 (1+\xi)^2(1+\xi+x\xi_x) \mathfrak{\hat F}(\xi) \zeta_\tau \,dx {\nonumber}\\
		& ~~~~~~  + \tilde\alpha^{l_2-1} \int x^2 \lbrack (1+\xi)^2(1+\xi+x\xi_x) \rbrack_\tau \mathfrak{\hat F}(\xi)\cdot\zeta_\tau \,dx {\nonumber}\\
		& ~~~~~~ + \tilde\alpha^{l_2-1} \int x^2(1+\xi)^2(1+\xi+x\xi_x)\lbrack\mathfrak{\hat F}(\xi)\rbrack_\tau \cdot \zeta_\tau \,dx.
	\end{align*}
	Integration in the temporal variable of \eqref{Rene:id=003} and \eqref{Rene:id=004} yields the following inequalities, provided $ r_2 < 0, l_2+2 \leq 0 $, 
	\begin{equation}\label{Rene:ineq=101}
		\begin{aligned}
			& \mathcal{\hat E}_{\xi,2} + (-r_2) \mathcal{\hat D}_{\xi,21} + (1-\hat\omega) \mathcal{\hat D}_{\xi,22}
			\lesssim \int \hat L_1 \,d\tau + \int \hat L_2 \,d\tau + \int \hat L_3\,d\tau + \int \hat L_4\,d\tau + \mathcal{\hat E}_0, \\
			& \mathcal{\hat E}_{\zeta,2} + (- l_2/2 - 1)\int a_1 e^{(l_2-2) a_1 \tau} \int x^2 \bar\rho\zeta_\tau^2 \,dx \,d\tau + (1-\hat\omega) \mathcal{\hat D}_{\zeta,2} 
			\\
			& ~~~~ \lesssim \int\hat L_5\,d\tau + \int\hat L_6\,d\tau + \int\hat L_7\,d\tau + \int\hat L_8\,d\tau + \mathcal{\hat E}_0 . 
		\end{aligned}
	\end{equation}
	\todo{11:35}
	Similarly as before, we shall establish the following estimates concerning the right of \eqref{Rene:ineq=101}.
	\begin{align*}
		& \int \hat L_1 \,d\tau 
		\lesssim \epsilon \mathcal{\hat D}_{\xi,21} + C_\epsilon a_1^2 \sup_\tau e^{(r_2-r_1-3)a_1\tau} \cdot \mathcal{\hat D}_{\xi,11} , \\
		& \int \hat L_2 \,d\tau 
		\lesssim \epsilon \mathcal{\hat D}_{\xi, 2} {\nonumber} 
		+ C_\epsilon \hat \omega ( a_1^{-2} \sup_\tau e^{(r_2 - r_1 + 1)a_1\tau} + a_1^{-1}\sup_\tau e^{(r_2-r_1-1)a_1\tau} ) \cdot \mathcal{\hat D}_{\xi,11} {\nonumber}\\
		& ~~~~ ~~~~ + C_\epsilon \hat \omega ( a_1^{-1} \sup_\tau e^{(r_2-r_1-1)a_1\tau} + \sup_\tau e^{(r_2-r_1-3)a_1\tau} ) \cdot \mathcal{\hat D}_{\xi,12} , \\
		& \int \hat L_3 \,d\tau \lesssim (1+\hat\omega) \int a_1e^{(r_2-3)a_1\tau}\int x^2 \bar\rho (|\zeta| + (|\xi| + |x\xi_{x}|) ) (|\xi_{\tau\tau}| + |x\xi_{x\tau\tau}|) \,dx\,d\tau \lesssim \epsilon \mathcal{\hat D}_{\xi,2}   {\nonumber} \\
		& ~~~~ ~~~~ + C_\epsilon(1+\hat\omega) (a_1 \sup_\tau e^{(r_2-l_1-6)a_1\tau} + a_1^2 \sup_\tau e^{(r_2-l_1-8)a_1\tau}) \cdot \mathcal{\hat D}_{\zeta,12} {\nonumber}\\
		& ~~~~ ~~~~ + C_\epsilon(1+\hat\omega) \int(a_1e^{(r_2-4)a_1\tau} + a_1^2 e^{(r_2-6)a_1\tau}) \,d\tau  \times ( \int x^4 \bar\rho \xi^2 \,dx + \int x^4 \xi_x^2\,dx ), \\
		& \int \hat L_4 \,d\tau \lesssim (1+\hat\omega) \int e^{(r_2-3)a_1\tau} \int x^2 \bar\rho ( |\zeta_\tau| + |\xi_\tau| + |x\xi_{x\tau} )(|\xi_{\tau\tau}| + |x\xi_{x\tau\tau}|) \,dx\,d\tau \lesssim \epsilon \mathcal{\hat D}_{\xi,2}  {\nonumber} \\
		& ~~~~ + C_\epsilon (1+\hat\omega) \int (a_1^{-1} e^{(r_2-l_2-2)a_1\tau} + e^{(r_2-l_2-4)a_1\tau}) \,d\tau \cdot \mathcal{\hat E}_{\zeta,2} {\nonumber}\\
		& ~~~~ + C_\epsilon (1+\hat \omega)  \sup_\tau ( a_1^{-2} e^{(r_2-r_1-5)a_1\tau} {\nonumber}
		+ a_1^{-1} e^{(r_2-r_1-7)a_1\tau} ) \cdot \mathcal{\hat D}_{\xi,11} {\nonumber} \\
		& ~~~~ + C_\epsilon (1+\hat\omega) ( a_1^{-1} \sup_\tau e^{(r_2-r_1-7)a_1\tau} + \sup_\tau e^{(r_2-r_1-9)a_1\tau} ) \cdot \mathcal{\hat D}_{\xi,12}, \\
		& \int \hat L_5 \,d\tau \lesssim (1+\hat\omega) \int e^{l_2a_1\tau} \int x^2 |\zeta_{x\tau}| \cdot (|\zeta_x| + |\bar\theta_x|) (|\xi_\tau| + |x\xi_{x\tau}|) \,dx\,d\tau \lesssim \epsilon \mathcal{\hat D}_{\zeta,2} {\nonumber}\\
		& ~~~~ + C_\epsilon \hat\omega \sup_\tau e^{(l_2-l_1-2)a_1\tau} \cdot \mathcal{\hat D}_{\zeta,12} + C_\epsilon (1+\hat\omega) a_1^{-1} \sup_\tau e^{(l_2-r_1-1)a_1\tau} \cdot \mathcal{\hat D}_{\xi,11} {\nonumber}\\
		& ~~~~ + C_\epsilon (1+\hat\omega) \sup_\tau e^{(l_2-r_1-3)a_1\tau}\cdot \mathcal{\hat D}_{\xi,12},\\
		& \int \hat L_6 \,d\tau \lesssim (1+\hat\omega) \int a_1 e^{(l_2-2)a_1\tau} \int x^2 |\zeta_\tau|(|\zeta| + |\bar\theta| ) \cdot (|\xi_\tau| + |x\xi_{x\tau}|) \,dx\,d\tau \lesssim \epsilon \mathcal{\hat D}_{\zeta,2} {\nonumber} \\
		& ~~~~ + C_\epsilon (1+\hat\omega)\hat\omega a_1^2 \sup_\tau e^{(l_2-l_1-6)a_1\tau} \cdot \mathcal{\hat D}_{\zeta,12} + C_\epsilon (1+\hat\omega) a_1 \sup_\tau e^{(l_2-r_1-5)a_1\tau} \cdot \mathcal{\hat D}_{\xi,11} {\nonumber}\\
		& ~~~~ + C_\epsilon (1+\hat\omega) a_1^2 \sup_\tau e^{(l_2-r_1-7)a_1\tau} \cdot \mathcal{\hat D}_{\xi,12}, \\
		& \int \hat L_7 \,d\tau \lesssim (1+\hat\omega) \int e^{(l_2-2)a_1\tau} \int x^2 \bar\rho |\zeta_\tau| \biggl\lbrack |\xi_\tau|^2 + |x\xi_{x\tau}|^2  + |\xi_{\tau\tau}| + |x\xi_{x\tau\tau}| \\
		& ~~~~ + (|\xi_\tau| + |\xi_{x\tau}|)|\zeta_\tau| \biggr\rbrack \,dx\,d\tau {\nonumber} \lesssim \epsilon  \mathcal{\hat D}_{\zeta,2} + C_\epsilon (1+\hat\omega) a_1^{-1}\sup_\tau e^{(l_2 - r_1 - 5)a_1\tau} \cdot \mathcal{\hat D}_{\xi,11}\\
		& ~~~~ + C_\epsilon (1+\hat\omega) \sup_\tau e^{(l_2-r_1-7)a_1\tau} \cdot\mathcal{\hat D}_{\xi,12}  {\nonumber} + C_\epsilon (1+\hat\omega)a_1^{-1} \sup_\tau e^{(l_2-r_2-2)a_1\tau} \cdot \mathcal{\hat D}_{\xi,21}\\
		& ~~~~ + C_\epsilon (1+\hat\omega) \sup_\tau e^{(l_2-r_2-4)a_1\tau} \cdot \mathcal{\hat D}_{\xi,22} {\nonumber} + C_\epsilon \hat\omega \int e^{-2a_1\tau} \,d\tau \cdot \mathcal{\hat E}_{\zeta,2} , \\
		& \int \hat L_8 \,d\tau \lesssim (1+\hat\omega) \int e^{(l_2-1)a_1\tau} \int x^2 |\zeta_\tau| \biggl\lbrack a_1  (|\xi_\tau|^2 + |x\xi_{x\tau}|^2 ) + ( |\xi_\tau|^3 + |x\xi_{x\tau}|^3 ) {\nonumber} \\
		& ~~~~ +  (|\xi_\tau| + |x\xi_{x\tau}|) \cdot (|\xi_{\tau\tau}| + |x\xi_{x\tau\tau}|) \biggr\rbrack \,dx \,d\tau \lesssim \epsilon \mathcal{\hat D}_{\zeta,2} + C_\epsilon \hat\omega ( a_1 + a_1^{-1} ) \sup_\tau e^{(l_2-r_1-3)a_1\tau} \cdot \mathcal{\hat D}_{\xi, 11} {\nonumber}\\
		& ~~~~ + C_\epsilon \hat\omega (a_1^2 + 1 ) \sup_\tau e^{(l_2-r_1-5)a_1\tau} \cdot \mathcal{\hat D}_{\xi, 12} + C_\epsilon \hat\omega a_1^{-1} \sup_\tau e^{(l_2-r_2)a_1\tau} \cdot \mathcal{\hat D}_{\xi, 21} {\nonumber}\\
		& ~~~~ + C_\epsilon \hat\omega \sup_\tau e^{(l_2-r_2-2)a_1\tau} \cdot \mathcal{\hat D}_{\xi,22}.
	\end{align*}
	Here we have applied similar arguments as in \eqref{ine:102}, \eqref{ine:103} to obtain the above estimates. 
	\todo{13:00}
	Now we shall choose $ r_2, l_2 $ satisfying 
	\begin{equation*}\tag{\ref{constrain=002}}
		r_2 \leq r_1 - 1 < 0 , ~ l_2 + 2 \leq 0, ~ 0 \leq r_2-l_2 < 2.
	\end{equation*}
	Then the above estimates yield
	\begin{align*}
		& \mathcal{\hat E}_{\xi,2} + (1-\hat\omega) \mathcal{\hat D}_{\xi,2} \lesssim \epsilon \mathcal{\hat D}_{\xi,2} + C_\epsilon(1+\hat\omega)(a_1^{-2}+a_1^{-1})\mathcal{\hat E}_{\zeta,2} + C_\epsilon(1+\hat\omega) (a_1+a_1^2) \mathcal{\hat D}_{\zeta,1} {\nonumber}\\
		& ~~~~ + C_\epsilon (a_1^2+ a_1^{-2}  + 1 ) \mathcal{\hat D}_{\xi,1} + (1+\hat\omega) (a_1 +1 ) (\int x^4 \bar\rho\xi^2 \,dx + \int x^4 \xi_x^2 \,dx) +\mathcal{\hat E}_0,\\
		& \mathcal{\hat E}_{\zeta,2} + (1-\hat\omega)\mathcal{\hat D}_{\zeta,2} \lesssim \epsilon \mathcal{\hat D}_{\zeta,2} + C_{\omega} (1+\hat\omega) (1+ a_1^{-1}) \mathcal{\hat D}_{\xi,2}  + C_\epsilon (1+\hat\omega)(1 + a_1^{-1} + a_1^2 ) \mathcal{\hat D}_{\xi,1} {\nonumber}\\
		& ~~~~ + C_\epsilon (\hat\omega + \hat\omega a_1^2 ) \mathcal{\hat D}_{\zeta,1} + \hat\omega a_1^{-1} \mathcal{\hat E}_{\zeta,2} +  \mathcal{\hat E}_0. 
	\end{align*}
	Therefore, after choosing $ \epsilon > 0 $ small enough, for  sufficiently small $ \hat\varepsilon_0 > 0 $ and sufficiently large $ a_1 $, these estimates imply
	\begin{equation*}\tag{\ref{lm:1002}}
		\mathcal{\hat E}_{\xi,2} + \mathcal{\hat E}_{\zeta,2} + \mathcal{\hat D}_{\xi,2}  + \mathcal{\hat D}_{\zeta,2} \lesssim C(\hat\varepsilon_0,a_1,r_1,l_1,r_2,l_2)  \mathcal{\hat E}_0,
	\end{equation*}
	where we have applied \eqref{lm:902} and \eqref{Rene:ub=001}.
\end{pf}
\todo{13:30}


\subsection{Interior Estimates}\label{sec:interior-est-thermodynamic}
In the following, without loss of generality, we assume $ 0 < \hat\omega < \hat \varepsilon_0 < 1 $.  Similarly as before, we shall perform some interior estimates in this section. 
To start with,  multiply \subeqref{eq:RGS-ptb-ln-tm-vb}{1} with $ \chi x \xi_\tau $ and integrate the resulting equation in the spatial variable. After integration by parts, it holds,
\begin{equation}\label{Rene:id=006}
	\dfrac{4\mu}{3}\tilde\alpha^3 \int \chi \biggl\lbrace \dfrac{(\xi_\tau + x\xi_{x\tau})^2}{1+\xi+x\xi_x} + \dfrac{(1+\xi+x\xi_x)\xi_{\tau}^2}{(1+\xi)^2} \biggr\rbrace\,dx = \hat J_1 + \hat J_2 + \hat J_3 + \hat J_4,
\end{equation}
where
\begin{align*}
	& \hat J_1 : = - \dfrac{8\mu}{3} \tilde\alpha^3 \int \chi' \dfrac{x\xi_\tau^2}{1+\xi} \,dx - \tilde\alpha^3 \int \chi' x \mathfrak{\hat B} \xi_\tau\,dx {\nonumber} \\
	& ~~~~~~~~ + \int \chi' K\bar\rho x\xi_\tau \bigl( \dfrac{\zeta+\bar\theta}{(1+\xi)^2(1+\xi+x\xi_x)} - \dfrac{\bar\theta}{(1+\xi)^4} \bigr) \,dx, \\
	& \hat J_2 : = - \tilde\alpha\int \chi \dfrac{x^2 \bar\rho\xi_{\tau\tau}\xi_{\tau}}{(1+\xi)^2} \,dx - \tilde\alpha_\tau \int \chi \dfrac{x^2 \bar\rho\xi_\tau^2}{(1+\xi)^2} \,dx,\\
	& \hat J_3 : = \int \chi K \bar\rho \bigl( \dfrac{(\zeta+\bar\theta)(\xi_\tau+x\xi_{x\tau})}{(1+\xi)^2(1+\xi+x\xi_x)} - \dfrac{\bar\theta(\xi_\tau+x\xi_{x\tau})}{(1+\xi)^4} +\dfrac{4 \bar\theta x \xi_x \xi_\tau}{(1+\xi)^5}\bigr) \,dx.
\end{align*}
Similarly as before, we have the following estimates on $ \hat J_k$'s on the right of \eqref{Rene:id=006}.
\begin{align*}
	& \hat J_1 \lesssim e^{3a_1\tau} \int x^4 \bar\rho \xi_\tau^2 \,dx + e^{3a_1\tau} \int x^2 \lbrack (1+\xi)x\xi_{x\tau} - x\xi_x\xi_\tau\rbrack^2 \,dx\\
	& ~~~~ + \bigl( \int x^4 \bar\rho \xi_\tau^2\,dx \bigr)^{1/2}\biggl\lbrace \bigl( \int x^2 \bar\rho\zeta^2 \,dx\bigr)^{1/2} + \bigl(\int x^4 \bar\rho\xi^2 \,dx \bigr)^{1/2} + \bigl(\int x^4 \xi_x^2 \,dx \bigr)^{1/2}\biggr\rbrace,\\
	& \hat J_2 \lesssim e^{a_1\tau} \biggl\lbrace \bigl(\int x^4\bar\rho \xi_{\tau\tau}^2 \,dx \bigr)^{1/2} + a_1 \bigl( \int x^4 \bar\rho \xi_{\tau}^2 \,dx \bigr)^{1/2} \biggr\rbrace \times \bigl( \int \chi \xi_\tau^2 \,dx \bigr)^{1/2},\\
	& \hat J_3 \lesssim \biggl\lbrace \bigl( \int \chi \zeta^2 \,dx  \bigr)^{1/2} + \bigl( \int \chi ((\xi+x\xi_{x})^2 + \xi^2 )\,dx \bigr)^{1/2} \biggr\rbrace \times \bigl( \int \chi ((\xi_\tau + x\xi_{x\tau})^2 + \xi_\tau^2 ) \,dx \bigr)^{1/2}.
\end{align*}

We will derive from \eqref{Rene:id=006} the following lemma.
\begin{lm}\label{lemma:thermodynamic-interior-1}
	Under the same assumptions as in Lemma \ref{lemma:thermodynamic-high-order-est}, we have the following inequalities,
	
	\begin{gather}
		\int e^{(3+\mathfrak r)a_1\tau} \int \chi ((\xi_\tau+x\xi_{x\tau})^2 + \xi_\tau^2 ) \,dx \,d\tau  \leq C(\hat\varepsilon_0, a_1,r_1,r_2,l_1,l_2,\mathfrak r) \mathcal{\hat E}_0, \label{Rene:ineq=302}\\
		\int\chi(\xi^2 + x^2 \xi_x^2)\,dx
		\leq C(\hat\varepsilon_0, a_1,r_1,r_2,l_1,l_2,\mathfrak r) \mathcal{\hat E}_0 , \label{Rene:ineq=306}
	\end{gather}
	where
	\begin{equation}\label{constrain=003}
		- 3 < \mathfrak r \leq r_2 - 1 < 0.
	\end{equation}

	In addition, the following estimates hold,
	\begin{equation}\label{Rene:ineq=305}
		\begin{gathered}
			\int x^2 \lbrack (1+\xi)x\xi_{x\tau}-x\xi_x\xi_\tau\rbrack^2 \,dx  \leq C(\hat\varepsilon_0,a_1,r_1,l_1,r_2,l_2) e^{\mathfrak d_1 a_1\tau} \cdot \mathcal{\hat E}_{0},\\
			\int x^2 \zeta_x^2 \,dx \leq C(\hat\varepsilon_0,a_1,r_1,l_1,r_2,l_2) e^{\mathfrak d_2 a_1\tau } \cdot \mathcal{\hat E}_0,
		\end{gathered}
	\end{equation}
	with
	\begin{equation}\label{constrain=004}
		\mathfrak d_1 =  
		- (r_1 + r_2)/2 - 3/2 < 0, ~~	\mathfrak d_2 = 
		- (l_1 + l_2)/2 - 1 > 0.
	\end{equation}
	Moreover, as a corollary, we have
	\begin{equation}\label{lm:1101}
		\int \chi( x^2 \xi_{x\tau}^2 + \xi_\tau^2 ) \,dx \lesssim C(\hat\varepsilon_0,a_1,r_1,l_1,r_2,l_2) e^{\mathfrak d_3 a_1\tau} \cdot \mathcal{\hat E}_0,
	\end{equation}
	with
	\begin{equation}\label{lm:1102}
		\mathfrak d_3  = - r_2 - 2 < 0.
	\end{equation}
	
\end{lm}

\begin{pf}
	Now multiply \eqref{Rene:id=006} with $ e^{\mathfrak{r}a_1\tau} $ and apply the Cauchy's inequality to the resulting equation, It holds,
	\begin{equation}\label{Rene:ineq=301}
		\begin{aligned}
			& (1-\epsilon) e^{(3+\mathfrak r)a_1\tau} \int \chi ((\xi_\tau+x\xi_{x\tau})^2 +\xi_\tau^2) \,dx \lesssim C_\epsilon\biggl\lbrace e^{(\mathfrak r - 1)a_1\tau} \int x^4 \bar\rho\xi_{\tau\tau}^2 \,dx\\
			& ~~~~ + a_1^2 e^{(\mathfrak r - 1)a_1\tau} \int x^4 \bar\rho \xi_{\tau}^2 \,dx + e^{(\mathfrak r - 3)a_1\tau}\int \chi ((\xi+x\xi_x)^2 + \xi^2 )\,dx + e^{(\mathfrak r - 3)a_1\tau} \int \chi\zeta^2 \,dx \biggr\rbrace \\
			& ~~~~ + e^{(3 + \mathfrak r )a_1\tau} \int x^4 \bar\rho \xi_\tau^2 \,dx + e^{(3+\mathfrak r)a_1\tau } \int x^2 \lbrack (1+\xi)x\xi_{x\tau}-x\xi_x\xi_\tau\rbrack^2 \,dx \\
			& ~~~~ + e^{\mathfrak r a_1\tau} \bigl( \int x^4 \bar\rho \xi_\tau^2 \,dx \bigr)^{1/2} \biggl\lbrace \bigl( \int x^2 \bar\rho\zeta^2 \,dx\bigr)^{1/2} + \bigl(\int x^4 \bar\rho\xi^2 \,dx \bigr)^{1/2} + \bigl(\int x^4 \xi_x^2 \,dx \bigr)^{1/2}\biggr\rbrace.
		\end{aligned}	
	\end{equation}
	The inequality \eqref{ine:105} can be applied to $ \xi $. Therefore, similarly as in \eqref{ene:le=006}, we have the following inequality
	\begin{equation}\label{Rene:ineq=304}
		\begin{aligned}
			& e^{(3+\mathfrak r)a_1\tau} \int \chi ((\xi_\tau+x\xi_{x\tau})^2 + \xi_\tau^2 ) \,dx \lesssim e^{(\mathfrak r - 3)a_1\tau} \cdot \int e^{-(3 + \mathfrak r)a_1\tau}\,d\tau  \\
			& ~~~~ \times \int e^{(3+\mathfrak r)a_1\tau} \int \chi ((\xi_\tau+x\xi_{x\tau})^2 + \xi_\tau^2 ) \,dx \,d\tau + \hat R_1,
		\end{aligned}	
	\end{equation}
	where
	\begin{align*}
		& \hat R_1 : = \mathcal{\hat E}_0 \cdot e^{(\mathfrak r-3)a_1\tau} + e^{(\mathfrak r - 3)a_1\tau} \bigl(  \int \chi \zeta^2 \,dx + \int x^2 \bar\rho\zeta^2 \,dx \bigr) \\
		& ~~~~ + e^{(\mathfrak r -1)a_1\tau } \int x^4 \bar\rho\xi_{\tau\tau}^2 \,dx + ( a_1^2 e^{(\mathfrak r-1)a_1\tau} + e^{(3+\mathfrak r)a_1\tau} ) \int x^4 \bar\rho\xi_\tau^2 \,dx \\
		& ~~~~ + e^{(3 + \mathfrak r)a_1\tau } \int x^2 \lbrack (1+\xi)x\xi_{x\tau} - x\xi_x\xi_\tau\rbrack^2 \,dx + e^{(\mathfrak r-3)a_1\tau} \cdot \sup_\tau \biggl\lbrace \int x^4 \bar\rho \xi^2 \,dx + \int x^4 \xi_x^2 \,dx \biggr\rbrace\\
		& \lesssim \mathcal{\hat E}_0 \cdot e^{(\mathfrak r-3)a_1\tau} + a_1^{-1} e^{(\mathfrak r - l_1 - 3)a_1\tau} \cdot a_1 e^{l_1a_1\tau} \int x^2 \bar\rho\zeta^2 \,dx + e^{(\mathfrak r-l_1-5)a_1\tau} \cdot e^{(2+l_1)a_1\tau} \int x^2 \zeta_x^2 \,dx  \\
		& ~~~~ + a_1^{-1} e^{(\mathfrak r - r_2  + 1)a_1\tau } \cdot a_1 e^{(r_2-2)a_1\tau} \int x^4 \bar\rho\xi_{\tau\tau}^2 \,dx + ( a_1 e^{(\mathfrak r - r_1 - 2)a_1\tau} + a_1^{-1} e^{(\mathfrak r - r_1+2)a_1\tau} ) \\
		& ~~~~ \times  a_1 e^{(1+r_1)a_1\tau} \int x^4 \bar\rho\xi_\tau^2 \,dx + e^{(\mathfrak r - r_1)a_1\tau } \cdot e^{(3+r_1)a_1\tau} \int x^2 \lbrack (1+\xi)x\xi_{x\tau} - x\xi_x\xi_\tau\rbrack^2 \,dx \\
		& ~~~~ + e^{(\mathfrak r-3)a_1\tau} \cdot \sup_\tau \biggl\lbrace \int x^4 \bar\rho \xi^2 \,dx + \int x^4 \xi_x^2 \,dx \biggr\rbrace.
	\end{align*}
	Then by applying Gr\"onwall's inequality together with \eqref{lm:902}, \eqref{Rene:ub=001} and \eqref{lm:1002}, it admits
	\begin{equation*}\tag{\ref{Rene:ineq=302}}
		\begin{aligned}
			& \int e^{(3+\mathfrak r)a_1\tau} \int \chi ((\xi_\tau+x\xi_{x\tau})^2 + \xi_\tau^2 ) \,dx \,d\tau \lesssim e^{C\int e^{(\mathfrak r-3)a_1\tau} \,d\tau \cdot \int e^{-(3+\mathfrak r)a_1\tau} \,d\tau}\\
			& ~~~~ \times \int \hat R_1 \,d\tau \leq C(\hat\varepsilon_0, a_1,r_1,r_2,l_1,l_2,\mathfrak r) \mathcal{\hat E}_0,
		\end{aligned}
	\end{equation*}
	provided
	\begin{equation*}\tag{\ref{constrain=003}}
		- 3 < \mathfrak r \leq r_2 - 1 < 0.
	\end{equation*}
	Consequently,
	\begin{equation*}\tag{\ref{Rene:ineq=306}}
		\begin{aligned}
			& \int \chi (x^2\xi_x^2 + \xi^2) \,dx \lesssim \int \chi((\xi+x\xi_x)^2 + \xi^2)\,dx \lesssim \mathcal{\hat E}_0 \\
			& ~~~~~~ + \int e^{-(3+\mathfrak r)a_1\tau} \,d\tau \cdot \int e^{(3+\mathfrak r)a_1\tau} \int \chi ((\xi_\tau+x\xi_{x\tau})^2 + \xi_\tau^2 ) \,dx \,d\tau \\
			& ~~~~ \leq C(\hat\varepsilon_0, a_1,r_1,r_2,l_1,l_2,\mathfrak r) \mathcal{\hat E}_0. 
		\end{aligned}
	\end{equation*}
	\todo{16:30}
	
	On the other hand, from \eqref{Rene:id=001} and \eqref{Rene:id=002}, 
	\begin{equation}\label{Rene:ineq=303}
		\begin{aligned}
			& e^{(3+r_1)a_1\tau} \int x^2 \lbrack (1+\xi)x\xi_{x\tau} - x\xi_x \xi_\tau\rbrack^2 \,dx  \lesssim a_1 e^{(1+r_1)a_1\tau} \int x^4 \bar\rho \xi_\tau^2 \,dx \\
			& ~~~~ + e^{(1+r_1)a_1\tau} \bigl(\int x^4 \bar\rho\xi_{\tau\tau}^2 \,dx \bigr)^{1/2} \bigl(\int x^4 \bar\rho\xi_\tau^2 \,dx\bigr)^{1/2} + \hat I_1 + \hat I_2 + \hat I_3,\\
			& e^{(2+l_1)a_1\tau} \int x^2 \zeta_x^2 \,dx \leq a_1 e^{l_1a_1\tau} \int x^2 \bar\rho \zeta^2 \,dx + e^{l_1  a_1\tau} \bigl(\int x^2 \bar\rho\zeta_\tau^2 \,dx \bigr)^{1/2} \bigl( \int x^2 \bar\rho\zeta^2 \,dx \bigr)^{1/2} \\
			& ~~~~  + \hat I_4 + \hat I_5 + \hat I_6. 
	\end{aligned}\end{equation}
	The right of the above inequalities can be estimated as in the following. \todo{17:00}
	\begin{align*}
		& \hat I_1 
		\lesssim e^{r_1a_1\tau} \bigl( \int x^2 \bar\rho \zeta^2\,dx \bigr)^{1/2} \times\bigl( \int x^4 \bar\rho \xi_\tau^2 \,dx + \int x^2 \lbrack (1+\xi)x\xi_{x\tau}-x\xi_x\xi_\tau\rbrack^2 \,dx \bigr)^{1/2} \\
		& ~~~~\lesssim \epsilon e^{(3+r_1)a_1\tau} \int x^2 \lbrack (1+\xi)x\xi_{x\tau} - x\xi_x\xi_\tau\rbrack^2 \,dx + C_\epsilon  e^{(r_1 - l_1 - 3)a_1\tau} \mathcal{\hat E}_{\zeta, 1}\\
		& ~~~~ + e^{(r_1/2 - 1/2 - l_1/2)a_1\tau} \cdot \mathcal{\hat E}_{\xi,1}^{1/2} \mathcal{\hat E}_{\zeta,1}^{1/2}, \\
		& \hat I_2, \hat I_3 
		\lesssim \epsilon e^{(3+r_1)a_1\tau} \int x^2 \lbrack (1+\xi)x\xi_{x\tau}-x\xi_x\xi_\tau\rbrack^2 \,dx + C_\epsilon e^{(r_1-3)a_1\tau} \\
		& ~~~~\times \sup_\tau \biggl\lbrace \int x^4 \bar\rho \xi^2 \,dx + \int x^4\xi_x^2 \,dx \biggr\rbrace + e^{(r_1/2-1/2)a_1\tau} \sup_\tau \biggl\lbrace (\int x^4 \bar\rho\xi^2 \,dx)^{1/2} \\
		&~~~~ + ( \int x^4 \xi_x^2 \,dx)^{1/2} \biggr\rbrace \cdot \mathcal{\hat E}_{\xi,1}^{1/2},\\
		& \hat I_4 
		\lesssim e^{(2+l_1)a_1\tau} \bigl( \int x^4 \bar\rho \xi^2\,dx + \int x^4 \xi_x^2 \,dx \bigr)^{1/2} \times \bigl(\int x^2 \zeta_x^2\,dx \bigr)^{1/2} \lesssim \epsilon e^{(2+l_1) a_1\tau} \int x^2 \zeta_x^2 \,dx\\
		& ~~~~ + C_\epsilon e^{(2+l_1) a_1\tau} \sup_\tau \biggl\lbrace \int x^4 \bar\rho \xi^2 \,dx + \int x^4 \xi_x^2\,dx \biggr\rbrace,\\
		& \hat I_5 
		\lesssim e^{l_1a_1\tau} \bigl( \int x^2 \bar\rho\zeta^2 \,dx \bigr)^{1/2} \bigl( \int x^4 \bar\rho\xi_\tau^2 \,dx + \int x^2 \lbrack (1+\xi)x\xi_{x\tau}-x\xi_x\xi_\tau\rbrack^2 \,dx \bigr)^{1/2}\\
		& ~~~~ \lesssim e^{(l_1/2 - r_1/2 - 1/2 )a_1\tau} \mathcal{\hat E}_{\zeta,1}^{1/2} \mathcal{\hat E}_{\xi,1}^{1/2} + e^{l_1/2a_1\tau} \mathcal{\hat E}_{\zeta,1}^{1/2} \bigl(\int x^2 \lbrack (1+\xi)x\xi_{x\tau}-x\xi_x\xi_\tau\rbrack^2 \,dx\bigr)^{1/2},\\
		& \hat I_6 
		\lesssim e^{(1+l_1)a_1\tau} \hat\omega \bigl( \int x^2 \zeta_x^2 \,dx \bigr)^{1/2} \bigl( \int x^4 \bar\rho \xi_\tau^2\,dx + \int x^2 \lbrack (1+\xi)x\xi_{x\tau} - x\xi_x\xi_\tau\rbrack^2 \,dx \bigr)^{1/2} \\
		& ~~~~ \lesssim \epsilon e^{(2+l_1)a_1\tau} \int x^2 \zeta_x^2 \,dx + C_\epsilon \hat\omega e^{(l_1-r_1-1)a_1\tau} \mathcal{\hat E}_{\xi,1} + C_\epsilon \hat\omega e^{l_1 a_1\tau}\\
		& ~~~~ \times \int x^2 \lbrack (1+\xi)x\xi_{x\tau} - x\xi_x\xi_\tau\rbrack^2 \,dx.
	\end{align*}
	Therefore, from \eqref{Rene:ineq=303}, together with \eqref{lm:902}, \eqref{Rene:ub=001} and \eqref{lm:1002}, we have the following estimates
	\begin{equation*}\tag{\ref{Rene:ineq=305}}
		\begin{gathered}
			\int x^2 \lbrack (1+\xi)x\xi_{x\tau}-x\xi_x\xi_\tau\rbrack^2 \,dx \leq C(\hat\varepsilon_0,a_1,r_1,l_1,r_2,l_2) e^{\mathfrak d_1 a_1\tau} \cdot \mathcal{\hat E}_{0},\\
			\int x^2 \zeta_x^2 \,dx \leq C(\hat\varepsilon_0,a_1,r_1,l_1,r_2,l_2) e^{\mathfrak d_2 a_1\tau } \cdot \mathcal{\hat E}_0,
		\end{gathered}
	\end{equation*}
	where
	\begin{equation*}\tag{\ref{constrain=004}}
		\begin{aligned}
			\mathfrak d_1 :=  &
			- (r_1 + r_2)/2 - 3/2 < 0,\\
			\mathfrak d_2 := & 
			- (l_1 + l_2)/2 - 1 > 0.
		\end{aligned}	
	\end{equation*}
	\todo{17:40}
	Combining \eqref{Rene:ineq=304} and \eqref{Rene:ineq=305} yields the following,
	\begin{align*}
		& \int \chi( ( \xi_\tau + x \xi_{x\tau})^2 + \xi_\tau^2 ) \,dx \lesssim  e^{-6a_1\tau} \cdot\int e^{-(3+\mathfrak r) a_1\tau} \,d\tau\\
		& ~~~~ \times \int e^{(3+\mathfrak r)a_1\tau} \int \chi ((\xi_\tau + x\xi_{x\tau})^2 + \xi_\tau^2) \,dx\,d\tau + \hat R_2,
	\end{align*}
	where
	\begin{align*}
		& \hat R_2 : = e^{-(3+\mathfrak r)a_1\tau} \hat R_1 \lesssim \mathcal{\hat E}_0\cdot e^{-6a_1\tau} + e^{-6a_1\tau} \int x^2 \zeta_x^2 \,dx + \int x^2 \lbrack(1+\xi)x\xi_{x\tau} - x\xi_x \xi_\tau\rbrack^2\,dx \\
		& ~~~~ + e^{(-r_2 - 2)a_1\tau} \mathcal{\hat E}_{\xi,2} + (a_1^2 + 1)e^{(-r_1 - 1)a_1\tau} \mathcal{\hat E}_{\xi,1} + e^{-6a_1\tau} \cdot \sup_\tau \biggl\lbrace \int x^4 \bar\rho \xi^2 \,dx + \int x^4 \xi_x^2 \,dx \biggr\rbrace.
	\end{align*}
	Consequently,
	\begin{equation*}\tag{\ref{lm:1101}}
		\int \chi( x^2 \xi_{x\tau}^2 + \xi_\tau^2 ) \,dx \lesssim e^{\mathfrak d_3 a_1\tau} \cdot \mathcal{\hat E}_0 ,
	\end{equation*}
	where
	\begin{equation*}\tag{\ref{lm:1102}}
		\mathfrak d_3 := \max\lbrace -6, -6 + \mathfrak d_2, \mathfrak d_1, -r_2 - 2, -r_1 - 1 \rbrace = - r_2 - 2 < 0.
	\end{equation*}
\end{pf}

In the following two lemmas, we shall derive some  regularity estimates of $ \xi, \zeta $ around the symmetric center $ x=0 $. 

First, we will need the following lemma concerning the estimate of the interior energy functional. 

\todo{14 Oct 2017}

\begin{lm}\label{lemma:thermodynamic-interior-2}Under the same assumption as in Lemma \ref{lemma:thermodynamic-interior-1}, define the interior energy and dissipation functionals
	\begin{equation}\label{lm:1201}
		\begin{aligned}
			& \mathcal{\hat E}_{\xi,3} : = e^{(r_3-2)a_1\tau} \int \chi x^2 \bar\rho\xi_{\tau\tau}^2 \,dx, \\
			& \mathcal{\hat D}_{\xi,3} : = 
			\int e^{r_3a_1\tau} \int \chi( x^2 \xi_{x\tau\tau}^2 + \xi_{\tau\tau}^2) \,dx \,d\tau,
		\end{aligned}
	\end{equation}
	where $ r_3 $ satisfies 
	\begin{equation}\label{constrain=005}
		r_3 \leq r_2 - 2 < 0.
	\end{equation}
	For $ \hat\omega < \hat \varepsilon_0 $ small enough, we have the following interior energy estimate,
	\begin{equation}\label{lm:1202}
		\mathcal{\hat E}_{\xi,3} + \mathcal{\hat D}_{\xi,3} \leq C(\hat\varepsilon_0, a_1,r_1,r_2,l_1,l_2,\mathfrak r, r_3)\mathcal{\hat E}_0.
	\end{equation} 
	
\end{lm}

\begin{pf}
	
	Multiply \subeqref{eq:1st-RGS-ptb-ln-tm-vb}{1} with $ \tilde\alpha^{r_3} \chi x \xi_{\tau\tau}  $ and integrate the resulting equation in the spatial variable. After integration by parts, it follows, 
	\begin{equation}\label{Rene:id=005}
		\dfrac{d}{d\tau} \hat E_{\xi, 3} + \hat D_{\xi,3} = \hat K_1 + \hat K_2 + \hat K_3	+ \hat K_4 + \hat K_5 + \hat K_6,
	\end{equation}
	where
	\begin{align*}
		& \hat E_{\xi, 3} := \dfrac{\tilde\alpha^{r_3-2}}{2} \int \chi x^2 \bar\rho \xi_{\tau\tau}^2 \,dx, \\
		& \hat D_{\xi,3} : = - \dfrac{r_3}{2} \tilde\alpha^{r_3-3}\tilde\alpha_\tau \int \chi x^2 \bar\rho \xi_{\tau\tau}^2 \,dx {\nonumber} \\
		& ~~~~~~ + \dfrac{4}{3}\mu \tilde\alpha^{r_3} \int \chi \biggl\lbrace \dfrac{(1+\xi)^2x^2 \xi_{x\tau\tau}^2}{1+\xi+x\xi_x} + \dfrac{\lbrack2(1+\xi+x\xi_x)^2 - (1+\xi)^2\rbrack\xi_{\tau\tau}^2}{1+\xi+x\xi_x} \biggr\rbrace \,dx , \\
		& \hat K_1 : = - (\tilde\alpha^{r_3-3}\tilde\alpha_{\tau\tau}-3\tilde\alpha^{r_3-4}\tilde\alpha_\tau^2) \int \chi x^2 \bar\rho \xi_{\tau} \xi_{\tau\tau} \,dx, \\
		& \hat K_2 : = \dfrac{4}{3}\mu\tilde\alpha^{r_3} \int \chi \biggl\lbrace (x(1+\xi)^2\xi_{\tau\tau}) _x\bigl( \dfrac{(\xi_\tau+x\xi_{x\tau})^2}{(1+\xi+x\xi_x)^2} - \dfrac{\xi_\tau^2}{(1+\xi)^2} \bigr) {\nonumber} \\
		& ~~~~~~ + 3(1+\xi)^2\xi_{\tau\tau} \bigl( - \dfrac{2 x\xi_{x\tau} \xi_\tau}{(1+\xi)^2} + \dfrac{2x\xi_x \xi_{\tau}^2}{(1+\xi)^3} \bigr) + 6 x(1+\xi)\xi_\tau\xi_{\tau\tau} \bigl(\dfrac{\xi_\tau}{1+\xi}\bigr)_x {\nonumber}\\
		& ~~~~~~ - 2 (x(1+\xi)\xi_\tau\xi_{\tau\tau})_x \bigl( \dfrac{\xi_\tau + x \xi_{x\tau}}{1+\xi+x\xi_x} - \dfrac{\xi_\tau}{1+\xi} \bigr) \biggr\rbrace \,dx  ,   \\
		& \hat K_3 : = - \tilde\alpha^{r_3} \int \chi' \biggl\lbrace x(1+\xi)^2\xi_{\tau\tau}\mathfrak{\hat B}_\tau + 2 x(1+\xi)\xi_{\tau}\xi_{\tau\tau}\mathfrak{\hat B} + \dfrac{4}{3}\mu x(1+\xi)\xi_{\tau\tau}^2  \biggr\rbrace \,dx, \\
		& \hat K_4 : = - 3 \tilde\alpha^{r_3-4} \tilde\alpha_\tau \int \chi \biggl\lbrace \dfrac{K\bar\rho(\zeta+\bar\theta)\lbrack x(1+\xi)^2\xi_{\tau\tau}\rbrack_x}{(1+\xi)^2(1+\xi+x\xi_x)} - K\bar\rho\bar\theta \bigl( \dfrac{x\xi_{\tau\tau}}{(1+\xi)^2}\bigr)_x \biggr\rbrace \,dx ,\\
		& \hat K_5 : = \tilde\alpha^{r_3-3} \int \chi \biggl\lbrace (\xi_{\tau\tau}+x\xi_{x\tau\tau}) \bigl\lbrack \dfrac{K\bar\rho(\zeta+\bar\theta)}{1+\xi+x\xi_x} - \dfrac{K\bar\rho\bar\theta}{(1+\xi)^2}\bigr\rbrack_\tau {\nonumber}\\
		& ~~~~~~ + \xi_{\tau\tau} \bigl\lbrack \dfrac{2K \bar\rho(\zeta+\bar\theta) x\xi_x}{(1+\xi)(1+\xi+x\xi_x)} + \dfrac{2K\bar\rho\bar\theta x\xi_x}{(1+\xi)^3} \bigr\rbrack_\tau \biggr\rbrace \,dx ,\\
		& \hat K_6 : = - 3 \tilde\alpha^{r_3-4}\tilde\alpha_\tau \int \chi' x \biggl\lbrace \dfrac{K\bar\rho(\zeta+\bar\theta)\xi_{\tau\tau}}{1+\xi+x\xi_x} - \dfrac{K\bar\rho\bar\theta \xi_{\tau\tau}}{(1+\xi)^2} \biggr\rbrace \,dx {\nonumber} \\
		& ~~~~~~ + \tilde\alpha^{r_3-3}\int \chi' x \xi_{\tau\tau} \biggl\lbrace \dfrac{K\bar\rho(\zeta+\bar\theta)}{1+\xi+x\xi_x} - \dfrac{K\bar\rho\bar\theta}{(1+\xi)^2} \biggr\rbrace_\tau\, dx.
	\end{align*}
	Then for $ r_3 < 0 $, integration in the temporal variable of \eqref{Rene:id=005} yields the following,
	\begin{equation}\label{Rene:ineq=201}
		\mathcal{\hat E}_{\xi,3} + \mathcal{\hat D}_{\xi,3} \lesssim \int \hat K_1\,d\tau + \int \hat K_2\,d\tau + \int \hat K_3\,d\tau + \int \hat K_4\,d\tau + \int \hat K_5\,d\tau + \int \hat K_6\,d\tau + \mathcal{\hat E}_0 .
	\end{equation}
	The right of \eqref{Rene:ineq=201} are estimated as in the following.
	\begin{align*}
		& \int \hat K_1 \,d\tau 
		\lesssim \epsilon \mathcal{\hat D}_{\xi,3} + C_\epsilon a_1^{3} \sup_\tau e^{(r_3-r_1 - 5)a_1\tau} \cdot \mathcal{\hat D}_{\xi,11},\\
		& \int \hat K_2 \,d\tau 
		\lesssim \epsilon \mathcal{\hat D}_{\xi, 3} {\nonumber} + C_\epsilon \hat\omega {\int e^{r_3 a_1\tau} \int \chi ( \xi_\tau^2 + x^2 \xi_{x\tau}^2 )\,dx \,d\tau}, \\
		& \int \hat K_3\,d\tau 
		\lesssim (1+\hat\omega) \mathcal{\hat D}_{\xi,21}^{1/2} \biggl\lbrace a_1^{-1} \sup_\tau e^{(r_3-r_2+2)a_1\tau} \cdot \mathcal{\hat D}_{\xi,21}^{1/2} + a_1^{-1/2} \sup_\tau e^{(r_3-r_2+1)a_1\tau} \cdot \mathcal{\hat D}_{\xi,22}^{1/2} {\nonumber}\\
		& ~~~~~~ ~~~~ + \hat\omega a_1^{-1} \sup_\tau e^{(r_3-r_2/2-r_1/2 + 1/2)a_1\tau} \cdot \mathcal{\hat D}_{\xi,11}^{1/2} + \hat\omega a_1^{-1/2} \sup_\tau e^{(r_3-r_2/2-r_1/2 - 1/2)a_1\tau} \cdot \mathcal{\hat D}_{\xi,12}^{1/2} \biggr\rbrace, \\ 
		& \int \hat K_4 \,d\tau \lesssim (1+\hat\omega) \int a_1 e^{(r_3-3)a_1\tau} \int \chi \bar\rho ( |\zeta| + |\xi| + |x\xi_x| ) ( |\xi_{\tau\tau}| + |x\xi_{x\tau\tau}|) \,dx\,d\tau {\nonumber}\\
		& ~~~~~~ \lesssim \epsilon \mathcal{\hat D}_{\xi,3} + C_\epsilon (1+\hat\omega) a_1^2 \sup_\tau e^{(r_3-l_1-8)a_1\tau} \cdot \mathcal{\hat D}_{\zeta,12} {\nonumber} \\ 
		& ~~~~~~ ~~~~ + C_\epsilon (1+\hat\omega) \int a_1^2 e^{(r_3-6)a_1\tau} \,d\tau \cdot \sup_{\tau}  {\int \chi ( \xi^2 + x^2 \xi_x^2 ) \,dx} , \\
		& \int \hat K_5 \,d\tau \lesssim (1+\hat\omega) \int e^{(r_3-3)a_1\tau} \int \chi \bar\rho (|\zeta_\tau| + |\xi_{\tau}| + |\xi_{x\tau}| )(|\xi_{\tau\tau}| + |x\xi_{x\tau\tau}|) \,dx\,d\tau {\nonumber}\\
		& ~~~~~~ \lesssim \epsilon \mathcal{\hat D}_{\xi,3} + C_\epsilon(1+\hat\omega) \sup_\tau e^{(r_3-l_2-6)a_1\tau} \cdot \mathcal{\hat D}_{\zeta,2} {\nonumber}\\
		& ~~~~~~ ~~~~ + C_\epsilon (1+\hat\omega) {\int e^{(r_3-6)a_1\tau}  \int \chi  ( \xi_\tau^2 + x^2 \xi_{x\tau}^2 )\,dx \,d\tau}, \\
		& \int \hat K_6 \,d\tau \lesssim (1+\hat\omega) \int  e^{(r_3-3)a_1\tau} \int x^4 \bar\rho (a_1|\zeta| + a_1|\xi| + a_1|x\xi_{x}| + |\zeta_\tau| + |\xi_\tau| + |x\xi_{x\tau}| )|\xi_{\tau\tau}| \,dx\,d\tau {\nonumber} \\
		& ~~~~~~ \lesssim (1+\hat\omega) \mathcal{\hat D}_{\xi, 21}^{1/2}\biggl\lbrace a_1^{1/2} \sup_\tau e^{(r_3-r_2/2-l_1/2 - 3)a_1\tau} \cdot \mathcal{\hat D}_{\zeta,12}^{1/2} + a_1^{1/2} \bigl(\int e^{(2r_3-r_2-4)a_1\tau} \,d\tau\bigr)^{1/2} {\nonumber} \\
		& ~~~~~ ~~~~ \times \sup_\tau \bigl( \int x^4\bar\rho\xi^2 \,dx + \int x^4 \xi_x^2 \,dx\bigr)^{1/2} + a_1^{-1/2} \sup_\tau e^{(r_3-r_2/2-l_2/2 - 2)a_1\tau} \cdot \mathcal{\hat D}_{\zeta,2}^{1/2} {\nonumber}\\
		& ~~~~~ ~~~~ + a_1^{-1} \sup_\tau e^{(r_3-r_2/2-r_1/2 - 5/2)a_1\tau} \cdot \mathcal{\hat D}_{\xi,11}^{1/2} + a_1^{-1/2} \sup_\tau e^{(r_3-r_2/2-r_1/2-7/2)a_1\tau} \cdot \mathcal{\hat D}_{\xi,12}^{1/2} \biggr\rbrace.
	\end{align*}
	Therefore from \eqref{Rene:ineq=201}, it holds the following,
	\begin{align*}
		& \mathcal{\hat E}_{\xi,3} + (1-\hat\omega) \mathcal{\hat D}_{\xi,3} \lesssim C(\hat\varepsilon_0, a_1) \biggl\lbrace \mathcal{\hat D}_{\xi,1} + \mathcal{\hat D}_{\xi,2} + \mathcal{\hat D}_{\zeta,1} + \mathcal{\hat D}_{\zeta,2}  \\
		& ~~~~ + \int e^{(r_3- \mathfrak r - 3)a_1\tau} \cdot e^{(3+\mathfrak r)a_1\tau} \int\chi(\xi_\tau^2 + x^2 \xi_{x\tau}^2 )\,dx \,d\tau + \int e^{(r_3-6)a_1\tau} \,d\tau \cdot \int \chi(\xi^2 + x^2 \xi_x^2 )\,dx \\
		& ~~~~ + \int e^{(r_3 - \mathfrak r - 9)a_1\tau}\cdot e^{(3+\mathfrak r)a_1\tau} \int \chi (\xi_\tau^2 + x^2 \xi_{x\tau}^2 )\,dx\,d\tau \biggr\rbrace \\
		& ~~~~ + \int e^{(2r_3-r_2-4)a_1\tau}\,d\tau \cdot\sup_\tau \bigl(\int x^4 \bar\rho \xi^2\,dx + \int x^4 \xi_x^2\,dx \bigr) + \mathcal{\hat E}_0 \lesssim \mathcal{\hat E}_0,
	\end{align*}
	provided
	\begin{equation*}\tag{\ref{constrain=005}}
		r_3 \leq r_2 - 2 < 0.
	\end{equation*}
	Then for $ \hat \varepsilon_0 $ small enough, we finish the proof of the lemma.
\end{pf}

Now we can derive the following lemma concerning the elliptic estimates of \eqref{eq:RGS-ptb-ln-tm-vb}. Similarly as in the isentropic case, we rewrite the system in the following in terms of the relative entropy functional, $\mathcal{\hat G}  = \mathfrak H(\xi) = \log(1+\xi)^2(1+\xi+x\xi_x) $. That is,
\todo{12:12}
\begin{equation}\label{eq:RGS-alternative}
	\begin{cases}
		\dfrac{4}{3}\mu \mathcal{\hat G}_{x\tau} + \tilde\alpha^{-3} \dfrac{K\bar\rho(\zeta+\bar\theta)}{(1+\xi)^2(1+\xi+x\xi_x)} \mathcal{\hat G}_{x} = \tilde\alpha^{-3} \dfrac{(K\bar\rho \zeta )_x}{(1+\xi)^2(1+\xi+x\xi_x)} \\
		~~~~~~ + \tilde\alpha^{-3} (K\bar\rho\bar\theta)_x \bigl( \dfrac{1}{(1+\xi)^2(1+\xi+x\xi_x)} - \dfrac{1}{(1+\xi)^4} \bigr) + \dfrac{1}{(1+\xi)^2} \bigl( \tilde\alpha^{-2}x\bar\rho\xi_{\tau\tau} \\
		~~~~~~ + \tilde\alpha^{-3}\tilde\alpha_\tau x\bar\rho\xi_\tau \bigr), \\
		\bigl\lbrack \dfrac{x^2(1+\xi)^2}{1+\xi+x\xi_x}\zeta_x\bigr\rbrack_x = \bigl\lbrack \bigl( 1 -\dfrac{(1+\xi)^2}{1+\xi+x\xi_x}\bigr)x^2\bar\theta_x \bigr\rbrack_x + K \tilde\alpha^{-2}  \dfrac{\bar\rho(\zeta+\bar\theta)(x^3(1+\xi)^2\xi_\tau)_x}{(1+\xi)^2(1+\xi+x\xi_x)} \\
		~~~~~~ + 3K \tilde\alpha^{-2} x^2 \bar\rho\zeta_\tau - \tilde\alpha^{-1} x^2 (1 + \xi)^2 (1 + \xi + x\xi_x)\cdot\mathfrak{\hat F}(\xi).
	\end{cases}
\end{equation}

\begin{lm}\label{lemma:thermodynamic-interior-3}Under the same assumption as in Lemma \ref{lemma:thermodynamic-interior-2}, supposed
	\begin{equation}\label{constrain=006} l_2 + 2 \geq 0, \end{equation}
	we have
	\begin{equation}{\label{Rene:L^2est-3}}
		\begin{gathered}
			\norm{\xi_x}{\stnorm{\infty}{2}}^2 + \norm{x\xi_{xx}}{\stnorm{\infty}{2}}^2 + \norm{\zeta_x}{\stnorm{\infty}{2}}^2 \leq C(\hat\varepsilon_0,a_1,r_1,l_1,r_2,l_2,\mathfrak r,r_3) \mathcal{\hat E}_0,  \\
			\norm{\xi_{x\tau}}{\stnorm{\infty}{2}}^2 + \norm{x\xi_{xx\tau}}{\stnorm{\infty}{2}}^2 \lesssim C(\hat\varepsilon_0,a_1,r_1,l_1,r_2,l_2,\mathfrak r,r_3) e^{(-r_3-2)a_1\tau} \mathcal{\hat E}_0.
		\end{gathered}
	\end{equation}
	In the meantime, we have the following
	\begin{equation}\label{lm:1301}
		\norm{x\zeta_{xx}}{\stnorm{\infty}{2}}^2 \leq C(\hat\varepsilon_0,a_1,r_1,l_1,r_2,l_2,\mathfrak r,r_3) (\mathcal{\hat E}_0 + \mathcal{\hat E}_0^2) .
	\end{equation}
\end{lm}

\begin{pf}
	From \subeqref{eq:RGS-alternative}{2}, by noticing \begin{equation*}
		\bigl(1 - \dfrac{(1+\xi)^2}{1+\xi+x\xi_x} \bigr)_x  = - \dfrac{4 (1+\xi)^3 \xi_x}{(1+\xi)^2(1+\xi+x\xi_x)} + \dfrac{(1+\xi)^4}{(1+\xi)^2(1+\xi+x\xi_x)} \mathcal{\hat G}_{x},
	\end{equation*}
	it holds, 
	\begin{align*}
		& \int \dfrac{1}{x^2(1+\xi)^2(1+\xi+x\xi_x)} \biggl\lbrace \bigl\lbrack \dfrac{x^2(1+\xi)^2}{1+\xi+x\xi_x} \zeta_x \bigr\rbrack_x \biggr\rbrace^2 \,dx \lesssim (1+\hat\omega) \int x^2 (\xi^2 + x^2 \xi_x^2 )\,dx \\
		& ~~~~ + (1+\hat\omega)  \int \mathcal{\hat G}_{x}^2 \,dx + (1+\hat\omega)  e^{-4a_1\tau} \int x^2 ( \xi_\tau^2 + x^2 \xi_{x\tau}^2 )\,dx \\
		& ~~~~ + (1+\hat\omega)  e^{-4a_1\tau} \int x^2 \bar\rho\zeta_\tau^2 \,dx + \hat\omega e^{-2a_1\tau} \int x^2 (\xi_\tau^2 + x^2 \xi_{x\tau}^2)\\
		& \lesssim \int \mathcal{\hat G}_{x}^2 \,dx + ( e^{(\mathfrak d_1 - 4) a_1\tau} + \hat\omega e^{(\mathfrak d_1 - 2) a_1\tau}) \cdot \mathcal{\hat E}_0 + ( e^{(- r_1 - 5)b\tau} + \hat\omega e^{(-r_1 - 3 )b\tau}) \cdot \mathcal{\hat E}_{\xi,1} \\
		& ~~~~ + e^{(-l_2-2)a_1\tau} \cdot\mathcal{\hat E}_{\zeta,2} + \mathcal{\hat E}_0 \lesssim \int \mathcal{\hat G}_x^2 \,dx + \mathcal{\hat E}_0,
	\end{align*}
	provided \eqref{constrain=006} holds, 
	where we have used \eqref{lm:902}  \eqref{Rene:ub=001}, \eqref{Rene:ineq=005},  \eqref{lm:1002}, \eqref{Rene:ineq=305} and the fact $ \bar\theta_x(0) = 0 $ which implies $ \abs{(x^2 \bar\theta_x)_x}{} \lesssim x^2 \norm{\bar\theta_{xx}}{L_x^\infty}, \abs{\bar\theta_x}{} \lesssim x \norm{\bar\theta_{xx}}{L_x^\infty} $.
	Meanwhile, direct calculation shows,
	\begin{align*}
		& \int \dfrac{1}{x^2(1+\xi)^2(1+\xi+x\xi_x)} \biggl\lbrace \bigl\lbrack \dfrac{x^2(1+\xi)^2}{1+\xi+x\xi_x} \zeta_x \bigr\rbrack_x \biggr\rbrace^2 \,dx = \int \biggl\lbrack \dfrac{x (1+\xi)}{\sqrt{1+\xi+x\xi_x}} \bigl(\dfrac{\zeta_x}{1+\xi+x\xi_x}\bigr)_x \\
		& ~~~~ +  2 \sqrt{1+\xi+x\xi_x} \dfrac{\zeta_x}{1+\xi+x\xi_x}  \biggr\rbrack^2 \,dx = \int  \dfrac{x^2(1+\xi)^2}{1+\xi+x\xi_x} \biggl\lbrack \bigl(\dfrac{\zeta_x}{1+\xi+x\xi_x}\bigr)_x\biggr\rbrack^2  \\
		& ~~~~ + 2  (1+\xi+x\xi_x) \bigl\lbrack \dfrac{\zeta_x}{1+\xi+x\xi_x} \bigr\rbrack^2 \,dx + 2 x(1+\xi) \bigl\lbrack \dfrac{\zeta_x}{1+\xi+x\xi_x}\bigr\rbrack^2 \biggr|_{x=R_0}.
	\end{align*}
	Therefore, combining the above inequalities yields,
	\begin{gather}
		\int \zeta_x^2 \,dx \lesssim \int \mathcal{\hat G}_x^2 \,dx + \mathcal{\hat E}_0, {\label{Rene:L^2est-1}} \\
		\int x^2 \zeta_{xx}^2 \,dx \lesssim \int \zeta_x^2 \,dx + (\norm{x\zeta_x}{L_x^\infty}^2 + 1) \int \mathcal{\hat G}_x^2 \,dx  + \mathcal{\hat E}_0. {\label{Rene:L^2est-2}}
	\end{gather}
	Here we have applied the following identity
	\begin{equation*}
		\bigl( \dfrac{\zeta_x}{1+\xi+x\xi_x} \bigr)_x 
		= \dfrac{2(1+\xi)\xi_x \zeta_x + (1+\xi)^2 \zeta_{xx}}{(1+\xi)^2(1+\xi+x\xi_x)} - \dfrac{(1+\xi)^2\zeta_x}{(1+\xi)^2(1+\xi+x\xi_x)} \mathcal{\hat G}_x.
	\end{equation*}
	
	On the other hand, after multiplying \subeqref{eq:RGS-alternative}{1} with $ \mathcal{\hat G}_x $ and integrating the resultant, one can derive, similarly as in the isentropic case, 
	\begin{align*}
		& \dfrac{d}{d\tau} \dfrac{2\mu}{3} \int \mathcal{\hat G}_x^2 \,dx + \tilde\alpha^{-3} \int \dfrac{K\bar\rho(\zeta+\bar\theta)}{(1+\xi)^2(1+\xi+x\xi_x)} \mathcal{\hat G}_x^2  \,dx \lesssim \epsilon e^{-3a_1\tau} \int \mathcal{\hat G}_x^2 \,dx + C_\epsilon \biggl\lbrace e^{-3a_1\tau}\int \zeta_{x}^2\,dx \\
		& ~~~~~~ + e^{-3a_1\tau} \int (\xi^2 + x^2\xi_x^2) \,dx + a_1 e^{(-2-r_1)a_1\tau} \cdot a_1 e^{(1+r_1)a_1\tau}\int x^4\bar\rho \xi_\tau^2 \,dx \\
		& ~~~~~~ + a_1^2 e^{(-4-r_1)a_1\tau} \cdot e^{(3+r_1)a_1\tau}\int x^2( (1+\xi)x\xi_{x\tau}-x\xi_x\xi_\tau)^2\,dx  \biggr\rbrace \\
		& ~~~~ +  (a_1^{-1}e^{(-r_2-2)a_1\tau} + e^{(-r_2-4)a_1\tau})  \cdot \int \mathcal{\hat G}_x^2\,dx  +  a_1 e^{(r_2-2)a_1\tau} \int x^4 \bar\rho \xi_{\tau\tau}^2\,dx\\
		& ~~~~ +  e^{r_2a_1\tau} \int x^2 ((1+\xi)x\xi_{x\tau\tau}-x\xi_x\xi_{\tau\tau})^2\,dx.
	\end{align*}
	Then by noticing \eqref{lm:902}, \eqref{Rene:ub=001}, \eqref{lm:1002}, \eqref{Rene:ineq=306} \eqref{Rene:L^2est-1}, applying the Gr\"onwall's inequality then yields,
	\begin{equation}\label{Rene:L^2est-03}
		\int \mathcal{\hat G}_x^2 \,dx + \int \zeta_x^2 \,dx  \lesssim \mathcal{\hat E}_0.
	\end{equation}
	Meanwhile, from \subeqref{eq:RGS-alternative}{1}, \todo{14:45}
	\begin{align*}
		& \int \mathcal{\hat G}_{x\tau}^2 \,dx \lesssim e^{-6a_1\tau} \int \mathcal{\hat G}_x^2 \,dx + e^{-6a_1\tau} \int \zeta_x^2 \,dx + e^{-6a_1\tau} \int (\xi^2 + x^2 \xi_x^2) \,dx + e^{-4a_1\tau} \int x^2 \bar\rho \xi_{\tau\tau}^2\,dx\\
		& ~~~~ + a_1^2 e^{-4a_1\tau} \int x^2 \bar\rho \xi_\tau^2 \,dx \lesssim  (e^{(-r_3-2)a_1\tau} + e^{(-r_2-2)b\tau} + a_1^2 e^{(-r_1-5)a_1\tau} + a_1^2 e^{(\mathfrak d_3 -4)a_1\tau} + e^{-6a_1\tau}) \cdot \mathcal{\hat E}_0,
	\end{align*}
	where \eqref{lm:902}, \eqref{Rene:ub=001} \eqref{lm:1002}, \eqref{Rene:ineq=306}, \eqref{lm:1101}, \eqref{lm:1202}, \eqref{Rene:L^2est-03} are applied. Consequently, 
	\begin{equation}\label{Rene:L^2est-04}
		\int \mathcal{\hat G}_{x\tau}^2 \,dx \lesssim e^{(-r_3-2)a_1\tau} \mathcal{\hat E}_0. 
	\end{equation}
	Then together with Lemma \ref{lm:estimates-of-relative-entropy}, \eqref{Rene:L^2est-03} and \eqref{Rene:L^2est-04} imply
	\begin{equation*}{\tag{\ref{Rene:L^2est-3}}}
		\begin{gathered}
			\norm{\xi_x}{\stnorm{\infty}{2}}^2 + \norm{x\xi_{xx}}{\stnorm{\infty}{2}}^2 + \norm{\zeta_x}{\stnorm{\infty}{2}}^2 \lesssim \mathcal{\hat E}_0,  \\
			\norm{\xi_{x\tau}}{\stnorm{\infty}{2}}^2 + \norm{x\xi_{xx\tau}}{\stnorm{\infty}{2}}^2 \lesssim e^{(-r_3-2)a_1\tau} \mathcal{\hat E}_0. 
		\end{gathered}
	\end{equation*}
	\todo{15:26}
	Consequently, from \eqref{Rene:L^2est-2}
	\begin{equation*}
		\norm{x\zeta_{xx}}{\stnorm{\infty}{2}}^2 \lesssim (1+\norm{x\zeta_x}{L_x^\infty}^2) \mathcal{\hat E}_0 \lesssim \mathcal{\hat E}_0 + \mathcal{\hat E}_0^2 + \norm{\zeta_x}{\stnorm{\infty}{2}}^2 + \dfrac{1}{2} \norm{x\zeta_{xx}}{\stnorm{\infty}{2}}^2,
	\end{equation*}
	where the following pointwise estimate is used
	\begin{equation*}
		\abs{x\zeta_x}{2} \lesssim \norm{x\zeta_x}{L_x^2}^2 + \int \abs{x\zeta_x(\zeta_x + x\zeta_{xx})}{}\,dx  \lesssim (1+1/\epsilon) \norm{\zeta_x}{L_x^2}^2 + \epsilon \norm{x\zeta_{xx}}{L_x^2}^2.
	\end{equation*}
	This finishes the proof.
\end{pf}


\subsection{Pointwise Bounds and Asymptotic Stability}\label{sec:pointwise-bounde-thermodynamic}

In this section, we shall derive the following
\begin{equation}\label{thermodynamic-pointwise-est-001}
	\hat \omega^2 \leq C \mathcal{\hat E}(\tau) + C\mathcal{\hat E}_0 + C \mathcal{\hat E}^2(\tau) + C\mathcal{\hat E}^2_0 \leq C ( \mathcal{\hat E}_0 + \mathcal{\hat E}_0^2).
\end{equation}
This will be sufficient to establish the asymptotic stability theory of the linearly expanding homogeneous solution of the thermodynamic model for the radiation gaseous star. Indeed, applying the continuity arguments as in Section \ref{sec:pointwise-bounde-isentropic} with small enough $ \hat \varepsilon_0 $ and initial energy $ \mathcal{\hat E}_0 $ will finish the proof of Theorem \ref{theorem3}.

What is left is to show the inequality \eqref{thermodynamic-pointwise-est-001}. First of all, similar calculation as in \eqref{lm:602}, \eqref{lm:603}, \eqref{lm:604}, \eqref{lm:605} yields, together with \eqref{Rene:L^2est-3},
\begin{equation}\label{thermodynamic-pointwise-est-002}
	\begin{aligned}
		& \norm{x\xi_x}{\stnorm{\infty}{\infty}}^2 + \norm{\xi}{\stnorm{\infty}{\infty}}^2 \lesssim C \mathcal{\hat E}(\tau) + C\mathcal{\hat E}_0,\\
		& \norm{x\xi_{x\tau}}{\stnorm{\infty}{\infty}}^2 + \norm{\xi_\tau}{\stnorm{\infty}{\infty}}^2 \lesssim (e^{-(1+r_1)a_1\tau} + e^{-((r_1+r_2)/2+3/2)a_1\tau})\mathcal{\hat E}(\tau)   \\
		& ~~~~~~ + ( e^{-(r_2/2+1)a_1\tau}+ e^{-((r_1+r_2)/4+3/4)a_1\tau} + e^{-(r_1/2+1/2)} ) e^{-(r_3/2+1)} \mathcal{\hat E}(\tau)^{1/2} \cdot \mathcal{\hat E}_0^{1/2}  \\
		& ~~~~ \lesssim e^{-((r_2+r_3)/2+2)a_1\tau} (\mathcal{\hat E}(\tau) + \mathcal{\hat E}_0).
	\end{aligned}
\end{equation}
Moreover, we have the following lemma.
\begin{lm}\label{lm:theromodynamic-pointwise-bound}
	Under the same assumption as in Lemma \ref{lemma:thermodynamic-interior-3}, 
	we have the following 
	\begin{equation}\label{lm:1401}
		\begin{aligned}
			& \norm{x\xi_x}{\stnorm{\infty}{\infty}}^2 + \norm{\xi}{\stnorm{\infty}{\infty}}^2 + \norm{x\xi_{x\tau}}{\stnorm{\infty}{\infty}}^2 + \norm{\xi_\tau}{\stnorm{\infty}{\infty}}^2 + \norm{\zeta/\sigma}{\stnorm{\infty}{\infty}}^2 \\
			& ~~~~~~ \leq C(\hat\varepsilon_0,a_1,r_1,l_1,r_2,l_2,\mathfrak r,r_3) ( \mathcal{\hat E}(\tau) + \mathcal{\hat E}_0 + \mathcal{\hat E}^2(\tau) + \mathcal{\hat E}^2_0).
		\end{aligned}
	\end{equation}
	This finishes the first inequality in \eqref{thermodynamic-pointwise-est-001}.
\end{lm}

\begin{pf}
	The pointwise bound of $ x\xi_x, \xi, x\xi_{x\tau},\xi_\tau $ is a direct consequence of  \eqref{thermodynamic-pointwise-est-002}. What is left is the estimate of $ \norm{\zeta/\sigma}{\stnorm{\infty}{\infty}} $. This is a direct consequence of \eqref{Rene:L^2est-3} and \eqref{lm:1301} by applying the poincar\'e inequality, the embedding theory and the fundamental theorem of calculus. 
\end{pf}

The second part of the inequality \eqref{thermodynamic-pointwise-est-001} is a direct consequence of Lemma \ref{lemma:thermodynamic-basic-energy-est}, Lemma \ref{lemma:thermodynamic-high-order-est}, Lemma \ref{lemma:thermodynamic-interior-1}, Lemma \ref{lemma:thermodynamic-interior-2} and Lemma \ref{lemma:thermodynamic-interior-3}.
\begin{lm}\label{lemma:thermodynamic-total-energy-estimate}
	Under the same assumption as in Lemma \ref{lemma:thermodynamic-interior-3}, we have  the total energy inequality $ \mathcal{\hat E}(\tau) + \mathcal{\hat D}(\tau)\leq C(\hat\varepsilon_0,a_1,r_1,l_1,r_2,l_2,\mathfrak r,r_3) (\mathcal{\hat E}_0 + \mathcal{\hat E}_0^2) $. 
\end{lm}


\end{document}